\documentclass[12pt]{article}
\usepackage{epsfig}
\begin{document}
\title{Shape Analysis, Lebesgue Integration and Absolute Continuity
Connections}
\author{Javier Bernal\\
{\small \sl National Institute of Standards and Technology,} \\
{\small \sl Gaithersburg, MD 20899, USA}}
\date{\ }
\maketitle
\begin{abstract}
As shape analysis of the form presented in Srivastava and Klassen's textbook
``Functional and Shape Data Analysis" is intricately related to
Lebesgue integration and absolute continuity, it is advantageous to
have a good grasp of the latter two notions. Accordingly, in these notes
we review basic concepts and results about Lebesgue integration and
absolute continuity. In particular, we review fundamental results
connecting them to each other and to the kind of shape analysis, or more
generally, functional data analysis presented in the aforemetioned textbook,
in the process shedding light on important aspects of all three notions.
Many well-known results, especially most results about Lebesgue integration
and some results about absolute continuity, are presented without proofs.
However, a good number of results about absolute continuity and most
results about functional data and shape analysis are presented with proofs.
Actually, most missing proofs can be found in Royden's ``Real Analysis"
and Rudin's ``Principles of Mathematical Analysis" as it is on these
classic textbooks and Srivastava and Klassen's textbook
that a good portion of these notes are based. However,
if the proof of a result does not appear in the aforementioned textbooks,
nor in some other known publication, or if all by itself it could be of
value to the reader, an effort has been made to present it accordingly.
\end{abstract}
\section{\large Introduction}
The concepts of Lebesgue integration and absolute continuity play a major
role in the theory of shape analysis or more generally in the theory of
functional data analysis of the form presented in~\cite{srivastava}. In
fact, well-known connections between Lebesgue integration and absolute
continuity are of great importance in the development of functional data
and  shape analysis of the kind in~\cite{srivastava}. Accordingly,
understanding functional data and  shape analysis as presented
in~\cite{srivastava} requires understanding the basics of Lebesgue
integration and absolute continuity, and the connections between them.
It is the purpose of these notes to provide a way to do exactly~that.
\par
In Section~2, we review fundamental concepts and results about Lebesgue
integration. Then, in Section~3, we review fundamental concepts and results
about absolute continuity, some results connecting it to
Lebesgue integration. Finally, in Section~4, we shed light on some important
aspects of functional data and shape analysis of the type in Srivastava and
Klassen's textbook \cite{srivastava}, in the process illustrating its
dependence on Lesbesgue integration, absolute continuity and the connections
between them. Accordingly, without page numbers, a table of contents for
these notes would be roughly as follows:
\\
1. Introduction\\
2. Lebesgue Integration\\
Algebras of sets, Borel sets, Cantor set\\
Outer measure\\
Measurable sets, Lebesgue measure\\
Measurable functions, Step functions, Simple functions\\
The Riemann integral\\
The Lebesgue integral\\
The $L^p$ Spaces\\
3. Absolute Continuity and its Connections to Lebesgue Integration\\
4. Functional Data and Shape Analysis and its Connections to Lebesgue Integration
and Absolute Continuity\\
Summary\\
Acknowledgements\\
References\\
Index of Terms
\par
The material in these notes about Lebesgue integration and absolute continuity is
mostly based on Royden's ``Real Analysis"~\cite{royden} and Rudin's
``Principles of Mathematical Analysis"~\cite{rudin}. The fundamental
ideas on functional data and shape analysis are mostly from
Srivastava and Klassen's ``Functional and Shape Data Analysis"~\cite{srivastava}.
An index of terms has been included at the end of the notes.
\section{\large Lebesgue Integration}
{\bf \large Algebras of sets, Borel sets, Cantor set}\\ \smallskip\\
{\bf Definition 2.1:} A collection $\cal A$ of subsets of a set $X$ is called an
{\bf algebra} on~$X$ if for $A$, $B$ in $\cal A$, $A\cup B$
is in $\cal A$, and for $A$ in $\cal A$, $\tilde{A}=X\setminus A$ is
in~$\cal A$.\\ \smallskip \\
{\bf Observation 2.1:} From De Morgan's laws if $\cal A$ is an algebra, then
for $A$, $B$ in $\cal A$, $A\cap B$ is in~$\cal A$.
\\ \smallskip \\
{\bf Definition 2.2:} An algebra $\cal A$ is called a {\bf {\boldmath $\sigma$}-algebra}
if the union of every countable collection of sets in $\cal A$ is
in~$\cal A$.\\ \smallskip \\
{\bf Observation 2.2:} From De Morgan's laws if $\cal A$ is a $\sigma$-algebra, then
the intersection of a countable collection of sets in $\cal A$ is in~$\cal A$.
\\ \smallskip \\
{\bf Definition 2.3:} A set of real numbers $O$ is said to be {\bf open} if for each
$x\in O$ there is $\delta>0$ such that each number $y$ with $|x-y|<\delta$
belongs to~$O$. A set of real numbers $F$ is said to be {\bf closed} if its complement
in {\bf R} is open, i.e., ${\bf R}\setminus F$ is open, where {\bf R} is the set of real
numbers. The collection of {\bf Borel} sets is the smallest $\sigma$-algebra
on the set {\bf R} of real numbers which contains all open sets of real numbers.
\\ \smallskip \\
{\bf Observation 2.3:} The collection of Borel sets contains in particular all closed
sets, all open intervals, all countable unions of closed sets, all countable
intersections of open sets, etc.
\\ \smallskip \\
{\bf Proposition 2.1:} Every open set of real numbers is the union of a countable
collection of disjoint open intervals. Proof in~\cite{royden}.\\ \smallskip \\
{\bf Proposition 2.2 (Lindel\"{o}f):} Given a collection $\cal C$ of open sets of real
numbers, then there is a countable subcollection $\{ O_i\}$ of $\cal C$ with
\[ \cup_{O\in \cal C}\, O = \cup_{i=1}^{\infty}\, O_i. \]
Proof in~\cite{royden}.\\ \smallskip \\
{\bf Definition 2.4:} A set of real numbers $F$ is said to be {\bf compact} if every open
cover of $F$ contains a finite subcover, i.e., if $\cal C$ is a collection of open
sets of real numbers such that $F\subseteq  \cup_{O\in \cal C}\, O$, then
there is a finite subcollection $\{O_i,\ i=1,\ldots,n\}$ of $\cal C$ with
$F\subseteq \cup_{i=1}^n\, O_i$.\\ \smallskip\\
{\bf Proposition 2.3 (Heine-Borel):} A set of real numbers $F$ is compact if and only if it
is closed and bounded. Proof in~\cite{royden} and \cite{rudin}.\\ \smallskip \\
{\bf Proposition 2.4:} Given a collection ${\cal K}$ of closed sets of real numbers such that
at least one of the sets is bounded and the intersection of every finite subcollection
of ${\cal K}$ is nonempty, then $\cap_{F\in \cal K}\, F\not = \emptyset$.
Proof in~\cite{royden} and~\cite{rudin}.\\ \smallskip \\
{\bf Definition 2.5:} A number $x$ is said to be a {\bf limit point} of a set of
real numbers $E$ if every open set that contains $x$ contains $y\not= x$, $y$ in~$E$.
\\ \smallskip\\
{\bf Definition 2.6:} A set of real numbers $E$ is said to be {\bf perfect} if it is
closed and if every number in $E$ is a limit point of $E$.
\\ \smallskip\\
{\bf Proposition 2.5:} A set is closed if and only if every limit point of the set is
a point of the set. A nonempty perfect set is uncountable. Proofs in~\cite{rudin}.
\\ \smallskip\\
{\bf Corollary 2.1:} Every interval is uncountable, thus the set of real numbers
is uncountable.
\\ \smallskip\\
{\bf Observation 2.4:}
Let $E_1$ be the union of the intervals $[0,\frac{1}{3}]$, $[\frac{2}{3},1]$
that are obtained by removing the open middle third of the interval $[0,1]$.
Let $E_2$ be the union of the intervals $[0,\frac{1}{9}]$, $[\frac{2}{9},\frac{3}{9}]$,
$[\frac{6}{9},\frac{7}{9}]$, $[\frac{8}{9},1]$ that are obtained by removing the
open middle thirds of the intervals $[0,\frac{1}{3}]$ and $[\frac{2}{3},1]$.
Continuing this way, a sequence of compact sets $E_n$ is obtained with $E_n\supset E_{n+1}$
for every positive integer~$n$. The set $\cap_{n=1}^{\infty}\, E_n$, called the
{\bf Cantor set}, is compact, nonempty, perfect thus uncountable,
and contains no interval.  Proofs in~\cite{rudin}. \\ \smallskip\\
{\bf Definition 2.7:} The {\bf extended real numbers} consist of the real numbers
together with the two symbols $-\infty$ and~$+\infty$. The definition of $<$ is
extended by declaring that if $x$ is a real number, then $-\infty < x < \infty$.
The operation $\infty -\infty$ is left undefined, the operation $0\cdot(\pm\infty)$
is defined to be $0$, while other definitions are extended: If $x$ is a real number, then\\
$x+\infty = \infty,\ \ \ \ $ $x-\infty = -\infty,\ \ \ \ $ $x/+\infty = x/-\infty = 0$,\\
$x\cdot\infty = \infty,\ \ \ \ $ $x\cdot -\infty =  -\infty\ \ \ \ $ if $x>0$,\\
$x\cdot\infty = -\infty,\ \ \ \ $ $x\cdot -\infty =  \infty\ \ \ \ $ if $x<0$.\\
Finally\\
$\infty + \infty = \infty,\ \ \ \ $ $-\infty -\infty = -\infty,\ \ \ \ $
$\infty\cdot (\pm\infty) = \pm\infty,\ \ \ \ $
$-\infty\cdot (\pm\infty) = \mp\infty.$
\\ \medskip\\
{\bf \large Outer measure}\\ \smallskip\\
{\bf Definition 2.8:} Given a set $A$ of real numbers, the {\bf outer measure}
$m^*A$ of $A$ is the extended real number defined by
\[ m^*A = \inf_{A\subseteq \cup I_n} \sum l(I_n), \]
where the $\{ I_n\}$ are countable collections of open intervals that cover~$A$,
and $l(I_n)$ is the length of the interval~$I_n$.
\\ \smallskip\\
{\bf Observation 2.5:} $m^*$ is a set function, $m^* \emptyset = 0$, $m^* A\leq m^* B$
if $A\subseteq B$, and the outer measure of a set consisting of a single point is zero.
\\ \smallskip\\
{\bf Proposition 2.6:} $m^*(I) = l(I)$ if $I$ is an interval.
Proof in~\cite{royden}.\\ \smallskip \\
{\bf Proposition 2.7:} {\bf Countable subadditivity} of $m^*$:
$m^*(\cup A_n) \leq \sum m^*A_n$ for any countable collection $\{ A_n\}$
of sets of real numbers. Proof in~\cite{royden}.\\ \smallskip \\
{\bf Corollary 2.2:} $m^* A = 0$ if $A$ is countable.
\\ \smallskip\\
{\bf Observation 2.6:} The Cantor set is an example of an uncountable set with outer
measure zero. Proof in~\cite{rudin}. \\ \smallskip\\
{\bf Proposition 2.8:} $m^*$ is translation invariant, i.e., $m^*(E+y) = m^*E$ for any
set $E$ of real numbers and any number $y$, where $E+y = \{x+y: x\in E\}$.
\\ \smallskip\\
{\bf Proof:} If $\{I_n\}$ is a countable collection of open intervals that covers $E$,
then $\{I_n + y\}$ is a countable collection of open intervals that covers $E+y$.
Since $l(I_n) = l(I_n +y)$ for each $n$, then $m^*(E+y)\leq m^*E$.
Similarly, if $\{I_n\}$ is a countable collection of open intervals that covers $E+y$,
then $\{I_n -y\}$ covers $E$, $l(I_n) = l(I_n -y)$, and therefore $m^*E\leq m^*(E+y)$.
Thus, $m^*E = m^*(E+y)$.
\\ \smallskip\\
{\bf Proposition 2.9:} For any set $A$ of real numbers and any $\epsilon >0$, there is an
open set $O$ with $A\subseteq O$ and $m^*O\leq m^*A + \epsilon$
($m^*O < m^*A + \epsilon$ if $m^*A < \infty$). In addition, there is a
set $G$ that is the intersection of a countable collection of open sets with
$A\subseteq G$ and $m^*G = m^*A$. \\ \smallskip\\
{\bf Proof:} From the definition of the outer measure there is
a countable collection $\{ I_n\}$ of open intervals that covers~$A$ with
$\sum l(I_n)\leq m^*A + \epsilon$. With $O=\cup I_n$, then $O$ is open, $A\subseteq O$,
and $m^*O\leq \sum l(I_n)\leq m^*A +\epsilon$, which proves the first part.
Now let $k>0$ be an integer. Then from the first part there is an open set $O_k$,
$A\subseteq O_k$, with $m^*O_k \leq m^*A + \frac{1}{k}$. With $G=\cap_{k=1}^\infty O_k$,
then for any integer $n>0$ we have $A\subseteq G \subseteq O_n$ and therefore
$m^*A\leq m^*G \leq m^*O_n \leq m^*A + \frac{1}{n}$. Letting $n\rightarrow\infty $,
then $m^*A\leq m^*G \leq m^*A $. Thus $m^*A = m^*G$, which proves the second part.
\\ \medskip\\
{\bf \large Measurable sets, Lebesgue measure}\\ \smallskip\\
{\bf Definition 2.9 (Carath\'{e}odory's criterion):} A set $E$ of real numbers is said to be
{\bf (Lebesgue) measurable} if for every set $A$ of real numbers, then 
\[ m^*A = m^*(A\cap E) + m^*(A\cap \tilde {E}),\]
$\tilde{E}$ the complement of $E$ in {\bf R}, i.e., $\tilde{E} = {\bf R}\setminus E$,
{\bf R} the set of real numbers.
\\ \smallskip\\
{\bf Observation 2.7:} Clearly $\tilde {E}$ is measurable if and only if $E$ is, and
$\emptyset$ and the set {\bf R} of real numbers are measurable. As for an example of a
nonmeasurable set, a rather complex one is presented in~\cite{royden}. Finally, note
that it is always true that $m^*A \leq m^*(A\cap E) + m^*(A\cap \tilde {E})$,
thus $E$ is measurable if and only if for every set $A$ we have 
$m^*A \geq m^*(A\cap E) + m^*(A\cap \tilde {E})$.
\\ \smallskip\\
{\bf Proposition 2.10:} If $m^*E = 0$, then $E$ is measurable.
\\ \smallskip\\
{\bf Proof:} For any set $A$, since $A \supseteq A\cap \tilde {E}$, $m^*E = 0$, and
$E \supseteq A\cap E$, then
\[m^*A \geq m^*(A\cap \tilde {E}) =
m^*(A\cap \tilde {E})+m^*(E) \geq m^*(A\cap \tilde {E})+m^*(A\cap E).\]
{\bf Proposition 2.11:} The collection of measurable sets is a $\sigma$-algebra
on~{\bf R}. Proof in~\cite{royden}. \\ \smallskip\\
{\bf Definition 2.10:} The {\bf Lebesgue measure} $m$ is the set function obtained
by restricting the set function $m^*$ to the collection of (Lebesgue) measurable~sets.
\\ \smallskip\\
{\bf Proposition 2.12:} {\bf Countable subadditivity} of $m$:
$m(\cup A_n) \leq \sum mA_n$ for any countable collection $\{ A_n\}$
of measurable sets (Proposition~2.7).\\
{\bf Countable additivity} of $m$: $m(\cup A_n) = \sum mA_n$ if the sets in $\{ A_n\}$ as above
are pairwise disjoint. Proof in~\cite{royden}. \\ \smallskip \\
{\bf Observation 2.8:} It is in the proof of countable additivity of the Lebesgue
measure $m$ that
%the definition of a measurable set based on
Carath\'{e}odory's criterion plays a major role. On the other hand, the importance
of the countable additivity of $m$ is immediately apparent in the proofs of the
two parts of the following very useful proposition.
\\ \smallskip\\
{\bf Proposition 2.13 (Nested sequences of measurable sets Lemma):}\\
1. Given a countable collection of measurable sets $\{E_n\}$ with
$E_{n+1}\subseteq E_n$ for each $n$, $mE_1<\infty$, then
$ m(\cap_{i=1}^{\infty}\, E_i) = \lim_{n\rightarrow\infty} mE_n. $\\
2. Given a countable collection of (not necessarily measurable) sets $\{E_n\}$
with $E_{n+1}\supseteq E_n$ for each $n$, then
$ m^*(\cup_{i=1}^{\infty}\, E_i) = \lim_{n\rightarrow\infty} m^*E_n. $\\
Proof in~\cite{royden} for the first part. In~\cite{rana} for the second part using
Proposition~2.9. As mentioned above, in the proofs of both parts the countable additivity
of~$m$ (Proposition~2.12) is~used.
\\ \smallskip\\
{\bf Proposition 2.14:} Every Borel set is measurable.  Proof in~\cite{royden}.
\\ \smallskip\\
{\bf Observation 2.9:} A rather complex example of a measurable set that is not Borel
is presented in~\cite{burk}. \\ \smallskip\\
{\bf Proposition 2.15 (Equivalent conditions for a measurable set):} Let $E$ be a set
of real numbers. Then the following five conditions are equivalent:\\
i. $E$ is measurable.\\
ii. Given $\epsilon >0$, there is an open set $O$, $O\supseteq E$ with
$m^*(O\setminus E)<\epsilon$.\\
iii. Given $\epsilon >0$, there is a closed set $F$, $F\subseteq E$ with
$m^*(E\setminus F)<\epsilon$.\\
iv. There is a set $G$ that is the intersection of a countable collection of
open sets, $G\supseteq E$ with $m^*(G\setminus E)=0$.\\
v. There is a set $F$ that is the union of a countable collection of
closed sets, $F\subseteq E$ with $m^*(E\setminus F)=0$. \\ \smallskip\\
{\bf Proof:} We only prove i $\Leftrightarrow$ ii. Proofs of all cases in~\cite{royden2}.\\
i $\Rightarrow$ ii: $E$ is measurable.\\
Case 1: $mE < \infty$.\\
There exists $O$ open such that $E\subseteq O$, $mO < mE +\epsilon$ (Proposition~2.9).
$mE < \infty$ then implies $mO-mE < \epsilon$. Accordingly,
\begin{eqnarray*}
\epsilon & > & m(O\cap E) + m(O\cap\tilde{E})-mE \\
& = & mE + m(O\setminus E) - mE \\
& = & m(O\setminus E).
\end{eqnarray*}
Case 2: $mE = \infty$.\\
For each integer $n>0$ let $E_n=E\cap((n-1,n]\cup (-n,-n+1])$. Then each $E_n$ is
measurable, $mE_n <\infty$, and $E=\cup E_n$. From Case 1 above it follows that
there is $O_n$ open such that $E_n\subseteq O_n$ and $m(O_n\setminus E_n)<\epsilon/2^n$.
With $O=\cup\, O_n$, then $O$ is open, $E\subseteq O$ and $O\setminus
E\subseteq \cup(O_n\setminus E_n)$. Thus,
\[ m(O\setminus E) \leq m(\cup (O_n\setminus E_n)) \leq \sum m(O_n\setminus E_n)<
\sum \epsilon/2^n =\epsilon. \]
ii $\Rightarrow$ i: Given $\epsilon >0$, there is an open set $O$, $O\supseteq E$ with
$m^*(O\setminus E)<\epsilon$.\\
 $O$ open implies $O$ is measurable. Thus, for any set $A$,
we have $m^*A = m^*(A\cap O) + m^*(A\cap \tilde{O})$. Accordingly,
\begin{eqnarray*}
m^*A + \epsilon & > & m^*A + m^*(O\cap\tilde{E})\\
\ & = & m^*(A\cap O) + m^*(A\cap \tilde{O})+m^*(O\cap\tilde{E})\\
\ & \geq & m^*(A\cap E) + m^*(A\cap \tilde{O}\cap\tilde{E})+m^*(O\cap\tilde{E}\cap A)\\
\ & = & m^*(A\cap E) + m^*(A\cap\tilde{E}).
\end{eqnarray*}
$\epsilon$ arbitrary implies $m^*A \geq m^*(A\cap E) + m^*(A\cap\tilde{E})$.
Thus, $E$ is measurable.\\ \smallskip \\
{\bf Corollary 2.3:} Every (Lebesgue) measurable set is the union of a
Borel set and a set of (Lebesgue) measure zero.\\ \smallskip \\
{\bf Proposition 2.16:} The translate of a measurable set is measurable, i.e., if $E$ is
a measurable set, then $E+y$ is measurable for any number~$y$.\\ \smallskip\\
{\bf Proof:} For any set $A$, setting $B=A-y$, $F=E+y$, and noting $(B\cap E)+y = A\cap F$,
$(B\cap\tilde{E})+y= A\cap\tilde{F}$, by Proposition~2.8, then
\[ m^*A=m^*B = m^*(B\cap E) + m^*(B\cap\tilde{E}) = m^*(A\cap F) + m^*(A\cap \tilde{F}). \]
Thus, $F=E+y$ is measurable.\\ \medskip\\
{\bf \large Measurable functions, Step functions, Simple functions}\\ \smallskip\\
{\bf Definition 2.11:} Let $f$ be an extended real-valued function defined on a (Lebesgue)
measurable set. Then $f$  is said to be {\bf (Lebesgue) measurable} if the set
$\{x\, |\, f(x)>a\}$ is (Lebesgue) measurable for every real number~$a$.\\ \smallskip\\
{\bf Proposition 2.17 (Equivalent conditions for a measurable function):} Let $f$ be
an extended real-valued function of measurable domain.
Then the following conditions are equivalent:\\
i. $\{x\, |\, f(x)>a\}$ is measurable for every real number~$a$.\\
ii. $\{x\, |\, f(x)\geq a\}$ is measurable for every real number~$a$.\\
iii. $\{x\, |\, f(x)<a\}$ is measurable for every real number~$a$.\\
iv. $\{x\, |\, f(x)\leq a\}$ is measurable for every real number~$a$.\\
Proof in~\cite{royden} and \cite{rudin}.\\ \smallskip\\
{\bf Observation 2.10:} Conditions ii, iii, iv can be used instead of condition~i to
define a measurable function. We note that if a real-vaued function $f$ defined on a
closed or open interval $I$ is {\bf continuous}~\cite{rudin}, then $f$ is measurable since
the set in condition~i is relatively open in~$I$~\cite{rudin}. Also the restriction of
a measurable function to a measurable subset of its domain, is measurable.\\ \smallskip\\
{\bf Definition 2.12:} Given real numbers $a$, $b$, $a<b$, and an integer $n>0$, by a
{\bf partition} or {\bf subdivision} of~$[a,b]$ we mean a finite set of points
$P= \{\xi_0, \xi_1,\ldots,\xi_n\}$ with $a=\xi_0<\xi_1 < \ldots < \xi_n = b$.
A function $\psi: [a,b]\rightarrow$ {\bf R} is called a {\bf step function} (on $[a,b]$)
if for an integer $n>0$, there are numbers $c_i$, $i=1,\ldots,n$, and a partition or
subdivision $P= \{\xi_0, \xi_1,\ldots,\xi_n\}$ of $[a,b]$, such that $\psi(a)=\psi(\xi_0) =
c_1$, $\psi(x)=c_i$, $\xi_{i-1}<x\leq\xi_i$, $i=1,\ldots,n$. \\ \smallskip\\
{\bf Definition 2.13:} Given sets $A$, $X$ of real numbers, $A\subseteq X$,
the {\bf characteristic function} $\chi_A$ of $A$ on $X$, $\chi_A: X\rightarrow \{0,1\}$, is
defined by
\[ \chi_A(x) = \left\{ \begin{array}{ll}
            1 & x\in A\\
            0 & x\in X\setminus A.
                    \end{array}
            \right. \]
{\bf Definition 2.14:} Given a measurable set $E$, an integer $n>0$, and for $i=1,\ldots,n$,
nonzero numbers $c_i$, measurable sets $E_i\subseteq E$, and characteristic functions
$\chi_{E_i}$ of $E_i$ on~$E$, a function $\varphi: E \rightarrow$ {\bf R} defined by
$\varphi(x) = \sum_{i=1}^n c_i\chi_{E_i}(x)$ for $x$ in $E$ is called a
{\bf simple function} on $E$. \\ \smallskip\\
{\bf Observation 2.11:} Step functions are measurable and a function is simple if and
only if it is measurable and assumes only a finite number of values. We note that
the representation of a simple function $\varphi$ is not unique. However it does have
a so-called {\bf canonical representation}: for some integer $m>0$, $\varphi(x) =
\sum_{i=1}^m a_i\chi_{A_i}(x)$ for $x$ in $E$, where $\{a_1,\ldots,a_m\}$ is the set
of distinct nonzero values of $\varphi$, and $A_i=\{x\in E:\,\varphi(x) = a_i\}$,
$i=1,\ldots,m$. This representation is characterized by the fact that the $a_i$'s
are distinct and nonzero, and the $A_i$'s are pairwise disjoint.\\ \smallskip\\
{\bf Proposition 2.18:} If $f$ is a measurable function, then the function $|f|$
is measurable.  Proof in~\cite{rudin}.\\ \smallskip\\
{\bf Proposition 2.19:} Let $f$ and $g$ be measurable real-valued (not extended)
functions defined on the same domain~$X$, $c$ a constant, and $F$
a continuous real-valued function on~{\bf R}$^2$. Then the function $h$ defined
by $h(x) = F(f(x),g(x))$, $x\in X$, is measurable. In particular $f+g$, $fg$,
$cf$, $f+c$, $f-g$ are measurable; $f/g$ is measurable if $g\not=0$ on~$X$.
Proofs in~\cite{royden} and~\cite{rudin}.
\\ \smallskip\\
{\bf Observation 2.12:} For extended real-valued measurable functions $f$ and $g$,
the function $fg$ is still measurable.
However in order for $f+g$ to be measurable, $f+g$ must be given the same value at
points where it is undefined, unless these points form a set of measure zero in which
case it makes no difference which values $f+g$ is given on the set.\\ \smallskip\\
{\bf Proposition 2.20:} Let $\{ f_n\}$ be a sequence of measurable functions defined on
the same domain~$X$ and let $N>0$ be an integer. The functions defined for each
$x\in X$ by $\sup_{n\leq N} f_n(x)$, $\sup_n f_n(x)$,
$\limsup_{n\rightarrow\infty} f_n(x)$, $\inf_{n\leq N} f_n(x)$, $\inf_n f_n(x)$,
$\liminf_{n\rightarrow\infty} f_n(x)$, are all measurable. Proofs in~\cite{royden}
and~\cite{rudin}.\\ \smallskip\\
{\bf Definition 2.15:} Given a set of real numbers $X$, a property associated with points
in $X$ is said to hold {\bf almost everywhere} ({\bf a.e.} for short) in $X$ if the
set of points in $X$ where the property fails has measure zero.\\ \smallskip\\
{\bf Proposition 2.21:} Let $f$ and $g$ be functions defined on the same measurable
domain~$X$. If $f$ is measurable and $f=g$ a.e., then $g$ is measurable.\\ \smallskip\\
{\bf Proof:} Let $a$ be any real number. $f=g$ a.e. means the set
$E=\{ x:f(x)\not= g(x)\}$ has measure zero. Thus $E$ and the set $\{ x\in E: g(x)>a\}
\subseteq E$ are measurable (Proposition~2.10). $X\setminus E$ must then be measurable and
the set $\{ x\in X\setminus E: f(x)>a\}$ is then measurable as $f$ is measurable.
Thus, since \[ \{x: g(x)>a\} = \{ x\in X\setminus E: f(x)>a\}\cup \{ x\in E: g(x)>a\} \]
it follows that $\{x:g(x)>a\}$ is measurable and therefore $g$ is measurable.
\\ \smallskip\\
{\bf Observation 2.13:} The definition of a measurable function and 1 of Proposition~2.13
allow the following proposition to be true.
\\ \smallskip\\
{\bf Proposition 2.22 (Egoroff's Theorem):} Let $\{ f_n\}$ be a sequence of measurable
functions on a measurable set~$E$, $mE<\infty$, that
converge a.e. to a real-valued (not extended) function $f$ on~$E$,
i.e., there is a set $B\subseteq E$ such that $mB=0$, $f_n\rightarrow f$ pointwise
on~$E\setminus B$. Then for every $\delta >0$ there is a closed set $F\subseteq E$ such
that $m(E\setminus F)<\delta$ and $f_n\rightarrow f$ uniformly on~$F$.\\ \smallskip\\
{\bf Proof:} Let $m>0$ be an integer. For each integer $n>0$ let
\[ G_n^m = \{ x\in E: |f_n(x)-f(x)| \geq 1/m\},\]
and for each integer $N>0$ set
\[ E_N^m = \cup_{n=N}^{\infty}\, G_n^m =
\{ x\in E: |f_n(x)-f(x)|\geq 1/m \ \mathrm{for\ some\ } n\geq N\}. \]
We note $f$ is measurable (Proposition~2.20 and Proposition~2.21). It follows that each
$E_N^m$ is measurable and of finite measure, $E_{N+1}^m \subseteq E_N^m$,
and for each $x\in E\setminus B$
there must be some $N$ for which $x\not\in E_N^m$, since $f_n(x)\rightarrow f(x)$.
Thus $\cap E_N^m\subseteq B$ and must therefore have measure zero.
It follows then that $\lim_{N\rightarrow\infty} mE_N^m=0$ (1 of Proposition~2.13).
Hence there exists $N$ such that
\[ mE_N^m 
=m(\{ x\in E: |f_n(x)-f(x)|\geq 1/m \ \mathrm{for\ some\ } n\geq N\})
< \delta/2^{m+1}.\]
Letting $A^m=E_N^m$, then $A^m$ is measurable, $mA^m < \delta/2^{m+1}$ and
\[E\setminus A^m = \{ x\in E: |f_n(x)-f(x)| < 1/m\ \mathrm{for\ all\ } n\geq N\}. \]
Now let $A=\cup_{m=1}^{\infty} A^m$.
Then $mA \leq \sum_{m=1}^{\infty}mA^m<\sum_{m=1}^{\infty}\delta/2^{m+1}=\delta/2.$\\
Given $\epsilon>0$ choose integer $m>0$ with $1/m<\epsilon$. For some $N$ and for all
$x\in E\setminus A$ then $x\in E\setminus A^m$ and $|f_n(x)-f(x)|<1/m<\epsilon$
for all~$n\geq N$.\\ Thus, $f_n\rightarrow f$ uniformly on~$E\setminus A$.\\
Finally, let $G=E\setminus A$. Clearly $G$ is measurable
and~$m(E\setminus G)=mA<\delta/2$.\\
Thus, there exists a closed set $F$, $F\subseteq G$ with $m(G\setminus F)<\delta/2$
(Proposition~2.15). It then follows that $m(E\setminus F) =
m(E\setminus G)+ m(G\setminus F)<\delta/2 + \delta/2 = \delta$ and $f_n\rightarrow f$
uniformly on~$F$.\\ \smallskip\\
Proof also in \cite{halmos} for spaces and measures more general than {\bf R} and
the Lebesgue measure.\\ \smallskip\\
{\bf Observation 2.14:} The assumption $mE<\infty$ is necessary in the above proposition.
To see this, let $E=$ {\bf R} and 
\[ f_n(x) = \left\{ \begin{array}{ll}
            1 & x > n\\
            0 & x \leq n.
                    \end{array}
            \right. \]
Clearly $f_n\rightarrow 0$ pointwise on $E$. However, for any integer $N>0$, $0<\epsilon<1$,
the set $\{x\in E: f_N(x) \geq \epsilon\}=(N,\infty)$ is of infinite measure. Thus, the
uniform convergence of $f_n$ to $0$ as proposed in the proposition can not occur.
\\ \smallskip\\
{\bf Definition 2.16:} Given a function $f$ defined on a set $E$,
the {\bf positive part $f^+$} of $f$ and the {\bf negative part $f^-$} of $f$
are the functions defined respectively by
$f^+(x) = \max\{f(x),0\}$, and $f^-(x) = \max\{-f(x),0\}$, $x\in E$.
\\ \smallskip\\
{\bf Proposition 2.23 (Approximation of a measurable function by simple functions):}
Let $f$ be a real-valued (not extended) measurable function on a measurable
set~$E$. Then there exists a sequence of simple functions $\{s_n\}$ on $E$
such that $s_n\rightarrow f$ pointwise on~$E$. Since $f = f^+ - f^-$,
and $f^+$ and $f^-$ are measurable because $f$ is, then $\{s_n\}$
can be chosen so that $s_n= s_n^u - s_n^l$, where $s_n^u$ and $s_n^l$ are simple
functions on $E$ such that $s_n^u\rightarrow f^+$, $s_n^l\rightarrow f^-$ pointwise
on $E$, and $s_n^u$ and $s_n^l$ increase monotonically to $f^+$ and~$f^-$,
respectively. If $f$ is bounded, then $\{s_n\}$ can also be chosen to converge
uniformly to~$f$ on~$E$. Proof in~\cite{rudin}.
\\ \smallskip\\
{\bf Proposition 2.24 (Lusin's Theorem):} Let $f$ be a real-valued (not extended) measurable
function on a measurable set $E$. Then given $\epsilon>0$, there exists a closed set
$F\subseteq E$ with $m(E\setminus F)<\epsilon$ such that $f|_F$ is continuous.
\\ \smallskip\\
{\bf Proof:} First we prove the proposition for $f$ a simple function on~$E$.
Accordingly, for some integer $n>0$, assume $f(x) = \sum_{i=1}^{n} c_i\chi_{E_i}(x)$
for $x$ in $E$, $E_i\subseteq E$ for each~$i$, the canonical representation of~$f$.
In addition let $E_0 = E\setminus \cup_{i=1}^{n}E_i$. Clearly the $E_i$'s are
pairwise disjoint. Given $\epsilon >0$,
since each $E_i$ is measurable there exists a closed set $F_i\subseteq E_i$ such
that $m(E_i\setminus F_i) < \epsilon/(n+1)$, $i=0,\ldots,n$ (Proposition~2.15).
Accordingly, $F=\cup_{i=0}^n\,F_i$ is closed, and
\[m(E\setminus F) = m(\cup_{i=0}^n\, E_i\setminus \cup_{i=0}^n\, F_i)
=m(\cup_{i=0}^n(E_i\setminus F_i)) = \sum_{i=0}^n m(E_i\setminus F_i)< \epsilon.\]
Now, to show $f|_F$ is continuous we show that if $\{x_k\}$, $x$ in $F$ are such
that $x_k\rightarrow x$, then $f(x_k)\rightarrow f(x)$.\\
We note that for some unique $j$, $0\leq j\leq n$, it must be that $x\in F_j$.
Thus, since $f$ is constant on $F_j$, it then suffices to show that for some integer
$K>0$, $x_k$ is in $F_j$ for $k\geq K$. If this is not the case and since there is
a finite number of $F_i$'s then for some $l$, $0\leq l\leq n$, $l\not = j$, there is
a subsequence $\{x_{k_m}\}$ of $\{x_k\}$ all contained in~$F_l$. But then
$x_{k_m}\rightarrow x$ so that $x$ is in $F_l$, a contradiction. \\ \smallskip\\
Now we prove the proposition for a general $f$.\\
Case 1: $mE<\infty.$\\
Let $s_n$ be simple functions such that $s_n\rightarrow f$ pointwise on $E$
(Proposition~2.23). Given $\epsilon>0$, as established above for simple functions
on $E$, for each $n$ there exists a closed set $F_n\subseteq E$ with
$m(E\setminus F_n) <\epsilon/2^{n+1}$ such that $s_n|_{F_n}$ is continuous.
In addition, since $mE<\infty$, there exists a closed set $F_0\subseteq E$
such that $m(E\setminus F_0)<\epsilon/2$ and $s_n\rightarrow f$ uniformly
on~$F_0$ (Proposition~2.22~(Egoroff's Theorem)).\\
Finally let $F=\cap_{n=0}^{\infty}\,F_n$. Then $F$ is closed and
\[m(E\setminus F) = m(\cup_{n=0}^{\infty}\,(E\setminus F_n))\leq \sum_{n=0}^{\infty}
m(E\setminus F_n)<\sum_{n=0}^{\infty}\epsilon/2^{n+1}=\epsilon. \]
Since $s_n|_{F_n}$ is continuous so must be $s_n|_F$. And since $s_n \rightarrow f$
uniformly on $F\subseteq F_0$ then $f|_F$ must be continuous (proof in \cite{rudin}).
\\ \smallskip\\
Case 2: $mE = \infty$.\\
For each integer $n>0$, let $\tilde{E}_n=(n-1,n]\cup (-n,-n+1]$ and $E_n=E\cap\tilde{E}_n$.
Then each $E_n$ is measurable, $mE_n <\infty$, and $E=\cup_{n=1}^{\infty} E_n$. From
Case~1 above, given $\epsilon>0$, it follows that there is a closed set $F_n\subseteq E_n$
with $m(E_n\setminus F_n)<\epsilon/2^n$ such that $f|_{F_n}$ is continuous.
Let $F=\cup_{n=1}^{\infty}\, F_n$. Then
\begin{eqnarray*}
m(E\setminus F) &=& m(\cup_{n=1}^{\infty}\, E_n\setminus \cup_{n=1}^{\infty}\, F_n)
= m(\cup_{n=1}^{\infty}(E_n\setminus F_n))\\
&=& \sum_{n=1}^{\infty} m(E_n\setminus F_n) < \sum_{n=1}^{\infty}\epsilon/2^n=\epsilon.
\end{eqnarray*}
We show $F$ is closed. Let $x$ be a limit point of~$F$ so that for $\{x_k\}$ in $F$,
then $x_k\rightarrow x$. Clearly for some integer $j>0$ it must be that $x\in \tilde{E}_j$.
It suffices to show that for some integer $K>0$, $x_k$ is in $F_j$ for~$k\geq K$ so that
$x$ is also~in $F_j\subseteq F$. With $\tilde{E}_n=\emptyset$ for $n\leq 0$, from the
definition of the $\tilde{E}_n$'s a neighborhood of the point $x$ exists that does not
intersect $\tilde{E}_i$, $i>j+1$ or $i<j-1$. Thus, with $F_0=\emptyset$, for $k$ large
enough the points $x_k$ can only be in $F_{j-i}$, $F_j$ and $F_{j+1}$. If it is not the
case that $K$ as described exists, then there must be a subsequence $\{x_{k_m}\}$ of $\{x_k\}$,
all of it contained in either $F_{j-1}$ or~$F_{j+1}$. But then $x_{k_m} \rightarrow x$ so that
$x$ is in either $F_{j-1}$ or $F_{j+1}$, a contradiction. (Actually showing that $x$ is in the
union of $F_{j-1}$, $F_j$ and $F_{j+1}$, would have sufficed).\\
Now, to show $f|_F$ is continuous we show that if $\{x_k\}$, $x$ in $F$ are such
that $x_k\rightarrow x$, then $f(x_k)\rightarrow f(x)$.\\
We note that for some unique $j>0$, it must be that $x\in F_j$. Thus,
since $f$ is continuous on $F_j$, it then suffices to show that for some integer
$K>0$, $x_k$ is in $F_j$ for $k\geq K$. If this is not the case, again with $F_0=\emptyset$,
an argument, similar to the one used above for proving $F$ is closed, can be used to get the
same contradiction that $x$ is in either $F_{j-1}$ or $F_{j+1}$.
\\ \medskip\\
{\bf \large The Riemann integral}\\ \smallskip\\
{\bf Definition 2.17:} Let $[a,b]$ be an interval and $f$ a bounded real-valued function
defined on~$[a,b]$. Given an integer $n>0$, and $P = \{ \xi_0, \xi_1,\ldots,\xi_n\}$,
a partition or subdivision of~$[a,b]$ (Definition 2.12), for $i=1,\ldots,n$, we define
\begin{eqnarray*}
m_i &=& \inf f(x),\ \xi_{i-1} \leq x\leq \xi_i,\\
M_i &=& \sup f(x),\ \xi_{i-1} \leq x\leq \xi_i,\ \mathrm{and}\\
L(P,f) &=& \sum_{i=1}^n\,m_i(\xi_i-\xi_{i-1}),\\
U(P,f) &=& \sum_{i=1}^n\,M_i(\xi_i-\xi_{i-1}).
\end{eqnarray*}
Then we define the {\bf lower Riemann integral} and the
{\bf upper Riemann integral} of $f$ over $[a,b]$, respectively, by
\[ {\cal R} \underline{\int_a^b} f(x)dx = \sup_P L(P,f), \]
\[ {\cal R} \overline{\int_a^b} f(x)dx = \inf_P U(P,f), \]
where the infimum and supremum are taken over all partitions $P$ of~$[a,b]$.\\
If the two are equal, $f$ is said to be {\bf Riemann integrable} over $[a,b]$,
and the common value is then called the {\bf Riemann integral} of $f$ over $[a,b]$
and denoted by \[ {\cal R} \int_a^b f(x)dx. \]
{\bf Observation 2.15:} Since $f$ is bounded, there exist numbers $m$ and $M$, such
that $m\leq f(x)\leq M$, $x\in [a,b]$. Thus, for every partition $P$, it must be
that $m(b-a)\leq L(P,f)\leq U(P,f)\leq M(b-a)$, so that the lower and upper Riemann
integrals of $f$ over $[a,b]$ are finite numbers.
\\ \smallskip\\
{\bf Proposition 2.25:}
$ {\cal R} \underline{\int_a^b} f(x)dx \leq {\cal R} \overline{\int_a^b} f(x)dx $.
Proof in~\cite{rudin}.\\ \smallskip\\
{\bf Observation 2.16:}
Let $\psi: [a,b]\rightarrow$ {\bf R} be a step function so that for an integer $n>0$,
there are numbers $c_i$, $i=1,\ldots,n$, and a partition $P= \{\xi_0, \xi_1,\ldots,\xi_n\}$
of $[a,b]$, such that
$\psi(a)=\psi(\xi_0) = c_1$, $\psi(x)=c_i$, $\xi_{i-1}<x\leq\xi_i$, \mbox{$i=1,\ldots,n$}.
Clearly $L(P,\psi)=U(P,\psi)$ and since $ L(P,\psi) \leq {\cal R} \underline{\int_a^b}
\psi(x)dx \leq {\cal R} \overline{\int_a^b} \psi(x)dx \leq U(P,\psi) $ (Proposition~2.25),
it must be that $\psi$ is Riemann integrable over $[a,b]$ and
${\cal R} \int_a^b \psi(x)dx = \sum_{i=1}^n c_i(\xi_i-\xi_{i-1})$.\\
From this it is then apparent that
\[ {\cal R} \underline{\int_a^b} f(x)dx = \sup_P L(P,f) =
\sup_{\psi\leq f} {\cal R} \int_a^b \psi(x)dx, \]
\[ {\cal R} \overline{\int_a^b} f(x)dx = \inf_P U(P,f) =
\inf_{\psi\geq f} {\cal R} \int_a^b \psi(x)dx, \]
where the $\psi$'s are all possible step functions on~$[a,b]$ satisfying the given conditions.
\\ \smallskip \\
{\bf Definition 2.18:} Given an interval $[a,b]$, let
$P= \{\xi_0, \xi_1,\ldots,\xi_n\}$ be a partition  of $[a,b]$. The number
\[ \mu(P) = \max_{i=1,\ldots,n} (\xi_i-\xi_{i-1}) \]
is called the {\bf mesh} of~$P$.\\
Let $f$ be a real-valued function defined on $[a,b]$.
Given a partition $P= \{\xi_0, \xi_1,\ldots,\xi_n\}$ of $[a,b]$,
a {\bf Riemann sum} of $f$ with respect to $P$ is a sum of the form
\[ S(P,f) = \sum_{i=1}^n f(t_i)(\xi_i-\xi_{i-1}),\]
where the choice of points
$t_1,\ldots,t_n$, $\xi_{i-1} \leq t_i \leq \xi_i$, $i=1,\ldots,n$, is arbitrary.
The Riemann sums of $f$ are said to converge to a finite number~$I$ as
\mbox{$\mu(P)\rightarrow 0$}, i.e., \[ I = \lim_{\mu(P)\rightarrow 0} S(P,f), \]
if given $\epsilon>0$, there exists $\delta>0$ such that for every partition $P$ with
mesh \mbox{$\mu(P)<\delta$} it must be that
\[ |S(P,f) - I| < \epsilon \]
(obviously for every choice of points
$t_1,\ldots,t_n$, $\xi_{i-1} \leq t_i \leq \xi_i$, $i=1,\ldots,n$).
\\ \smallskip \\
{\bf Proposition 2.26 (Riemann sums of $f$ that converge implies $f$ is bounded):}
If $\lim S(P,f)$ exists as $\mu(P)\rightarrow 0$, then $f$ is bounded
on~$[a,b]$. Proof in~\cite{pugh} and~\cite{richardson}.
\\ \smallskip\\
{\bf Proposition 2.27 (Riemann sums of $f$ converge if and only $f$ is Riemann
integrable):} Let $[a,b]$ be an interval and $f$ a bounded real-valued function
defined on $[a,b]$. Then $f$ is Riemann integrable over $[a,b]$ if and only~if
\[ I = \lim_{\mu(P)\rightarrow 0} S(P,f) \]
exists. If this is the case, then $I$~equals~${\cal R} \int_a^b f(x)dx$.
Proof in~\cite{rudin} and~\cite{wade}.
\\ \smallskip\\
{\bf Observation 2.17:} Given an interval $[a,b]$ and a set $A\subseteq [a,b]$, ideally the
characteristic function $\chi_A$ of $A$ on $[a,b]$, $\chi_A: [a,b]\rightarrow \{0,1\}$,
defined by
\[ \chi_A(x) = \left\{ \begin{array}{ll}
            1 & x\in A\\
            0 & x\in [a,b]\setminus A
                    \end{array}
            \right. \]
should be (Riemann) integrable over $[a,b]$, especially if $A$ is measurable, and its integral
over $[a,b]$ should equal the (outer) measure of~$A$. However, if $A$ is the set of rational
numbers in $[a,b]$, which is measurable with~$mA=0$, we see that
${\cal R} \underline{\int_a^b} \chi_A(x)dx = 0$ and
${\cal R} \overline{\int_a^b} \chi_A(x)dx = b-a$, not the ideal situation.\\ \medskip\\
{\bf \large The Lebesgue integral}\\ \smallskip\\
{\bf Definition 2.19:} Given a measurable set $E$,
let $\varphi(x) = \sum_{i=1}^m a_i\chi_{A_i}(x)$ be the canonical representation
of a simple function $\varphi$ on $E$, where for some integer $m>0$,
$\{a_1,\ldots,a_m\}$ is the set of distinct nonzero values of $\varphi$,
and $A_i=\{x\in E:\,\varphi(x) = a_i\}$, $i=1,\ldots,m$. We define the
{\bf Lebesgue integral} of $\varphi$ over $E$ as the extended real number
\[ \int_E \varphi(x)dx = \sum_{i=1}^m a_i mA_i. \]
{\bf Observation 2.18:} A consequence of the following two propositions is that if
$\varphi(x)= \sum_{i=1}^n c_i\chi_{E_i}(x)$ is any representation of a simple function
$\varphi$ on a measurable set $E$, then the Lebesgue integral of $\varphi$ over~$E$
(Definition 2.19) can be computed directly from the representation, i.e.,
by computing~$\sum_{i=1}^n c_i mE_i$.  \\ \smallskip\\
{\bf Proposition 2.28:} Let $\varphi(x)= \sum_{i=1}^n c_i\chi_{E_i}(x)$ be a representation
of a simple function $\varphi$ on a measurable set~$E$, with $E_i\cap E_j=\emptyset$ for
$i\not=j$ (not necessarily the canonical representation of $\varphi$). Then 
\[ \sum_{i=1}^n c_i mE_i = \int_E \varphi(x)dx.\]
Proof in~\cite{royden} for $E_i$'s of finite measure.
Same proof for the general case.
\\ \smallskip\\
{\bf Proposition 2.29:} Let $\varphi$ and $\psi$ be simple functions on a measurable
set~$E$. Then for any real numbers $a$ and $b$ we must have
\[ \int_E (a\varphi + b\psi)(x)dx = a\int_E\varphi(x)dx + b\int_E\psi(x)dx, \]
and, if $\varphi \geq \psi$ a.e., then $\int_E\varphi(x)dx \geq \int_E\psi(dx)$.
\\ \smallskip \\
Proof in~\cite{royden} using Proposition~2.28 for $E_i$'s of finite
measure. Proof essentially the same for the general case. \\ \smallskip\\
{\bf Corollary 2.4:} Let $\varphi(x)= \sum_{i=1}^n c_i\chi_{E_i}(x)$ be any representation
of a simple function $\varphi$ on a measurable set~$E$, the $E_i$'s not necessarily
pairwise disjoint. Then
\[ \sum_{i=1}^n c_i mE_i = \int_E \varphi(x)dx.\]
{\bf Proof:} Apply the first part of Proposition~2.29 to
$\varphi(x) = \sum_{i=1}^n c_i\chi_{E_i}(x)$. \\ \smallskip\\
{\bf Definition 2.20:} Given a measurable set $E$, let $f$ be a measurable nonnegative
function on~$E$. We define the {\bf Lebesgue integral} of $f$ over $E$ as the
extended real number 
\[\int_E f(x)dx = \sup_{\varphi\leq f} \int_E \varphi(x)dx, \]
where the $\varphi$'s are all possible simple functions on~$E$ satisfying the given condition.
\\ \smallskip\\
{\bf Definition 2.21:} Given a measurable set $E$, let $f$ be a measurable function on~$E$.
With $f^+$ and $f^-$ as the positive and negative parts of $f$ (Definition~2.16),
we define the {\bf Lebesgue integral} of $f$ over $E$ as the extended real number
\[\int_E f(x)dx = \int_E f^+(x)dx - \int_E f^-(x)dx, \]
if at least one of the integrals $\int_E f^+(x)dx$, $\int_E f^-(x)dx$ (Definition~2.20)
is finite.\\
If $\int_E f(x)dx$ is finite, then $f$ is said to be {\bf Lebesgue integrable} over~$E$.
\\ \smallskip\\
{\bf Proposition 2.30:} Let $f$ and $g$ be Lebesgue integrable functions over a measurable
set~$E$, and $c$ a real number. Then\\
i. $cf$ is Lebesgue integrable over~$E$ with $\int_E cf(x)dx = c\,\int_E f(x)dx$.\\
ii. $f+g$ is Lebesgue integrable over~$E$ with
\[\int_E (f+g)(x)dx = \int_E f(x)dx + \int_E g(x)dx.\]
iii. If $f\leq g$ a.e., then $\int_E f(x)dx \leq \int_E g(x)dx$.\\
iv. If $A,B\subseteq E$ are disjoint measurable sets, then
\[ \int_{A\cup B} f(x)dx = \int_A f(x)dx + \int_B f(x)dx. \]
Proofs~in \cite{royden} and~\cite{rudin}.\\ \smallskip\\
{\bf Proposition 2.31:} A measurable function $f$ is Lebesgue integrable over $E$ if
and only if $|f|$ is Lebesgue integrable over $E$, in which case
\[ |\int_E f(x)dx|\leq \int_E |f(x)|dx. \]
Also, if $0\leq f\leq g$ on~$E$ and $g$ is Lebesgue integrable over~$E$, then $f$ is
Lebesgue integrable over~$E$. In particular, if $|f|\leq g$ and $g$ is Lebesgue
integrable, then $|f|$, and therefore $f$, is Lesbegue integrable over~$E$.
\\ \smallskip\\
{\bf Proof:} The first part follows from $f=f^+ -f^-$, $|f| = f^+ + f^-$,
and~iv of Proposition~2.30. The inequality from $f\leq |f|$, $-f\leq |f|$, and~i and~iii of
Proposition~2.30. The rest from~Definition~2.20.
\\ \smallskip\\
{\bf Observation 2.19:}
Let $f$ be a measurable function on a measurable set $E$ with $mE$ finite,
and let $a$, $b$ be real numbers such that \mbox{$a\leq f(x)\leq b$} for $x\in E$.
By looking at $\int_E f^+(x)dx$ and $\int_E f^-(x)dx$ for the different possible
signs of $a$ and $b$, then it is evident that $a\,mE\leq\int_E f(x)dx\leq b\,mE$.
Accordingly, if $f$ is a measurable and bounded function on a measurable set
$E$ with $mE$ finite, since then for some $M>0$, $-M\leq f(x) \leq M$
for $x\in E$, it must be that $-M\,mE\leq\int_E f(x)dx\leq M\,mE$,
and therefore $f$ is Lebesgue integrable.
However, there is more to this situation as the following proposition shows.\\ \smallskip\\
{\bf Proposition 2.32 (Integrable equivalent to measurable):}
Let $f$ be a bounded function defined on a measurable set $E$ with
$mE$ finite. Let \[
L(f) = \sup_{\varphi\leq f} \int_E \varphi(x)dx,
\ \ \ \ 
U(f)= \inf_{\varphi\geq f} \int_E \varphi(x)dx,
\] 
where the $\varphi$'s are all possible simple functions on~$E$ satisfying the given
conditions. Then $L(f)=U(f)$ if and only if $f$ is measurable. Whenever $L(f)=U(f)$ then
$f$ is Lebesgue integrable and $\int_E f(x) dx=L(f)=U(f)$. Proof in~\cite{royden}.
\\ \smallskip\\
{\bf Proposition 2.33 (Riemann integrable implies Lebesgue integrable):} Let $f$ be
a bounded function on interval $[a,b]$. If $f$ is Riemann integrable
over $[a,b]$, then $f$ is measurable and Lebesgue integrable over $[a,b]$~with
\[ \int_{[a,b]} f(x)dx = {\cal R} \int_a^b f(x)dx. \]
{\bf Proof:} Since step functions are simple functions, then
\[ \sup_{\psi\leq f} {\cal R} \int_a^b \psi(x)dx \leq 
 \sup_{\varphi\leq f} \int_{[a,b]} \varphi(x)dx \leq 
 \inf_{\varphi\geq f} \int_{[a,b]} \varphi(x)dx \leq 
 \inf_{\psi\geq f} {\cal R} \int_a^b \psi(x)dx, \] 
where the $\psi$'s and the $\varphi$'s are all possible step functions and simple
functions on~$[a,b]$, respectively, satisfying the given conditions. Since $f$ is Riemann
integrable over $[a,b]$, then all the inequalities above are equalities so that
$f$ must be measurable and Lebesgue integrable over $[a,b]$ with
\[ \int_{[a,b]} f(x)dx =
 \sup_{\varphi\leq f} \int_{[a,b]} \varphi(x)dx =
 \inf_{\varphi\geq f} \int_{[a,b]} \varphi(x)dx =
{\cal R} \int_a^b f(x)dx \]
by Proposition~2.32. \\ \smallskip\\
{\bf Proposition 2.34:} Let $f$ be a measurable function on a measurable set~$E$.\\
i. If $f\geq 0$ on~$E$ and $\int_E f(x)dx = 0$, then $f=0$ a.e. on~$E$.\\
ii. If $f$ is Lebesgue integrable over $E$, then $f$ is finite a.e. on~$E$.\\ \smallskip\\
{\bf Proof:} For each integer $n>0$ let $E_n = \{x\in E: f(x)>1/n \}$.\\
For each $n$, we note $mE_n = 0$ or else $\int_E f(x)dx >0$.\\
Let $A = \{x\in E: f(x) \not=0\}$. Then $A=\cup_{n=1}^{\infty} E_n$. Thus,
$mA = m(\cup_{n=1}^{\infty} E_n) \leq \sum_{n=1}^{\infty} m(E_n) = 0$
so that $f=0$~a.e. on~$E$, which proves~i.\\
In order to prove ii, for each integer $n>0$ let $E_n = \{x\in E: |f(x)|\geq n \}$.\\
Then $n\cdot mE_n \leq \int_{E_n} |f(x)|dx \leq \int_E |f(x)|dx=C$, so that
$mE_n \leq C/n$.\\
Let $A = \{x\in E: |f(x)| = \infty\}$. Then $A=\cap_{n=1}^{\infty} E_n$.
Since for each $n$, $A\subseteq E_n$, then $mA\leq mE_n\leq C/n$ so that $mA=0$,
which proves~ii.
\\ \smallskip\\
{\bf Proposition 2.35 (Lebesgue's criterion for Riemann integrability):}
Let $f$ be a bounded function on~$[a,b]$. Then $f$ is Riemann integrable over $[a,b]$
if and only if $f$ is continuous a.e. on~$[a,b]$. Proof in~\cite{rudin}. It involves
Proposition~2.33 and i of Proposition~2.34.\\ \smallskip\\
{\bf Observation 2.20:} Function $\chi_A$ in Observation 2.17 with $A$ equal to the set
of rational numbers fails the continuity hypothesis of Proposition~2.35 and thus it is
not Riemann integrable over~$[a,b]$ as observed there. Actually, it can be easily shown
to be nowhere continuous on~$[a,b]$.
\\ \smallskip\\
{\bf Proposition 2.36 (Countable additivity of the Lebesgue integral):} Let $\{E_n\}$ be
a countable collection of pairwise disjoint measurable sets.
Let $E=\cup_{i=1}^{\infty}\, E_n$, and let $f$ be a measurable function on~$E$.
Assume either $f\geq 0$ on~$E$ or $f$ is Lebesgue integrable over~$E$. Then
\[ \int_E f(x)dx = \sum_{i=1}^{\infty} \int_{E_i} f(x) dx. \]
Proof in~\cite{rudin}.\\ \smallskip\\
{\bf Observation 2.21:} If $f$ is a measurable function on a set $E$ with
\mbox{$mE=0$}, then $\int_E f(x)dx = 0$. Also, if sets $A$, $E$ are measurable with
$A\subseteq E$, and $f$ is Lebesgue integrable over $E$, then it is Lebesgue integrable
over~$A$. From all this then, if $f$ and $g$ are functions on a
measurable set $E$, $f=g$ a.e. on~$E$, $f$ Lebesgue integrable over $E$,
then so is $g$ and $\int_E g(x)dx = \int_E f(x)dx$.
Finally, we note that since integrals over sets of measure zero are zero, throughout
these notes, if sets $F$ and $E$ are measurable with $F\subseteq E$, $mF=0$, and
$f$ is a function on $E\setminus F$, possibly not defined on part or all of~$F$,
$f$~Lebesgue integrable over $E\setminus F$, we say $f$ is Lebesgue
integrable over~$E$ with $\int_E f(x)dx = \int_{E\setminus F} f(x)dx$. This makes
sense as it is always possible to define $f$ arbitrarily for points in~$F$ so that
then $f$ is defined on all of~$E$ and $\int_E f(x)dx = \int_{E\setminus F} f(x)dx +
\int_{F} f(x)dx = \int_{E\setminus F}f(x)dx + 0 = \int_{E\setminus F}f(x)dx$.
\\ \smallskip\\
{\bf Proposition 2.37 (Lebesgue's Monotone Convergence Theorem):} Let $\{f_n\}$
be an increasing sequence of nonnegative measurable functions on a
measurable set~$E$. Let $f$ be defined by $f(x) = \lim_{n\rightarrow\infty}f_n(x)$
for~$x\in E$. Then
\[ \int_E f(x)dx = \lim_{n\rightarrow\infty} \int_E f_n(x)dx. \]
Proof in~\cite{royden} and~\cite{rudin}. It involves Proposition~2.13.\\ \smallskip\\
{\bf Corollary 2.5:} Let $\{f_n\}$ be a sequence of nonnegative measurable functions
on a measurable set~$E$. Let $f$ be defined by $f(x) = \sum_{n=1}^{\infty}\, f_n(x)$
for~$x\in E$. Then 
\[ \int_E f(x)dx = \sum_{n=1}^{\infty} \int_E f_n(x)dx. \]
{\bf Proof:} $\{h_n\}$ defined by $h_n(x) =\sum_{i=1}^{n}\,f_i(x)$ for $x\in E$ is
an increasing sequence of nonnegative measurable functions on~$E$.\\ \smallskip\\
{\bf Observation 2.22:}  Proposition~2.36 can now be proved more easily.
It suffices to prove it for $f\geq 0$ on~$E$. Let $f_n(x) = f(x)\cdot \chi_{E_n}(x)$
for $x\in E$. Then $f(x)=\sum_{n=1}^{\infty}\,f_n(x)$ for $x\in E$ and the result
follows from Corollary~2.5.\\ \smallskip\\
{\bf Observation 2.23:} The following proposition says that if a nonnegative function
is Lebesgue integrable over a measurable set, then the Lebesgue integral of the
function over a measurable subset of the set is arbitrarily small if the measure of
the subset is small enough. Later we will see that it can be used
to show that every indefinite integral is absolutely continuous
(indefinite integrals and absolute continuity defined in the next section).\\ \smallskip\\
{\bf Proposition 2.38 (Absolute continuity of the Lebesgue integral):} Let $f$ be a
nonnegative Lebesgue integrable function over a measurable set~$E$. Then given
$\epsilon>0$ there is $\delta>0$ such that for each measurable set $A\subseteq E$
with $mA<\delta$, then~$\int_A f(x)dx <\epsilon$. Proof in~\cite{royden}.
It involves Lebesgue's Monotone Convergence Theorem (Proposition~2.37).
\\ \smallskip\\
{\bf Proposition 2.39 (Fatou's Lemma):} Let $\{f_n\}$ be a sequence of nonnegative
measurable functions on a measurable set~$E$. Let $f$ be defined by
\mbox{$f(x)= \liminf_{n\rightarrow\infty} f_n(x)$} for $x\in E$. Then
\[ \int_E f(x)dx \leq \liminf_{n\rightarrow\infty} \int_E f_n(x)dx. \]
Proof in~\cite{royden} and~\cite{rudin}. It involves Proposition~2.37.\\ \smallskip\\
{\bf Proposition 2.40 (Lebesgue's Dominated Convergence Theorem):} 
Let $\{f_n\}$ be a sequence of measurable functions on a
measurable set~$E$ such that there is a function $f$ on~$E$ with
$f_n\rightarrow f$ pointwise a.e. on~$E$. If
there is a function $g$ that is Lebesgue integrable over~$E$ such that
$|f_n|\leq g$ on~$E$ for all~$n$, then
\[ \int_E f(x)dx = \lim_{n\rightarrow\infty} \int_E f_n(x)dx. \]
Proof in~\cite{royden} and~\cite{rudin}. It involves Proposition~2.39.\\ \smallskip\\
{\bf Corollary 2.6 (Bounded Convergence Theorem):}
Let $\{f_n\}$ be a sequence of measurable functions on a
measurable set~$E$ of finite measure such that there is a function $f$ on~$E$
with $f_n\rightarrow f$ pointwise a.e. on~$E$. If
there is a real number $M$ such that $|f_n|\leq M$ on~$E$ for all~$n$, then
\[ \int_E f(x)dx = \lim_{n\rightarrow\infty} \int_E f_n(x)dx. \]
{\bf \large\ \\ The $L^p$ Spaces}\\ \smallskip\\
{\bf Definition 2.22:} Given a real number $p>0$, the $L^p[0,1]$ or $L^p$ {\bf space} is
the space of measurable functions on~$[0,1]$ satisfying: the $p$-th power of the
absolute value of each function in the space is Lebesgue integrable over~$[0,1]$.
Thus, a measurable function~$f$ on~$[0,1]$ is in $L^p$ (the $L^p$ space)
if and only if 
\[ \int_{[0,1]} |f(x)|^p dx < \infty. \]
Writing $\int_0^1 |f(x)|^p dx$ instead of  $\int_{[0,1]} |f(x)|^p dx$ for $f$
in~$L^p$, we define
\[ ||f||_p = \{ \int_0^1 |f(x)|^p dx \}^{1/p}  \]
and call $||\cdot||_p$ the $L^p$ {\bf norm},
and $||f||_p$ the $L^p$ {\bf norm} of~$f$.\\
Finally, the $L^{\infty}[0,1]$ or $L^{\infty}$ {\bf space} is the space of measurable
functions on $[0,1]$ satisfying: each function in the space is bounded on~$[0,1]$
except possibly on a set of measure zero. Thus, a measurable function $f$
on~$[0,1]$ is in~$L^{\infty}$ (the $L^{\infty}$ space) if and only if the essential
supremum of $|f|$ on~$[0,1]$ is finite, i.e.,
\[ \mathrm{ess}\ \sup |f(t)| = \inf \{M: m(\{t:|f(t)|>M\}) = 0\} <\infty. \]
We also note $\mathrm{ess}\ \sup |f(t)| =
\inf\ \{\sup_{t\in [0,1]} |g(t)|: g=f\ \mathrm{a.e.}\}$. Defining
\[ ||f||_{\infty} = \mathrm{ess}\ \sup |f(t)| \]
we call $||\cdot||_{\infty}$ the $L^{\infty}$ {\bf norm}, 
and $||f||_{\infty}$ the $L^{\infty}$ {\bf norm} of~$f$. \\ \smallskip\\
{\bf Observation 2.24:} In the definition of the $L^p$ spaces, the interval~$[0,1]$
was chosen for simplicity. Given a real number $p>0$, if $f\in L^p$, then clearly
$cf\in L^p$ for any real number~$c$. In addition, if $f,\ g\in L^p$, since
$|f+g|^p\leq 2^p(|f|^p + |g|^p)$, then $f+g\in L^p$. Thus, $L^p$ is a linear space
and so is $L^{\infty}$.\\
Given $f$ in $L^p$, $0<p\leq\infty$, then the $L^p$ norm of $f$, i.e., $||f||_p$
(Definition~2.22), equals zero if and only if~$f=0$~a.e. on~$[0,1]$. Accordingly,
we think of the elements of $L^p$ as equivalent classes of functions, each
class composed of functions that are equal to one another~a.e. on~$[0,1]$, and as noted in
Observation~2.21, some functions undefined on subsets of $[0,1]$ of measure zero.
Thus, assuming there is no distinction between two functions in the same equivalence
class, we note that given $p$, $1\leq p\leq \infty$, then the $L^p$ norm $||\cdot||_p$
is indeed a norm since clearly $||cf||_p=c||f||_p$ for any real number~$c$, and
as will be seen below, if $f,\ g\in L^p$, then $||f+g||_p\leq ||f||_p + ||g||_p$.
\\ \smallskip\\
{\bf Proposition 2.41 (H\"{o}lder's inequality):} Given $p, q$, $1\leq p,q \leq\infty$,
with \mbox{$1/p + 1/q = 1$}, if $f\in L^p$ and $g\in L^q$, then $f\cdot g \in L^1$ and
\[ \int_0^1 |(f\cdot g)(x)|dx \leq ||f||_p \cdot ||g||_q\,, \]
with equality for $p, q$, $1<p,q<\infty$ if and only if $\alpha|f|^p = \beta|g|^q$~a.e.
for nonzero constants $\alpha$ and~$\beta$.
Proof in \cite{royden}, \cite{rudin2}.  Proof in \cite{rudin} for $p=q=2$. \\ \smallskip\\
{\bf Proposition 2.42 (Minkowski's inequality):} Given $p$, $1\leq p\leq \infty$, if
$f, g\in L^p$, then $f+g\in L^p$ and
\[ ||f+g||_p \leq ||f||_p + ||g||_p. \]
Proof in \cite{royden}, \cite{rudin2}.  Proof in \cite{rudin} for $p=q=2$.\\ \smallskip\\
{\bf Observation 2.25:} For $p=q=2$, H\"{o}lder's inequality becomes {\bf Schwarz's
inequality}:
\[ \int_0^1 |(f\cdot g)(x)|dx \leq ||f||_2 \cdot ||g||_2 = \{\int_0^1 |f(x)|^2dx\}^{1/2}
\cdot \{\int_0^1 |g(x)|^2dx\}^{1/2}. \]
Note all of the above inequalities (H\"{o}lder's, Minkowski's, Schwarz's), in which
all integrations are over $[0,1]$, can be generalized by integrating everywhere over
a measurable set instead. Proof in~\cite{rudin2}. \\ \smallskip\\
{\bf Definition 2.23:} Given a norm $||\cdot||$ on a linear space $X$, we say $X$ is
a {\bf normed linear space with norm}~$||\cdot||$. We say this especially if among
all the possible norms that can be defined on~$X$, our current intent is to
associate $X$ exclusively with~$||\cdot||$.\\
A sequence $\{x_n\}$ in a normed linear space with norm~$||\cdot||$ is said to
{\bf converge in norm} to an element~$x$ in the space if, given~$\epsilon>0$,
there is an integer $N>0$ such that for $n\geq N$, then~$||x_n-x||<\epsilon$.\\
A sequence $\{x_n\}$ in a normed linear space with norm~$||\cdot||$ is said to
be a {\bf Cauchy sequence} if, given~$\epsilon$, there is an integer $N>0$ such
that for $n, m\geq N$, then~$||x_n-x_m||<\epsilon$.\\
A normed linear space with norm $||\cdot||$ is called {\bf complete} if every
Cauchy sequence in the space converges in norm to an element of the space.
\\ \smallskip\\
{\bf Proposition 2.43 (Riesz-Fischer):} Given $p$, $1\leq p\leq \infty$, then
$L^p$ is complete. Moreover, given $\{f_n\} \rightarrow f$ in $L^p$, then a
subsequence of $\{f_n\}$ converges pointwise to $f$ a.e. on~$[0,1]$.
Proof of first part in~\cite{royden}, \cite{royden2}. It involves Proposition~2.37
(Lebesgue's Monotone Convergence Theorem), Proposition~2.39 (Fatou's Lemma),
Proposition~2.40 (Lebesgue's Dominated Convergence Theorem) and ii of Proposition~2.34.
Proof of last part in~\cite{royden2}.
\\ \smallskip\\
{\bf Proposition 2.44 (Density of simple and step functions in $L^p$ space):}
Given $p$, $1\leq p\leq \infty$, then the subspace of simple functions on $[0,1]$
in $L^p$ is dense in $L^p$. Given $p$, $1\leq p<\infty$, then the subspace of step
functions on $[0,1]$ is dense in~$L^p$. Proof in~\cite{royden2}.
%\begin{figure}
%\centerline{\epsfig{figure=nlnk.ps,width=3.5in,height=3.0in}}
%\caption{A $3-$layer neural network. The leftmost layer is the input
%\end{figure}
%A copy of \mbox{NEURBT} can be obtained from
%\verb+http://math.nist.gov/~JBernal+
\section{\large Absolute Continuity and its Connections to\\ Lebesgue Integration}
{\bf Definition 3.1:} Let $f$ be a real-valued function defined on an interval~$[a,b]$.
Given $x\in [a,b]$, if for some finite number~$I$,
\[ I = \lim_{t\rightarrow x}\frac{f(t)-f(x)}{t-x},\ \ \ a<t<b,\ \ \ t\not= x, \]
then $f$ is said to be {\bf differentiable} at $x$; a number $f'(x)$ is defined
and said to exist by setting $f'(x)$ equal to~$I$; and $f'$ is said to exist at~$x$.
Accordingly, $f'$ is a function associated with
$f$, called the {\bf derivative} of~$f$, whose domain of definition is the set of
points $x$ at which $f'$ exists. If $f'$ exists at every point of
a set $E\subseteq [a,b]$, we say $f$ is differentiable on~$E$ or $f'$ exists on $E$.\\
Note that given $x\in [a,b]$, if the limit defining $I$ above equals $\infty$ or $-\infty$
then the convention here is to say that $f$ is not differentiable at~$x$.
\\ \smallskip \\
{\bf Proposition 3.1 (Fundamental Theorem of calculus I):}
Let $f$ be Riemann integrable over an interval~$[a,b]$. If there is a function
$F$ differentiable on~$[a,b]$ such that $F'=f$ on~$[a,b]$, then
\[ {\cal R} \int_a^b f(x)dx = F(b)-F(a). \]
Proof in~\cite{apostol} and~\cite{rudin}.\\ \smallskip\\
{\bf Proposition 3.2 (Fundamental Theorem of calculus II):}
Let $f$ be Riemann integrable over an interval~$[a,b]$. Define a function $F$ by
\[ F(x) = {\cal R} \int_a^x f(t)dt,\ \ \  x\in [a,b]. \]
Then $F$ is continuous on $[a,b]$, and if $f$ is continuous
at~$x\in [a,b]$, then $F$ is differentiable at~$x$ with~$F'(x)=f(x)$.
Proof in~\cite{apostol} and~\cite{rudin}.\\ \smallskip\\
{\bf Corollary 3.1 (Differentiability of the Riemann integral - Fundamental Theorem
of calculus for continuous functions):}\\
i. If $f$ is Riemann integrable over $[a,b]$ and $F(x) = {\cal R} \int_a^x f(t)dt$,
$x\in [a,b]$, then $F'=f$ a.e. on~$[a,b]$.\\
ii. If $f$ is continuous on~$[a,b]$, then there is a differentiable function~$F$ on~$[a,b]$
such that $F'=f$ on~$[a,b]$, and ${\cal R} \int_a^x f(t)dt = F(x)$. If $G$ is any
differentiable function on~$[a,b]$ such that $G'=f$ on~$[a,b]$, then
$G-F= C$, $C$ a constant, and ${\cal R} \int_a^x f(t)dt = G(x)- G(a)$, $G(a) = C$.
\\ \smallskip\\
{\bf Proof:} i follows from Proposition~3.2 and Proposition~2.35 (Lebesgue's criterion). First
part of ii from Proposition~3.2. Proof in~\cite{rudin} that $G'-F'=0$ on~$[a,b]$ implies $G-F=C$
on~$[a,b]$, $C$ a constant. ${\cal R} \int_a^x f(t)dt = G(x)- G(a)$ from Proposition~3.1.
Clearly $G(a)=C$ as ${\cal R} \int_a^x f(t)dt = F(x)$.  \\ \smallskip\\
{\bf Definition 3.2:} Let $f$ be a real-valued function defined on an interval~$[a,b]$.
Given $x\in [a,b]$, if for some finite number~$I$,
$I = \lim_{t\rightarrow x} f(t)$, $a\leq t\leq x$, then a number $f(x^-)$ called the
{\bf left-hand limit} of $f$ at~$x$ is defined by setting $f(x^-)$ equal to~$I$.
Similarly, if for some finite number~$I$, $I = \lim_{t\rightarrow x} f(t)$,
$x\leq t\leq b$, then a number $f(x^+)$ called the {\bf right-hand limit} of $f$
at~$x$ is defined by setting $f(x^+)$ equal to~$I$. \\ \smallskip\\
{\bf Observation 3.1:} A function $f$ is continuous at $x\in [a,b]$ if and only if
$f(x^-)$ and $f(x^+)$ exist and $f(x)=f(x^-)=f(x^+)$.
\\ \smallskip\\
{\bf Proposition 3.3 (\bf Monotonic functions: continuity):}
Let $f$ be a monotonic real-valued function on an interval~$[a,b]$.
Then $f(x^-)$ and $f(x^+)$ exist for every point $x\in [a,b]$, and
the set of points of $[a,b]$ at which $f$ is discontinuous is at most countable.
Proof in~\cite{rudin}.\\ \smallskip \\
{\bf Corollary 3.2 (Monotonic surjective $f$ implies $f$ is continuous):}
If $f$ is monotonic from $[a,b]$ onto $[c,d]$, then $f$ is continuous on~$[a,b]$.
\\ \smallskip\\
{\bf Proof:} Assume $f$ is discontinuous at $x\in [a,b]$. Since $f(x^-)$ and $f(x^+)$
exist from Proposition~3.3, it must be that $f(x^-)\not=f(x^+)$ so that $y$ exists in
$[c,d]$ between $f(x^-)$ and $f(x^+)$, $y\not= f(x)$. But then $y$ can not be in the range
of~$f$ as $f$ is monotonic, which contradicts that the range of $f$ is all of~$[c,d]$.
\\ \smallskip\\
{\bf Observation 3.2:} A function $f$ is described below from $[0,1]$ into~$[0,1]$ that is
strictly increasing on $[0,1]$, discontinuous at each rational number in~$(0,1]$,
continuous at each irrational number in~$[0,1]$ and at~zero, $f(0)=0$, $f(1)=1$.
\\ \smallskip\\
Let $\{r_n\}_{n=1}^{\infty}$ be an enumeration of the rational numbers in~$(0,1]$.\\
Given $x\in (0,1]$, let $R(x)=\{n: r_n\leq x\}$, and set $R(0)=\emptyset$.\\
Define $f:[0,1]\rightarrow [0,1]$ by $f(0)=0$ and
\[ f(x)=\sum_{n\in R(x)} 1/2^n,\ \ \ x\in (0,1]. \]
Given $x,x' \in [0,1]$, $x<x'$, then $R(x)\subseteq R(x')$, and since there is a
rational number~$r$ such that $x<r<x'$, then $R(x)\not= R(x')$. Thus, it must be
that $f(x)<f(x')$ so that $f$ is strictly increasing on~$[0,1]$.\\ \smallskip\\
Since $R(1)$ includes every $n$ then $f(1)=1$.\\ \smallskip\\
Let $x$ be a rational number in $(0,1]$. We show $f$ is discontinuous at~$x$.\\
For some integer $k>0$, $x= r_k$. Thus, $k\in R(x)$ but $k\not\in R(x')$ for every
$x'\in [0,1]$, $x'<x$. $R(x')\subseteq R(x)$ then implies $f(x)-f(x')>1/2^k$.\\
Thus, $f$ is discontinuous at $x$ (a rational number in~$(0,1]$).\\ \smallskip\\
Let $x$ be an irrational number in $[0,1]$. We show $f$ is continuous at~$x$.\\
Given $\epsilon>0$, choose integer $N>0$ such that $1/2^N <\epsilon$, and let
\[ \delta=\min_{n\leq N}|x-r_n|. \]
Given $x'\in [0,1]$, $x'<x$, $|x-x'|<\delta$, then $R(x')\subseteq R(x)$, and if
$n\in R(x)\setminus R(x')$, it must be that $x'<r_n<x$ so that $|x-r_n|<\delta$ and
thus~$n>N$. Accordingly, $f(x)-f(x')\leq \sum_{n=N+1}^{\infty}1/2^n =1/2^N<\epsilon$.\\
Finally, given $x'\in [0,1]$, $x'>x$, $|x'-x|<\delta$, then $R(x)\subseteq R(x')$,
and if $n\in R(x')\setminus R(x)$, it must be that $x<r_n\leq x'$ so that
$|x-r_n|<\delta$ and thus~$n>N$. Accordingly,
$f(x')-f(x)\leq \sum_{n=N+1}^{\infty}1/2^n =1/2^N<\epsilon$.\\
Thus, $f$ is continuous at $x$ (an irrational number in~$[0,1]$) and at
zero by an argument similar to the one just used for the case~$x'>x$.
\\ \smallskip\\
{\bf Proposition 3.4 (\bf Monotonic functions: differentiability):}
Let $f$ be a monotonic real-valued function on an interval~$[a,b]$.
Then $f$ is differentiable a.e. on~$[a,b]$, and $f'$ is measurable. If, in addition,
$f$ is increasing on~$[a,b]$ (note $f'\geq 0$ where it exists), then $f'$ is Lebesgue
integrable over~$[a,b]$, and \[ \int_a^b f'(x) dx \leq f(b)- f(a), \]
where we write $\int_a^b f'(x)dx$ instead of $\int_{[a,b]} f'(x)dx$.
Proof in~\cite{bruckner} and \cite{royden}. It involves Proposition~2.39 (Fatou's Lemma)
and ii of Proposition~2.34.
\\ \smallskip\\
{\bf Definition 3.3:} Let $f$ be a real-valued function defined on an interval~$[a,b]$.
Given a partition $P= \{x_0, x_1,\ldots,x_n\}$ of $[a,b]$,
set $\Delta f_i = f(x_i)-f(x_{i-1})$, $i=1,\ldots,n$, and define
\[ V(f;a,b) = \sup_P \sum_{i=1}^n |\Delta f_i|, \]
the supremum taken over all partitions $P$ of $[a,b]$.\\
$f$ is said to be of {\bf bounded variation} on~$[a,b]$ if~$V(f;a,b)<\infty$.
\\ \smallskip\\
{\bf Proposition 3.5 (Jordan decomposition):} A function $f$ is of bounded variation
on $[a,b]$ if and only if it is the difference of two monotonically increasing
real-valued functions on~$[a,b]$. Proof in \cite{royden} and~\cite{rudin}.\\ \smallskip\\
{\bf Corollary 3.3:} If $f$ is of bounded variation on $[a,b]$ then $f$ is
differentiable a.e. on~$[a,b]$, and $f'$ is measurable and Lebesgue integrable
over~$[a,b]$.\\ \smallskip\\
{\bf Proof:} By Proposition~3.5, $f = f_1 - f_2$ on~$[a,b]$, where $f_1$ and $f_2$
are monotonically increasing on~$[a,b]$. Thus, by Proposition~3.4,
$f'$ is measurable and exists a.e. on~$[a,b]$. Since $|f'|\leq |f_1'| + |f_2'| =
f_1' + f_2'$ a.e. on~$[a,b]$, then again by Proposition~3.4,
\[ \int_a^b |f'(x)|dx \leq \int_a^b f_1'(x)dx+\int_a^b f_2'(x)dx
\leq f_1(b)-f_1(a) + f_2(b)-f_2(a), \]
and therefore $f'$ is Lesbegue integrable over $[a,b]$ (Proposition~2.31).
\\ \smallskip\\
{\bf Definition 3.4:} Given a Lebesgue integrable function $f$ over $[a,b]$, and a
real-valued function $F$ on~$[a,b]$ such that
\[ F(x) = F(a) + \int_a^x f(t)dt,\ \ \ x\in [a,b],\]
then the function $F$ is said to be an {\bf indefinite integral} of~$f$ over~$[a,b]$.
\\ \smallskip\\
{\bf Proposition 3.6 (Indefinite integral of $f$ zero everywhere, then $f$ is zero a.e.):}
If $f$ is Lebesgue integrable over $[a,b]$ and $\int_a^x f(t)dt = 0$ for all $x\in [a,b]$,
then $f = 0$ a.e. on $[a,b]$.  Proof in~\cite{royden}. It involves Proposition~2.15.
\\ \smallskip\\
{\bf Proposition 3.7 (Differentiability of the indefinite integral):}
Let $f$ be Lebesgue integrable over an interval~$[a,b]$, and $F$ a function such that
\[ F(x) = F(a) + \int_a^x f(t)dt,\ \ \  x\in [a,b], \]
i.e., an indefinite integral. Then $F'=f$ a.e. on~$[a,b]$.\\
Proof in~\cite{royden}. It involves Proposition~3.6, Corollary~2.6 (Bounded Convergence
Theorem), the inequality in Proposition~3.4, and i of Proposition~2.34.
\\ \smallskip\\
{\bf Definition 3.5:} A real-valued function $f$ defined on an interval $[a,b]$ is said
to be {\bf absolutely continuous} on $[a,b]$ if for every $\epsilon>0$ there is
$\delta>0$ such that
\[ \sum_{i=1}^n |f(x_i') - f(x_i)| < \epsilon \]
for any integer $n>0$ and any disjoint collection of open intervals
$(x_i,x_i') \subseteq [a,b]$, $i=1,\ldots,n$, with
\[ \sum_{i=1}^n (x_i' - x_i) < \delta. \]
{\bf Proposition 3.8 (Absolutely continuous $f$ is constant if $f'$ is zero a.e.):}
If $f$ is absolutely continuous on $[a,b]$ with $f'=0$ a.e. on~$[a,b]$, then
$f$ is constant on~$[a,b]$, i.e., $f(x)=f(a)$ for all~$x\in [a,b]$.
Proof in~\cite{royden}.\\ \smallskip\\
{\bf Observation 3.3:} Absolutely continuous $\Rightarrow$ {\bf uniformly continuous}~\cite{rudin}
$\Rightarrow$ continuous. Moreover, a continuous real-valued function of compact domain is
uniformly continuous~\cite{rudin}. Accordingly, a function $f$ called the {\bf Cantor function}
from $[0,1]$ onto $[0,1]$ that is continuous, thus uniformly continuous, but not absolutely
continuous is described below. This function $f$ is monotonically increasing on~$[0,1]$
and thus differentiable a.e. on~$[0,1]$. Actually, $f'=0$ at points not in the Cantor set
(described in Observation~2.4) and does not exist at points in it. Thus, $f'=0$ a.e.
on~$[0,1]$, $f$ is not constant on~$[0,1]$, hence $f$ can not be absolutely continuous
on~$[0,1]$ by Proposition~3.8. \\ \smallskip\\
For this purpose, we note that given $x\in [0,1]$, $x$ can be expressed in its ternary
expansion as $0.a_1a_2a_3\,$$\cdot\cdot\cdot$ so that $x=\sum_{n=1}^{\infty} a_n/3^n$,
$a_n\in \{0,1,2\}$. Note $x=1$ is then expressed as~$0.222\,$$\cdot\cdot\cdot$.
Similarly, given $y\in [0,1]$, $y$ can be expressed in its binary expansion as
$0.b_1b_2b_3\,$$\cdot\cdot\cdot$ so that $y=\sum_{n=1}^{\infty} b_n/2^n$,
$b_n\in \{0,1\}$. Note $y=1$ is then expressed as~$0.111\,$$\cdot\cdot\cdot$.\\ \smallskip\\
In Observation 2.4 the Cantor set was identified as $\cap_{n=1}^{\infty} E_n$,
where $E_1$ is the union of $[0,1/3]$ and $[2/3,1]$ obtained by removing
the open middle third of $[0,1]$, $E_2$ is the union of $[0,1/9]$, $[2/9,3/9]$, $[6/9,7/9]$,
$[8/9,1]$ obtained by removing the open middle thirds of $[0,1/3]$ and $[2/3,1]$, and so on.
Actually, with $E_0=[0,1]$, then at stage~$m$, open intervals of the form
$((3k-2)/3^m,(3k-1)/3^m)$, $k\in \{1,\ldots,3^{m-1}\}$, are removed from $E_{m-1}$, if
contained in it, to obtain~$E_m$. We note that endpoints of any such intervals have two
ternary expansions, and in what follows, only the expansion of any such point that contains
no 1's is considered. Fixing one of these removed open intervals, we note it is the open middle
third of a closed interval in $E_{m-1}$, all numbers in the closed interval in $E_{m-1}$ having
the same first $m-1$ digits in their ternary expansions, none of them equal to~1. Finally, we
note numbers in the removed open interval have 1 as the $m^{th}$ digit of their ternary
expansions, while numbers in the closed left and right thirds of the closed interval,
closed thirds that become part of $E_m$, have 0 and 2, respectively, as the $m^{th}$
digit of their ternary expansions. Thus, the Cantor set is exactly the set of numbers
in $[0,1]$ that have no 1's in their ternary expansions.
\\ \smallskip\\
An attempt can be made to identify the Cantor function as follows. Recalling that
$(1/3,2/3)$ was the open middle third that was removed from $[0,1]$ to obtain~$E_1$,
given $x$ in its closure, i.e., in $[1/3,2/3]$, set $f(x) = 1/2$. Again, recalling that
$(1/9,2/9)$ and $(7/9,8/9)$ were the open middle thirds that were removed from
$[0,1/3]$ and $[2/3,1]$, respectively, to obtain~$E_2$, given $x$ in the closure of
$(1/9,2/9)$, i.e., in $[1/9,2/9]$, set $f(x)=1/4$, and given $x$ in the closure of
$(7/9,8/9)$, i.e., in $[7/9,8/9]$, set $f(x)=3/4$. Accordingly, $f$ can be identified this
way at each stage of the contruction of the Cantor set but this is not enough as it
has not been identified for points in ``the limit" that are part of the Cantor~set.
\\ \smallskip \\
The Cantor function is properly identified as follows. Given $x\in [0,1]$ with ternary
expansion $0.a_1a_2a_3\,$$\cdot\cdot\cdot$ so that $x=\sum_{n=1}^{\infty} a_n/3^n$,
$a_n\in \{0,1,2\}$, let $N$ be the smallest $n$ such that $a_n$ equals~1. If such
an $n$ does not exist, i.e., $x$~is in the Cantor set, let~$N=\infty$.
With $b_n=a_n/2$ if $n<N$, $b_n=1$ if $n=N$, and $b_n=0$ if $n>N$, let $y$ be the
number in $[0,1]$ with binary expansion $0.b_1b_2b_3\,$$\cdot\cdot\cdot$
so that $y=\sum_{n=1}^{\infty} b_n/2^n = \sum_{n=1}^N b_n/2^n$, and set~$f(x)=y$.
The function $f$ identified this way is then called the Cantor function.\\ \smallskip\\
{\bf Proposition 3.9:} Let $f$ be the Cantor function. Then $f$ is continuous, thus
uniformly continuous, from $[0,1]$ onto~$[0,1]$. In addition, $f$ is monotonically
increasing on~$[0,1]$ and thus differentiable a.e. on~$[0,1]$. Actually, $f'=0$ at
points not in the Cantor set and does not exist at points in it.
\\ \smallskip\\
{\bf Proof:} Given $x_1$, $x_2 \in [0,1]$, $x_1<x_2$, we show $f(x_1)\leq f(x_2)$.\\
Let $0.a_1a_2a_3$$\cdot\cdot\cdot$, $0.c_1c_2c_3$$\cdot\cdot\cdot$ be $x_1$, $x_2$,
respectively, in their ternary expansions. 
Let $0.b_1b_2b_3$$\cdot\cdot\cdot$, $0.d_1d_2d_3$$\cdot\cdot\cdot$ be $f(x_1)$,
$f(x_2)$, respectively, in their binary expansions.\\
Let $N_1$ be the smallest $n$ such that $a_n=1$; $N_1=\infty$ if there is no such~$n$.\\
Let $N_2$ be the smallest $n$ such that $c_n=1$; $N_2=\infty$ if there is no such~$n$.\\
Let $N'$ be the smallest $n$ such that $a_n<c_n$.\\
If $N'>N_1$, since then $c_n=a_n$, $n=1,\ldots,N_1$, and, in particular,
$c_{N_1}=a_{N_1}=1$, it must be that $N_2=N_1$ so that $b_n=d_n$, $n=1,\ldots,N_1=N_2$,
and therefore~$f(x_1)=\sum_{n=1}^{N_1}b_n/2^n=\sum_{n=1}^{N_2}d_n/2^n=f(x_2)$.\\
Similarly if $N'>N_2$, and the case $N'=N_1=N_2$ can not be.\\
If $N'=N_1$ and $N_2>N'$, since $a_{N_1}=1$, it must be that $c_{N_1}=2$ so that
$b_n=d_n$, $n=1,\ldots,N_1-1$, $b_{N_1}=d_{N_1}=1$. Therefore,
$\sum_{n=1}^{N_1}b_n/2^n=\sum_{n=1}^{N_1}d_n/2^n$ thus
$f(x_1)=\sum_{n=1}^{N_1}b_n/2^n\leq \sum_{n=1}^{N_2}d_n/2^n=f(x_2)$.\\
If $N'=N_2$ and $N_1>N'$, since $c_{N_2}=1$, it must be that $a_{N_2}=0$ so that
$b_n=d_n$, $n=1,\ldots,N_2-1$, $b_{N_2}=0$, $d_{N_2}=1$. Therefore,
$\sum_{n=1}^{N_2}b_n/2^n<\sum_{n=1}^{N_2}d_n/2^n$ thus
$f(x_1)=\sum_{n=1}^{N_1}b_n/2^n\leq \sum_{n=1}^{N_2}d_n/2^n=f(x_2)$.\\
Finally, if $N_1>N'$, $N_2>N'$, since $a_{N'}<c_{N'}$, it must be that $a_{N'}=0$,
$c_{N'}=2$, so that
$b_n=d_n$, $n=1,\ldots,N'-1$, $b_{N'}=0$, $d_{N'}=1$. Therefore,
$\sum_{n=1}^{N'}b_n/2^n<\sum_{n=1}^{N'}d_n/2^n$ thus
$f(x_1)=\sum_{n=1}^{N_1}b_n/2^n\leq \sum_{n=1}^{N_2}d_n/2^n=f(x_2)$.
Thus, $f(x_1)\leq f(x_2)$ for all cases and therefore $f$ is monotonically increasing.
\\ \smallskip\\
Given $y\in [0,1]$, we show there is $x\in [0,1]$ with~$f(x)=y$.\\
Let $0.b_1b_2b_3$$\cdot\cdot\cdot$ be $y$ in its binary expansion.\\
For each $n$, let $a_n=2b_n$. Then for each $n$, $a_n$ is either zero or two.\\
Let $x$ be the point in $[0,1]$ which in its ternary expansion is
$0.a_1a_2a_3$$\cdot\cdot\cdot$.\\
Then $x$ is actually a point in the Cantor set and~$f(x)=y$.\\
Thus, $f$ is onto~$[0,1]$.\\ \smallskip\\
That $f$ is continuous, thus uniformly continuous, on $[0,1]$,
now follows from Corollary~3.2.
\\ \smallskip\\
Finally, given $x\in [0,1]$, if $x$ is not in the Cantor set, we show $f'(x)=0$. On
the other hand, if $x$ is in the Cantor set, we show that $f'(x)$ does not exist.\\
If $x$ is not in the Cantor set, its ternary expansion must contain 1 as one of
its digits. Then for some integer $m$, $m>0$, the $m^{th}$ digit of the expansion
equals~1 with no previous digits equal to~1. It follows that $x$ must be contained
in an open interval of the form $((3k-2)/3^m,(3k-1)/3^m)$, $k\in \{1,\ldots,3^{m-1}\}$.
Thus, it suffices to show $f$ is constant on any such interval. But this follows
immediately since all numbers in the interval have the same first $m$ digits in their
ternary expansions with 1 as the $mth$ digit and no previous digits equal to~1.\\
On the other hand, if $x$ is in the Cantor set, its ternary expansion consists of
0's and 2's. Given an integer $n>0$, define $x_n$ to be the number in~$[0,1]$
whose ternary expansion is exactly that of $x$ except at its $n^{th}$ digit. Its
$n^{th}$ digit is 0 if the $n^{th}$ digit of $x$ is 2, and it is 2 if that of $x$
is~0. It follows then that $|x_n-x|=2/3^n$ so that $x_n\rightarrow x$.
Also, $|f(x_n)-f(x)|=1/2^n$.  Thus, since $f$ is monotonically increasing,
$(f(x_n)-f(x))/(x_n-x)=|f(x_n)-f(x)|/|x_n-x|=(1/2)(3/2)^n\rightarrow\infty$
so that $f'(x)$ does not~exist.
\\Thus, $f'$ does not exist at points in the Cantor set and equals zero otherwise.
\\ \smallskip \\
{\bf Corollary 3.4:} The Cantor function is not absolutely continuous on~$[0,1]$.
\\ \smallskip\\
{\bf Proof:} Let $f$ be the Cantor function and assume it is absolutely continuous
on~$[0,1]$. By Proposition~3.9, $f'=0$ a.e.  on $[0,1]$. Thus, by Proposition~3.8,
$f$ must be constant on~$[0,1]$, i.e., $f(x)=f(0)$ for all $x\in [0,1]$. But this
is a contradiction as for instance $f(0)=0$ and~$f(1)=1$.
Thus, $f$ is not absolutely continuous on~$[0,1]$. \\ \smallskip\\
{\bf Observation 3.4:} For the sake of completeness, we analyze the nondifferentiability
of the Cantor function $f$ on the Cantor set.\\
Let $x_L$ and $x_R$ be points in the Cantor set that are the left and right endpoints
of an open interval $I$ removed at the $m^{th}$ stage of the construction of the Cantor~set.
It must then be that in their ternary expansions,~$x_L$ can be expressed as
$0.a_1a_2$$\cdot\cdot\cdot$$a_{m-1}0\overline{2}$ (a~bar on a digit means the digit is
infinitely repeated), and $x_R$ as $0.a_1a_2$$\cdot\cdot\cdot$$a_{m-1}2\overline{0}$,
the set $\{a_1,a_2,\ldots,a_{m-1}\}$ with elements equal to $0$ or $2$ if~$m>1$,
empty if~$m=1$. Given $x$ in the open interval $I$, it must be that in its ternary
expansion the $m^{th}$ digit is~1, and if $m>1$, then the first $m-1$ digits are also
$a_1,a_2,\ldots,a_{m-1}$. Define for each integer $n$, $n>0$, a number
$x_L^n$ in $I$ that in its ternary expansion the first $m$ digits are as described
above, and all other digits are~$0$ except the $(m+n)^{th}$ digit which is~$1$. Then
$\lim_{n\rightarrow\infty} x_L^n = x_L$ and $f(x_L^n)=f(x_L)$ for all~$n$ so that
$\lim_{n\rightarrow\infty} (f(x_L^n)-f(x_L))/(x_L^n-x_L)=0$.
Since for any sequence $\{x_n\}$ in $I$, with
$\lim_{n\rightarrow\infty} x_n = x_L$, then $f(x_n)=f(x_L)$ for all~$n$, it follows
that $(f(t)-f(x_L))/(t-x_L)$ has a limit as $t\rightarrow x_L$ from the right side
of~$x_L$ and it is~zero. Similarly, define for each integer $n$, $n>0$, a number
$x_R^n$ in $I$ that in its ternary expansion the first $m$ digits are as described
above, and all other digits are~$0$ except the $(m+1)^{th},\ldots,(m+n)^{th}$
digits which are~$2$. Then
$\lim_{n\rightarrow\infty} x_R^n = x_R$ and $f(x_R^n)=f(x_R)$ for all~$n$ so that
$\lim_{n\rightarrow\infty} (f(x_R^n)-f(x_R))/(x_R^n-x_R)=0$.
Since for any sequence $\{x_n\}$ in $I$, with $\lim_{n\rightarrow\infty} x_n = x_R$,
then $f(x_n)=f(x_R)$ for all~$n$, it follows that $(f(t)-f(x_R))/(t-x_R)$ has a limit
as $t\rightarrow x_R$ from the left side of~$x_R$ and it is~zero.\\
In the proof of Proposition 3.9, given any $x$ in the Cantor set, a sequence $\{x_n\}$
of points in the Cantor set was identifed with $x_n\rightarrow x$ and
$(f(x_n)-f(x))/(x_n-x)\rightarrow\infty$. We show that with $x_L$, $x_R$ as above,
then $(f(t)-f(x_L))/(t-x_L)$ has a limit as $t\rightarrow x_L$ from the left side
of~$x_L$ and it is~$\infty$, and $(f(t)-f(x_R))/(t-x_R)$ has a limit as
$t\rightarrow x_R$ from the right side
of~$x_R$ and it is also~$\infty$. Actually, we only show it for $x_R$ as the proof
for $x_L$ can be similarly accomplished. Accordingly, let $n\geq m$ be an integer
such that the ternary expansions of~$t$ and~$x_R$ coincide in the first $n$ digits and
the $(n+1)^{th}$ digit of~$t$ is~1 or~2. As mentioned above, all digits of~$x_R$ after
the $m^{th}$ digit equal~0. Thus, $f(t)-f(x_R)\geq 1/2^{n+1}$ and
$t-x_R \leq 2/3^{n+1} +2/3^{n+2} +\,\cdot\cdot\cdot=1/3^n$, so that
\[ \lim_{t\rightarrow x_R^+} \frac{f(t)-f(x_R)}{t-x_R} \geq
\lim_{n\rightarrow\infty}(1/2)(3/2)^n=\infty. \]
Finally, it is of interest to note that if $x$ is any point in the Cantor set,
then at stage~$m$ of the contruction of the Cantor set, $x$ is in a closed interval
$[a_m,b_m]\subset [0,1]$, where if $0.x_1x_2$$\cdot\cdot\cdot$ is~$x$ in its ternary
expansion, then $0.x_1\cdot\cdot\cdot x_m\overline{0}$ is $a_m$ in its ternary
expansion, and $0.x_1\cdot\cdot\cdot x_m\overline{2}$ is $b_m$ in its ternary
expansion. It follows that $b_m-a_m = \sum_{i=m+1}^{\infty}2/3^i=1/3^m$ and
$f(b_m)-f(a_m) = \sum_{i=m+1}^{\infty}1/2^i= 1/2^m$. Thus, with $a_m\leq x\leq b_m$,
we have \[ \lim_{m\rightarrow\infty} \frac{f(b_m)-f(a_m)}{b_m-a_m}=
\lim_{m\rightarrow\infty}(3/2)^m=\infty. \]
If $x= x_L$, $x_L$ as above, then for some $m$, $x=b_m$ and $\lim_{m\rightarrow\infty}
(f(x)-f(a_m))/(x-a_m)=\infty$, as expected. Similarly,
if $x= x_R$, $x_R$ as above, then for some $m$, $x=a_m$ and $\lim_{m\rightarrow\infty}
(f(x)-f(b_m))/(x-b_m)=\infty$, also as expected.
As for a point $x$ in the Cantor set that is not an endpoint of an open
interval removed at some stage of the construction of the Cantor set, it is easier
to see that $\lim_{m\rightarrow\infty} (f(x)-f(a_m))/(x-a_m)=\infty$,
and $\lim_{m\rightarrow\infty} (f(x)-f(b_m))/(x-b_m)=\infty$, by looking at the
ternary expansions of $x$, $a_m$ and $b_m$. Actually, we only show it for $\{a_m\}$
as the proof for $\{b_m\}$ can be similarly accomplished.  Accordingly, let $m>0$
be an integer such that the $(m+1)^{th}$ digit of $x$ in its ternary expansion,
i.e., $x_{m+1}$, equals~2. As mentioned above, the ternary expansions of $x$ and
$a_m$ coincide in the first $m$ digits and all digits of $a_m$ after the
$m^{th}$ digit equal~0. Thus, $f(x)-f(a_m)\geq 1/2^{m+1}$ and
$x-a_m \leq 2/3^{m+1} +2/3^{m+2} +\,\cdot\cdot\cdot=1/3^m$, so that
\[ \lim_{m\rightarrow\infty} \frac{f(x)-f(a_m)}{x-a_m} \geq
\lim_{n\rightarrow\infty}(1/2)(3/2)^m=\infty. \]
Note that $\lim_{m\rightarrow\infty} (f(x)-f(a_m))/(x-a_m)=$
$\lim_{m\rightarrow\infty} (f(x)-f(b_m))/(x-b_m)=\infty$
does not imply that $\lim_{t\rightarrow x} (f(x)-f(t))/(x-t)=\infty$.
\\ \smallskip\\
{\bf Observation 3.5:} A function on $[a,b]$ that is a finite linear combination of
absolutely continuous functions on~$[a,b]$ is absolutely continuous on~$[a,b]$.
The proof is analogous to the proof that a finite linear combination of continuous
functions is continuous. In addition, the product of two absolutely continuous
functions on $[a,b]$ is absolutely continuous on~$[a,b]$.\\ \smallskip\\
{\bf Proposition 3.10 (Absolutely continuous $f$ implies $f$ is of bounded variation):}
If $f$ is absolutely continuous on $[a,b]$, then $f$ is of bounded
variation on~$[a,b]$. Proof in~\cite{royden}. \\ \smallskip\\
{\bf Corollary 3.5:} If $f$ is absolutely continuous on $[a,b]$ then $f$ is
differentiable a.e. on~$[a,b]$, and $f'$ is measurable and Lebesgue integrable
over~$[a,b]$.\\ \smallskip\\
{\bf Proposition 3.11 (Absolute continuity of the indefinite integral):} If $F$
is an indefinite integral over~$[a,b]$, then $F$ is absolutely continuous on~$[a,b]$.
\\ \smallskip\\
{\bf Proof:} Assume (Definition 3.4) $F(x) = F(a) + \int_a^x f(t)dt$, $x\in [a,b]$, $f$
is Lebesgue integrable over~$[a,b]$. By Proposition~2.31, $|f|$ is Lebesgue integrable
over~$[a,b]$.
Then by Proposition~2.38, given $\epsilon>0$ there is $\delta>0$ such that for each
measurable set $A\subseteq [a,b]$ with $m(A)<\delta$, then~$\int_A |f(t)|dt<\epsilon$.\\
Given integer $n>0$ and disjoint open intervals
$(x_i,x_i') \subseteq [a,b]$, $i=1,\ldots,n$, with
$\sum_{i=1}^n (x_i' - x_i) < \delta$, let $A=\cup_{i=1}^n(x_i,x_i')$. Then $A$ is
measurable and $m(A)<\delta$. Thus,~$\int_A |f(t)|dt<\epsilon$. Accordingly, then
\begin{eqnarray*}
\sum_{i=1}^n |F(x_i') - F(x_i)| &=& \sum_{i=1}^n |\int_a^{x_i'}f(t)dt - \int_a^{x_i}f(t)dt|
= \sum_{i=1}^n |\int_{x_i}^{x_i'}f(t)dt|\\
&\leq& \sum_{i=1}^n \int_{x_i}^{x_i'}|f(t)|dt = \int_A |f(t)|dt<\epsilon.
\end{eqnarray*}
Thus, $F$ is absolutely continuous on $[a,b]$.\\ \smallskip\\
{\bf Proposition 3.12 (Equivalent conditions for an absolutely continuous function):}
Given a real-valued function $f$ on $[a,b]$,
then the following three conditions are equivalent:\\
i. $f$ is absolutely continuous on $[a,b]$.\\
ii. There exists a Lebesgue integrable function~$g$ over~$[a,b]$ such that\\
$f(x) = f(a) + \int_a^x g(t)dt$, $x\in [a,b]$.\\
(Note that then by Proposition~3.7, $f'=g$ a.e. on~$[a,b]$).\\
iii. $f'$ exists a.e. on $[a,b]$ and is Lebesgue integrable over~$[a,b]$, and\\
$f(x) = f(a) + \int_a^x f'(t)dt$, $x\in [a,b]$.\\ \smallskip\\
{\bf Proof:}\\
iii $\Rightarrow$ ii:\\ Take $g=f'$.\\
ii $\Rightarrow$ i:\\ This is Proposition 3.11.\\
i $\Rightarrow$ iii:\\ By Corollary 3.5, $f'$ exists a.e. on $[a,b]$ and is Lebesgue
integrable over~$[a,b]$.\\
Let $G(x) = \int_a^x f'(t)dt$, $x\in [a,b]$. Then $G$ is an indefinite integral
of~$f'$ over~$[a,b]$ and by Proposition~3.11, is absolutely continuous on~$[a,b]$, and so
is the function $h = f-G$ by Observation~3.5.\\
By Proposition 3.7, $G'=f'$ a.e. on~$[a,b]$. Thus $h'=0$ a.e. on~$[a,b]$, and by
Proposition~3.8, $h$ is constant on~$[a,b]$, i.e., $f-G=C$ on~$[a,b]$ for some
constant~$C$, i.e., $f(x)-\int_a^x f'(t)dt=C$, $x\in [a,b]$.\\
Since $C=f(a)$, it then follows that $f(x) = f(a) + \int_a^x f'(t)dt$, $x\in [a,b]$.
\\ \smallskip\\
{\bf Corollary 3.6 (Fundamental Theorem of Lebesgue integral calculus):}
Given real-valued functions $f$, $g$ on~$[a,b]$, $f$ absolutely continuous
on~$[a,b]$ and $f'=g$ a.e. on~$[a,b]$, then $f(x) = f(a) + \int_a^x g(t)dt$, $x\in [a,b]$.
\\ \smallskip\\
{\bf Proposition 3.13 (Fundamental Theorem of Lebesgue integral calculus
(Alternate form)):} Given a real-valued function $f$ on $[a,b]$, if $f'(x)$ exists
for every $x\in [a,b]$, and $f'$ is Lebesgue integrable over~$[a,b]$, then
$f(x) = f(a) + \int_a^x f'(t)dt$, $x\in [a,b]$. Proof in~\cite{rudin2}.\\ \smallskip\\
{\bf Corollary 3.7:} Given a real-valued function $f$ on $[a,b]$, if $f'$ exists everywhere
on~$[a,b]$ and $f'$ is Lebesgue integrable over~$[a,b]$, then $f$ is absolutely continuous
on~$[a,b]$.\\ \smallskip\\
{\bf Proof:} From Proposition 3.13 and then Proposition 3.12.\\ \smallskip\\
{\bf Proposition 3.14 (Change of variable for Riemann integral):} Given a strictly
monotonic continuous function $u$ from an interval~$[a,b]$ onto an interval~$[c,d]$
($u(a)=c$, $u(b)=d$ if $u$ is strictly increasing, $u(a)=d$, $u(b)=c$
if it is strictly decreasing), with $u'$ Riemann integrable over~$[a,b]$,
and a real-valued Riemann integrable function $f$ over~$[c,d]$, then
\[ {\cal R} \int_{u(a)}^{u(b)} f(x)dx = {\cal R} \int_a^b f(u(t))u'(t)dt, \]
with ${\cal R} \int_{u(a)}^{u(b)} f(x)dx = -{\cal R} \int_c^d f(x)dx$
if $u$ is decreasing.
Proof in \cite{apostol} and \cite{rudin}. Note that by Proposition 3.4, $u'$ exists
a.e. on~$[a,b]$.\\ \smallskip\\
{\bf Proposition 3.15 (Substitution rule for Riemann integral):} Given a function $u$
from an interval $[a,b]$ into an interval $I$ such that $u'(x)$ exists for every
$x\in [a,b]$ with $u'$ Riemann integrable over $[a,b]$, and a real-valued continuous
function $f$ on $I$, then
\[ {\cal R} \int_{u(a)}^{u(b)} f(x)dx = {\cal R} \int_a^b f(u(t))u'(t)dt, \]
with ${\cal R} \int_{u(a)}^{u(b)} f(x)dx = -{\cal R} \int_{u(b)}^{u(a)} f(x)dx$
if $u(b) < u(a)$.
\\ \smallskip\\ 
{\bf Proof:} By Proposition 3.2 (Fundamental Theorem of calculus II), with $e$ the
left endpoint of $I$, the function $F$ defined by
$F(x) = {\cal R} \int_e^x f(t) dt$, $x\in I$, satisfies $F'(x)=f(x)$ for every $x\in I$
since $f$ is continuous on~$I$. From the definition of~$F$, given $c,\,d\in I$, $c$ not
necessarily less than~$d$, then $F(d)-F(c) = {\cal R} \int_c^d f(x)dx$.
In particular,
\[ F(u(b))-F(u(a)) = {\cal R} \int_{u(a)}^{u(b)} f(x)dx.\]
Since $u$ is differentiable on $[a,b]$ and $F$ is differentiable on~$I$, the composite
function $F\circ u$ is differentiable on~$[a,b]$ and by the usual chain rule of calculus,
\mbox{$(F\circ u)'(t)$} $= F'(u(t))u'(t) = f(u(t))u'(t)$, for every $t\in [a,b]$.
Thus, since $f(u(t))u'(t)$ is clearly Riemann integrable over~$[a,b]$,
by Proposition~3.1 (Fundamental Theorem of calculus I), it must be that 
\begin{eqnarray*}
{\cal R}\int_a^b f(u(t))u'(t)dt &=& {\cal R}\int_a^b (F\circ u)'(t)dt
                                 =  (F\circ u)(b) - (F\circ u)(a)\\
                                &=& F(u(b))- F(u(a))
                                 = {\cal R} \int_{u(a)}^{u(b)} f(x)dx.
\end{eqnarray*}
{\bf Proposition 3.16 (Substitution rule for Lebesgue integral):} Given an absolutely
continuous function $u$ from an interval~$[a,b]$ into an interval~$I$, and a real-valued
continuous function~$f$ on~$I$, then
\[ \int_{u(a)}^{u(b)} f(x)dx = \int_a^b f(u(t))u'(t)dt, \]
with $\int_{u(a)}^{u(b)} f(x)dx = -\int_{u(b)}^{u(a)} f(x)dx = -\int_{[u(b),u(a)]} f(x)dx$
if $u(b) < u(a)$.
\\ \smallskip\\ 
{\bf Proof:} By Proposition 3.2 (Fundamental Theorem of calculus II), with $e$ the
left endpoint of $I$, the function $F$ defined by
$F(x) = {\cal R} \int_e^x f(t) dt$, $x\in I$, satisfies $F'(x)=f(x)$ for every $x\in I$
since $f$ is continuous on~$I$. From the definition of~$F$, given $c,\,d\in I$, $c$ not
necessarily less than~$d$, then $F(d)-F(c) = {\cal R} \int_c^d f(x)dx$.
In particular,
\[ F(u(b))-F(u(a)) = {\cal R} \int_{u(a)}^{u(b)} f(x)dx = \int_{u(a)}^{u(b)} f(x)dx, \]
where the last equation is by Proposition~2.33 (Riemann implies Lebesgue).\\
Since $u$ is differentiable a.e. on~$[a,b]$ and $F$ is differentiable on $I$, the
composite function $F\circ u$ is differentiable a.e. on~$[a,b]$. Indeed it is differentiable
exactly at the points where $u$ is differentiable. Thus, by the usual chain rule of calculus,
\mbox{$(F\circ u)'(t)$} $= F'(u(t))u'(t) = f(u(t))u'(t)$, for $t\in [a,b]$ at
which $u'$~exists.\\
Finally, we show $F\circ u$ is absolutely continuous on $[a,b]$ in order to use
Corollary~3.6 (Fundamental Theorem of Lebesgue integral calculus) with $F\circ u$ as
the absolutely continuous function in the hypothesis of the corollary.
For this purpose, since $f$ is
continuous, assume $|f|<M$ on $[a,b]$, for some $M>0$ . Given $\epsilon>0$, let $\delta>0$
correspond to $\epsilon/M$ in the definition of the absolute continuity of~$u$. Given
integer~$n>0$ and disjoint open intervals $(t_i,t'_i)\subseteq [a,b]$, $i=1,\ldots,n$,
with $\sum_{i=1}^n(t'_i-t_i)<\delta$, then
\begin{eqnarray*}
\sum_{i=1}^n|F\circ u(t'_i)-F\circ u(t_i)| &=&
\sum_{i=1}^n|{\cal R}\int_{u(t_i)}^{u(t'_i)}f(x)dx|\\
&<& \sum_{i=1}^n|{u(t_i)}-{u(t'_i)}|M
= M\sum_{i=1}^n|{u(t_i)}-{u(t'_i)}|\\
&<& M\epsilon/M=\epsilon.
\end{eqnarray*}
Thus, $F\circ u$ is absolutely continuous and by Corollary 3.6, it must be that
\begin{eqnarray*}
\int_a^b f(u(t))u'(t)dt &=& \int_a^b (F\circ u)'(t)dt
                                 =  (F\circ u)(b) - (F\circ u)(a)\\
                                &=& F(u(b))- F(u(a))
                                 = \int_{u(a)}^{u(b)} f(x)dx.
\end{eqnarray*}
{\bf Observation 3.6:} Note that in the proof of Proposition~3.16 above, while proving that
$F\circ u$ is absolutely continuous on $[a,b]$, we have actually proved that if $u$, $[a,b]$,
$I$, $e$ are as given there and $f$ is Lebesgue integrable over $I$ and bounded on $I$, and
the function $F$ is defined by $F(x) = \int_e^x f(t) dt$, $x\in I$, then $F\circ u$ is
absolutely continuous on~$[a,b]$. At the end of this section, results are presented for
carrying out a change of variable in Lebesgue integrals, useful in shape~analysis.
\\ \smallskip\\
{\bf Proposition 3.17 (Saks' inequality \cite{saks}):} Given a real-valued function $f$ on
$[a,b]$, a real number $r\geq 0$, and $E\subseteq [a,b]$ such that $|f'(x)|\leq r$ for each
$x\in E$, then \[ m^*(f(E))\leq r\,m^*(E).\]
{\bf Proof:} Given $\epsilon>0$, for every integer $n>0$, and any $y\in (a,b)$, define
\[E_n=\{x\in E:\,\mathrm{if\,} 0<|x-y|<1/n,\,\mathrm{then\,}
|f(x)-f(y)|<(r+\epsilon)|x-y|\}.\]
Since $E_n\subseteq E_{n+1}$ for all $n$ and $\cup_{n=1}^{\infty} E_n=E$, then
by 2 of Proposition~2.13, $m^*(E)= \lim_{n\rightarrow\infty} m^*(E_n)$. Similarly,
since $f(E_n)\subseteq f(E_{n+1})$ for all $n$ and $\cup_{n=1}^{\infty} f(E_n)=f(E)$,
then again by 2 of Proposition~2.13,
\[ m^*(f(E))= \lim_{n\rightarrow\infty} m^*(f(E_n)).\]
Given integer $n>0$, let $\{I_k\}$ be a countable collection of open intervals covering
$E_n$, i.e., $E_n\subseteq \cup_{k=1}^{\infty}I_k$, with
\[ \sum_{k=1}^{\infty} m(I_k) < m^*(E_n) +\epsilon. \]
Note $\{I_k\}$ can be chosen so that $m(I_k)<1/n$ for each~$k$. Then, for each~$k$,
given $x$, $x'$ $\in E_n\cap I_k$, $x\not=x'$, from the definition of $E_n$,
since $|x-x'|<1/n$, it must be that
\[ |f(x)-f(x')| < (r+\epsilon)\,|x-x'| < (r+\epsilon)\,m(I_k). \]
Thus,\[ m^* (f(E_n\cap I_k))\leq \sup_{x,x'\in E_n\cap I_k}\,|f(x)-f(x')|
\leq (r+\epsilon)\,m(I_k). \] Since
$E_n=\cup_{k=1}^{\infty}\,(E_n\cap I_k)$, then $f(E_n)=\cup_{k=1}^{\infty}\,f(E_n\cap I_k)$,
therefore,
\begin{eqnarray*}
m^*(f(E_n)) &\leq& \sum_{k=1}^{\infty}\,m^*(f(E_n\cap I_k))
\leq \sum_{k=1}^{\infty}\,(r+\epsilon)\,m(I_k)\\
&=& (r+\epsilon)\,\sum_{k=1}^{\infty}\,m(I_k) < (r+\epsilon)\,(m^*(E_n) +\epsilon).
\end{eqnarray*}
Thus, since as established above $m^*(E)= \lim_{n\rightarrow\infty} m^*(E_n)$ and
$ m^*(f(E))= \lim_{n\rightarrow\infty} m^*(f(E_n))$, it must be that
\begin{eqnarray*}
m^*(f(E))&=& \lim_{n\rightarrow\infty} m^*(f(E_n))\\
 &\leq& \lim_{n\rightarrow\infty} (r+\epsilon)\,(m^*(E_n)+\epsilon)\\
 &=& (r+\epsilon)\,(m^*(E)+\epsilon).
\end{eqnarray*}
Hence, since $\epsilon>0$ is arbitrary, it must be that $m^*(f(E)) \leq r\,m^*(E)$.
\\ \smallskip\\
{\bf Corollary 3.8:} Given a real-valued function $f$ on~$[a,b]$, let $E$ be a subset
of $[a,b]$ on which $f'=0$. Then~$m(f(E))=m^*(f(E))=0$. \\ \smallskip\\
{\bf Proof:} By Proposition 3.17, $m^*(f(E))\leq 0\cdot m^*(E)=0$. Thus, $m(f(E))=0$.
\\ \smallskip \\
{\bf Corollary 3.9:} Given a real-valued function $f$ on~$[a,b]$, and $E\subseteq [a,b]$
with $m(E)=0$, such that $f'$ exists on~$E$, then~$m(f(E))=m^*(f(E))=0$.
\\ \smallskip\\
{\bf Proof:} Without any loss of generality assume $E\subseteq (a,b)$.
Given $\epsilon >0$, let $\{I_k\}$ be a countable collection of open intervals
covering~$E$, i.e., $E\subseteq \cup_{k=1}^{\infty}I_k$, with
$\sum_{k=1}^{\infty} m(I_k) < \epsilon$, $I_k\subseteq (a,b)$ for each~$k$.\\
Given an integer $l\geq 0$, and an integer $k>0$, let
\[ E_{l} =\{x\in E:\, f'(x)\ \mathrm{exists\ and\ } |f'(x)|\leq l\} \]
and
\[ E_{lk} =\{x\in I_k:\, f'(x)\ \mathrm{exists\ and\ } |f'(x)|\leq l\}.\]
Note $E_l\subseteq \cup_{k=1}^{\infty}E_{lk}$
so that $f(E_l)\subseteq \cup_{k=1}^{\infty}f(E_{lk})$.\\
By Proposition 3.17, $m^*(f(E_{lk}))\leq l\,m^*(E_{lk})$ for each~$k$, and therefore,
\[ m^*(f(E_l))\leq\sum_{k=1}^{\infty} m^*(f(E_{lk}))\leq\sum_{k=1}^{\infty}l\, m^*(E_{lk})
\leq\sum_{k=1}^{\infty}l\,m(I_k)<l\,\epsilon. \]
Since $\epsilon>0$ is arbitrary, then $m^*(f(E_l))=0$.\\
Finally, note $f(E_{l+1})\supseteq f(E_l)$ for each~$l$,
and $f(E)=\cup_{l=0}^{\infty}\,f(E_l)$.\\
Thus, by 2 of Proposition 2.13, $m^*(f(E))=\lim_{\,l\rightarrow\infty} m^*(f(E_l))=0$.
\\ \smallskip\\
{\bf Corollary 3.10 (Saks' Theorem \cite{saks}):} Given a real-valued function $f$
on~$[a,b]$, and $E\subseteq [a,b]$ such that $f'$ exists on~$E$, if $f'=0$ a.e. on~$E$,
then~$m(f(E))=0$.\\ \smallskip\\
{\bf Proof:} Let $E_1$, $E_2$ be subsets of $E$, $E=E_1\cup E_2$, $f'=0$ on~$E_1$ and
$m(E_2)=0$. By Corollary~3.8, $m(f(E_1))=0$. By Corollary~3.9, $m(f(E_2))=0$.
Thus, $m(f(E)) \leq m(f(E_1))+m(f(E_2))= 0$. \\ \smallskip\\
{\bf Observation 3.7:} Let $f$ be the Cantor function and $C$ the Cantor set.
Then $f'=0$~a.e. on~$[0,1]$ and~$m(C)=0$. Since $f([0,1])=[0,1]$, by Saks' Theorem
(Corollary~3.10), it must be that $f$ is not differentiable at certain points in~$[0,1]$,
and since $f(C)=[0,1]$, by Corollary~3.9, it must be that $f$ is not differentiable
at certain points in~$C$. Of course we know $f$ is not differentiable at any point
in~$C$ and $f'=0$ on~$[0,1]\setminus C$ so that by Corollary~3.8, $m(f([0,1]\setminus C))=0$
which makes sense as $f([0,1]\setminus C)$ is countable.\\
The following proposition is the converse of Saks' Theorem (Corollary 3.10): $m(f(E))=0$
implies $f'=0$~a.e. on~$E$. Here again it is assumed $f'$ exists everywhere on~$E$.
However, almost everywhere (a.e.) will suffice.\\ \smallskip\\
{\bf Proposition 3.18 (Serrin-Varberg's Theorem \cite{serrin}):} Given $f$, a real-valued
function on~$[a,b]$, and $E\subseteq [a,b]$ such that $f'$ exists on $E$, if $m(f(E))=0$,
then $f'=0$~a.e. on~$E$.
\\ \smallskip\\
{\bf Proof:} Let $B=\{x\in E:\, |f'(x)|>0 \}$, and for every integer $n>0$,
and any $y\in (a,b)$, define
\[B_n=\{x\in B:\,\mathrm{if\,} 0<|x-y|<1/n,\,\mathrm{then\,}
|f(x)-f(y)|\geq |x-y|/n \}.\]
Clearly, we need to show $m(B)=0$, and since $B=\cup_{n=1}^{\infty}\,B_n$, then it suffices
to show $m(B_n)=0$ for each~$n$. However, since each $B_n$ can be covered by a countable
collection of intervals, each interval of length less than~$1/n$, it then suffices to show that if
$I$ is any interval of length less than~$1/n$ and $A=I\cap B_n$, then~$m(A)=0$.\\
For this purpose, given $\epsilon>0$, since $A\subseteq E$ so that $m(f(A))=0$, let $\{I_k\}$
be a countable collection of open intervals covering~$f(A)$, i.e.,
$f(A)\subseteq \cup_{k=1}^{\infty}\,I_k$, with $\sum_{k=1}^{\infty}\,m(I_k)<\epsilon$.
In addition, let $A_k= f^{-1}(f(A)\cap I_k)$. Then $A=\cup_{k=1}^{\infty} A_k$, and
since for each $k$, $A_k\subseteq A=I\cap B_n$, given $x$, $x'\in A_k$, $x\not=x'$, then
$|x-x'|<1/n$, and it must be that \[ |x-x'|\leq n |f(x)-f(x')|< n\,m(I_k).\]
Thus, \[ m^*(A_k) \leq \sup_{x,x'\in A_k} |x-x'| \leq n\,m(I_k), \] and then
\[ m^*(A)\leq\sum_{k=1}^{\infty} m^*(A_k)\leq \sum_{k=1}^{\infty}n\,m(I_k)=
n\sum_{k=1}^{\infty}m(I_k)<n\epsilon. \]
Hence, since $\epsilon >0$ is arbitrary, it must be that $m(A)=0$. \\ \smallskip\\
{\bf Corollary 3.11:} Given a real-valued function $f$ on $[a,b]$, and $E\subseteq [a,b]$
such that $f'$ exists on~$E$, if $f$  is constant on~$E$, then $f'=0$ a.e. on~$E$.
\\ \smallskip\\
{\bf Corollary 3.12 (Serrin-Varberg's Theorem (Alternate form)):}
Given a real-valued function $f$ on~$[a,b]$, and $E\subseteq [a,b]$
such that $f'$ exists a.e. on $E$, if $m(f(E))=0$, then $f'=0$ a.e.~on~$E$.
In particular, if $f$ is of bounded variation on $[a,b]$, and $E\subseteq [a,b]$
with $m(f(E))=0$, then $f'=0$ a.e. on~$E$. \\ \smallskip\\
{\bf Proof:}
Let $E_1$, $E_2$ be subsets of $E$, $E=E_1\cup E_2$, $f'$ exists on $E_1$
and~$m(E_2)=0$. Since $m(f(E_1))\leq m(f(E))=0$, then by Proposition~3.18, $f'=0$ a.e.
on~$E_1$. Thus, since $m(E_2)=0$, then $f'=0$ a.e. on~$E$.
If $f$ is of bounded variation on $[a,b]$, then by Corollary 3.3, $f$ is differentiable
a.e. on~$[a,b]$ and therefore on any subset~$E$ of~$[a,b]$.
\\ \smallskip\\
{\bf Proposition 3.19 (Measurability of the derivative of a measurable function):}
Let $f$ be a real-valued measurable function on~$[a,b]$, and $E$ a measurable subset
of~$[a,b]$. If $f'$ exists on~$E$, then $f'$ is a measurable function on~$E$.
Proof in~\cite{yeh}.\\ \smallskip\\
{\bf Proposition 3.20:} Let $f$ be a real-valued measurable function on~$[a,b]$, and $E$
a measurable subset of~$[a,b]$. If $f'(x)$ exists for each $x\in E$, then
\[ m^*(f(E))\leq \int_{E} |f'(x)|dx. \]
{\bf Proof:} Given $\epsilon>0$, for each integer $n>0$, define
\[ E_n=\{x\in E: (n-1)\,\epsilon\leq |f'(x)|<n\,\epsilon\}. \]
Clearly, the $E_n$'s are pairwise disjoint and $E=\cup_{n=1}^{\infty}E_n$ so that
$f(E)=\cup_{n=1}^{\infty}f(E_n)$.
By Proposition~3.19, $f'$~is measurable on~$E$. Hence, each set~$E_n$ must be
measurable. Since $f'(x)$ exists for every $x\in E_n$, by Proposition~3.17,
$m^*(f(E_n))\leq n\epsilon\,m(E_n)$ for each~$n$, and given an integer~$N>0$,
it must be that $\sum_{n=1}^N (n-1)\epsilon\,m(E_n)\leq \int_E |f'(x)|dx$
by the definition of the Lebesgue integral, so that
$\sum_{n=1}^{\infty} (n-1)\epsilon\,m(E_n)\leq \int_E |f'(x)|dx$.
Thus,
\begin{eqnarray*}
m^*(f(E)) &\leq& \sum_{n=1}^{\infty} m^*(f(E_n))\leq
\sum_{n=1}^{\infty}n\epsilon\,m(E_n)= \sum_{n=1}^{\infty}((n-1)\epsilon +\epsilon)m(E_n)\\
&=& \sum_{n=1}^{\infty}(n-1)\epsilon\,m(E_n)+ \sum_{n=1}^{\infty}\epsilon\,m(E_n)\leq
\int_E |f'(x)|dx +\epsilon\,m(E),
\end{eqnarray*}
where, in the last step, the countable additivity of $m$ on measurable sets is~used.
Since $\epsilon$ is arbitrary and $m(E)<\infty$, then $m^*(f(E))\leq \int_{E} |f'(x)|dx$.
\\ \smallskip\\
{\bf Proposition 3.21 (Absolutely continuous $f$ maps zero-measure sets to zero-measure
sets~\cite{saks} - Absolutely continuous $f$ maps measurable sets to measurable sets):} 
Let $f$ be an absolutely continuous function on~$[a,b]$. If $E\subseteq [a,b]$ with
$m(E)=0$, then~$m(f(E))=0$. In addition, given any measurable subset $E$ of~$[a,b]$,
then~$f(E)$ is measurable.
\\ \smallskip\\
{\bf Proof:} Without any loss of generality assume $E\subseteq (a,b)$.
Given $\epsilon>0$, let $\delta>0$ correspond to $\epsilon$ in the definition of the
absolute continuity of~$f$. Since $m(E)=0$, then by Proposition~2.1, there is a
collection $\{I_k\} = \{(a_k,b_k)\}$ of nonoverlapping open intervals covering~$E$,
i.e., $E\subseteq \cup_{k=1}^{\infty} I_k$, with
$\sum_{k=1}^{\infty}(b_k-a_k)=\sum_{k=1}^{\infty} m(I_k)<\delta$, $I_k\subseteq (a,b)$
for each~$k$.\\
For each $k$, since $f$ is continuous, let $c_k$ and $d_k$ be points in $[a_k,b_k]$
where $f$ attains its minimum and maximum, respectively. Assuming without any loss of
generality that $c_k<d_k$, then $\{(c_k,d_k)\}$ is a collection of nonoverlapping open
intervals, and since again $f$ is continuous, by the intermediate value theorem
\cite{rudin}, it follows that
\[ f(E)\subseteq f(\cup_{k=1}^{\infty}(a_k,b_k))=
\cup_{k=1}^{\infty}f((a_k,b_k)) \subseteq \cup_{k=1}^{\infty}f([c_k,d_k]). \]
Thus, \[ m^*(f(E))\leq \sum_{k=1}^{\infty} m(f([c_k,d_k]) =
\sum_{k=1}^{\infty}(f(d_k)-f(c_k)). \]
Finally, given an integer $N>0$, then it must be that
\[ \sum_{k=1}^N (d_k-c_k)< \sum_{k=1}^{\infty}(d_k-c_k)\leq
\sum_{k=1}^{\infty}(b_k-a_k))<\delta \]
so that $\sum_{k=1}^N(f(d_k)-f(c_k))<\epsilon$. 
$N$ arbitrary then implies
$m^*(f(E))\leq\sum_{k=1}^{\infty}(f(d_k)-f(c_k))\leq\epsilon$.
Thus, $m^*(f(E))=0$ since $\epsilon$ is arbitrary.\\ \smallskip\\
Assume now $E$ is a measurable subset of~$[a,b]$.\\
By v of Proposition~2.15, there is a set $F$ that is the union of a countable collection of
closed sets, $F\subseteq E$ with $m^*(E\setminus F)=0$, i.e., $E=K\cup(\cup_{n=1}^{\infty}F_n)$,
where $m^*(K)=0$, $F=\cup_{n=1}^{\infty}F_n$, and $F_n$ is closed for each~$n$, thus compact
(Proposition~2.3). Note then that
$f(E)=f(K\cup(\cup_{n=1}^{\infty}F_n))=f(K)\cup(\cup_{n=1}^{\infty}f(F_n))$.
Since $m(K)=0$, then $m(f(K))=0$ as just proved above. Thus, $f(K)$~is measurable. Also since
$f$~is continuous and $F_n$~is compact for each~$n$, it must be that $f(F_n)$ is compact for
each~$n$~\cite{rudin} and thus measurable (Proposition~2.14). It then follows that $f(E)$
is the union of a countable collection of measurable sets and therefore it must be measurable.
\\ \smallskip\\
{\bf Proposition 3.22 (Banach-Zarecki Theorem):} Let $f$ be a real-valued function on $[a,b]$. Then
$f$ is absolutely continuous on~$[a,b]$ if and only if it satisfies the following three
conditions:\\
i. $f$ is continuous on $[a,b]$.\\
ii. $f$ is of bounded variation on $[a,b]$.\\
iii. $f$ maps sets of measure zero to sets of measure zero.\\ \smallskip\\
{\bf Proof:} The necessity was established in Observation~3.3, Proposition~3.10,
and Proposition~3.21. For the sufficiency, assume $f$ satisfies all three conditions.
Given $[c,d]\subseteq [a,b]$, we show \[ |f(d)-f(c)|\leq \int_c^d |f'(x)|dx. \]
By condition ii, $f$ is of bounded variation on $[a,b]$, so that by Corollary~3.3, $f$~is
differentiable a.e. on~$[a,b]$ and therefore on~$[c,d]$. Accordingly, let $E_1$, $E_2$
be subsets of~$[c,d]$, $[c,d] = E_1\cup E_2$, $f'$ exists on~$E_1$ and $m(E_2)=0$.
By condition iii, it then must be that $m(f(E_2))=0$.\\
Assume without any loss of generality that $f(c)<f(d)$. By condition i, $f$ is continuous
on $[a,b]$, so that by the intermediate value theorem \cite{rudin}, given $y$, $f(c)<y<f(d)$,
there must be $x$, $c<x<d$, with $f(x)=y$. Thus, $[f(c),f(d)]\subseteq f([c,d])$, and since
$f$ is measurable on $[a,b]$ ($f$ is continuous on $[a,b]$) and $E_1$ is measurable
($E_1=[c,d]\setminus E_2$), by Proposition 3.20, we get
\begin{eqnarray*}
|f(d)-f(c)|&\leq& m(f([c,d]))=m(f(E_1)\cup f(E_2))\\
&\leq& m^*(f(E_1))+m(f(E_2))=m^*(f(E_1))+0\\
&=& m^*(f(E_1))\leq\int_{E_1}|f'(x)|dx=\int_c^d |f'(x)|dx.
\end{eqnarray*}
By Corollary 3.3, since $f$ is of bounded variation, it must be that $f'$ is integrable
over~$[a,b]$ and so is $|f'|$ by Proposition~2.31. Given $\epsilon>0$, by Proposition~2.38,
there is $\delta>0$ such that if $A$ is a measurable set with $m(A)<\delta$,
then~$\int_A|f'(x)|dx<\epsilon$. Accordingly, for any integer $n>0$ and any disjoint
collection of open intervals $(x_i,x_i') \subseteq [a,b]$, $i=1,\ldots,n$, with
$\sum_{i=1}^n (x_i' - x_i) < \delta$, let $A=\cup_{i=1}^n (x_i' - x_i)$.
Since $m(A)<\delta$,
then \[ \sum_{i=1}^n |f(x_i') - f(x_i)|\leq\sum_{i=1}^n\int_{x_i}^{x_i'}|f'(x)|dx
=\int_A |f'(x)|dx < \epsilon. \]
Thus, $f$ is absolutely continuous.
\\ \smallskip\\
{\bf Proposition 3.23 (Inverse function theorem):} Let $f$ be a strictly monotonic
continuous function on~$[a,b]$. Then $I=f([a,b])$ is a closed interval with endpoints
$f(a)$, $f(b)$, and $f^{-1}$, the inverse function of~$f$, exists on~$I$, and is
strictly monotonic and continuous on~$I$. Given $x_0\in [a,b]$ such that $f$ is
differentiable at~$x_0$ with~$f'(x_0)\not=0$, then $f^{-1}$ is differentiable
at $y_0 = f(x_0)$ with \[ (f^{-1})'(y_0) = 1/f'(x_0). \]
{\bf Proof:} Without any loss of generality, assume $f$ is increasing on~$[a,b]$.
Since $f$ is strictly increasing and continuous, $f$ is one-to-one and by the
intermediate value theorem \cite{rudin}, its range is $I=[f(a),f(b)]=f([a,b])$. Thus,
$f^{-1}$ exists on~$I$ and is strictly increasing from $I$ onto~$[a,b]$. By Corollary~3.2,
$f^{-1}$ is continuous on~$I$. With $y_0=f(x_0)$, $y=f(x)$, by the continuity
of~$f^{-1}$, if~$y\rightarrow y_0$, it must be that~$x\rightarrow x_0$. Thus,
\[ \frac{f^{-1}(y) -f^{-1}(y_0)}{y-y_0} =
\frac{x-x_0}{f(x)-f(x_0)}\rightarrow 1/f'(x_0) \]
as $y\rightarrow y_0$, since $f$ is differentiable at~$x_0$, and~$f'(x_0)\not=0$.
Hence, $(f^{-1})'(y_0)$ exists and equals~$1/f'(x_0)$.
\\ \smallskip\\
{\bf Proposition 3.24 (Zarecki's criterion for an absolutely continuous
inverse~\cite{cabada}):} Let $f$ be a monotonic continuous function on~$[a,b]$.
Then $I=f([a,b])$ is a closed interval with endpoints $f(a)$, $f(b)$ , and $f^{-1}$
exists and is absolutely continuous on~$I$ if and only if $\{x:f'(x)=0\}$ has
measure~zero. Whenever $f^{-1}$ is absolutely continuous on the closed interval~$I$,
then \[ (f^{-1})'= 1/(f'(f^{-1})) \mathrm{\ a.e.\ on\ } I. \]
{\bf Proof:} Without any loss of generality, assume $f$ is increasing on~$[a,b]$.
If~$f^{-1}$  exists or if $m^*(\{x:f'(x)=0\})=0$, then $f$ is strictly increasing.
Thus, assume $f$ is strictly increasing.
By Proposition~3.23, $I$ is closed, $I=[f(a),f(b)]$, and $f^{-1}$ exists on~$I$
and is strictly increasing and continuous on~$I$.\\
Assume $\{x:f'(x)=0\}$ has measure zero. As already established, $f^{-1}$ is continuous
on~$I$, and
since it is increasing, it is of bounded variation on~$I$. Thus, by Proposition~3.22,
it suffices to show that $f^{-1}$ maps sets of measure zero to sets of measure zero.
For this purpose, let $E\subseteq I$ be of measure zero, and $F=f^{-1}(E)$ so that~$f(F)=E$.
Since $f$ is increasing, it is of bounded variation on~$[a,b]$ as well, and by
Corollary~3.12, $f'=0$ a.e. on~$F$. Accordingly, let $F_1$, $F_2$ be subsets of~$F$,
$F=F_1\cup F_2$, $f'=0$~on $F_1$ and~$m(F_2)=0$. Since $F_1\subseteq \{x:f'(x)=0\}$,
then~$m(F_1)=0$. Thus,~$m(F)=0$.\\
Assume now $f^{-1}$ is absolutely continuous on~$I$. By Corollary~3.8, since $f'=0$ on
$E=\{x:f'(x)=0\}$, then~$m(f(E))=0$. Thus, by Proposition~3.21, $m(E)=m(f^{-1}(f(E)))=0$.\\
Finally, whenever $f^{-1}$ is absolutely continuous on~$I$,
define $F_1$, $F_2$, $F_3$, disjoint subsets of~$[a,b]$, $[a,b]=F_1\cup F_2\cup F_3$,
as follows.  $F_1=\{x:f'(x)\not=0\}$, $F_2=\{x:f'(x)=0\}$,
$F_3=\{x:f'(x) \mathrm{\ does\ not\ exist}\}$. Since $f$ is of bounded variation
on~$[a,b]$, $f$ is differentiable a.e. on~$[a,b]$, thus, $m(F_3)=0$. Also, since $f^{-1}$
is absolutely continuous, then $m(F_2)=0$ as just proved above. Hence, since by
Proposition~3.23, for each $x\in F_1$, $f^{-1}$ is differentiable at $y=f(x)$ with 
$(f^{-1})'(y) = 1/f'(x)$, then $(f^{-1})'= 1/(f'(f^{-1}))$ a.e. on~$I$.
\\ \smallskip\\
{\bf Proposition 3.25 (Composition of absolutely continuous functions):} Let $g$ be an
absolutely continuous monotonic function from an interval~$[a,b]$ into an
interval~$[c,d]$, and let $f$ be an absolutely continuous function on~$[c,d]$.
Then the function $h= f\circ g$ is absolutely continuous on~$[a,b]$.\\ \smallskip\\
{\bf Proof:} Let $\epsilon>0$ be given. It follows easily that since $f$ is absolutely
continuous, then there is $\rho>0$ such that for any integer $n>0$ and for any
nonempty subset $A$ of $\{1,\ldots,n\}$ it must be that
$\sum_{i\in A} |f(y_i') - f(y_i)| < \epsilon$ for disjoint open intervals
$(y_i,y_i') \subseteq [c,d]$, $i\in A$, with $\sum_{i\in A} (y_i' - y_i) < \rho$.
For $\rho$ as just described, since $g$ is also absolutely continuous, then there is
$\delta>0$ such that for any integer $n>0$ it must be that 
$\sum_{i=1}^n |g(x_i') - g(x_i)| < \rho$ for disjoint open intervals
$(x_i,x_i') \subseteq [a,b]$, $i=1,\ldots,n$, with $\sum_{i=1}^n (x_i' - x_i) < \delta$.
Accordingly, for any integer $n>0$ let $(x_i,x_i') \subseteq [a,b]$, $i=1,\ldots,n$,
be any collection of disjoint open intervals with $\sum_{i=1}^n (x_i' - x_i) < \delta$.
Setting $A=\{i:\,1\leq i\leq n,\,g(x_i')\not=g(x_i)\}$, if $A\not=\emptyset$, then
the collection of open intervals $(g(x_i),g(x_i'))$ , $i\in A$, if $g$ is increasing;
$(g(x_i'),g(x_i))$, $i\in A$, if $g$ is decreasing; must be pairwise disjoint by
the monoticity of~$g$ with $\sum_{i\in A}|g(x_i') - g(x_i)|<\rho$.
Thus, if $A\not=\emptyset$, it must be that
$S_A=\sum_{i\in A}|f(g(x_i')) - f(g(x_i))|<\epsilon$. Setting $S_A=0$ if~$A=\emptyset$,
and since $|f(g(x_i')) - f(g(x_i))|=0$ for $i\not\in A$, then
\[ \sum_{i=1}^n|f(g(x_i')) - f(g(x_i))|= S_A +0 =S_A<\epsilon.\]
Thus, $h=f\circ g$ is absolutely continuous on $[a,b]$.
\\ \smallskip\\
{\bf Proposition 3.26 (Chain rule \cite{serrin}):} Given real-valued functions $F$, $f$
on~$[c,d]$, $F'=f$ a.e. on~$[c,d]$, and a function $u:[a,b]\rightarrow [c,d]$,
$u$ and $F\circ u$ differentiable a.e. on~$[a,b]$, if $F$ maps zero-measure sets to
zero-measure sets, then
\[(F\circ u)'=  (f\circ u)u' \mathrm{\ a.e.\ on\ } [a,b].\]
{\bf Proof:} Let $A=\{x\in [c,d]: F'(x)=f(x)\}$, $B=[c,d]\setminus A$, and
$C=\{t\in [a,b]:u(t)\in B\}$. Clearly, $mB=0$. Letting $D=[a,b]\setminus C$, since $u$ is
differentiable a.e. on~$[a,b]$, then it is differentiable a.e. on~$D$. Since $u(D)\subseteq A$
and $F$ is differentiable on~$A$, the composite function $F\circ u$ is differentiable a.e.
on~$D$. Indeed it is differentiable exactly at the points in~$D$ where $u$ is
differentiable. Thus, by the usual chain rule of calculus,
$(F\circ u)'(t) = F'(u(t))u'(t) = f(u(t))u'(t)$, for $t\in D$ at which $u'$ exists,
i.e., $(F\circ u)'=  (f\circ u)u' \mathrm{\ a.e.\ on\ }D$.\\
Note that if $mC=0$, then the proof is complete.\\
Thus, assuming $mC\not=0$, we show $(F\circ u)'=  (f\circ u)u' \mathrm{\ a.e.\ on\ }C$.
Note $f\circ u$ is defined on~$C$. Since $u(C)\subseteq B$, then $m(u(C))=mB=0$, and since
$u'$ exists a.e.  on $C$, then by Corollary~3.12, $u'=0$ a.e. on~$C$ so that $(f\circ u)u'=0$
a.e. on~$C$.  In addition, since $F$ maps zero-measure sets to zero-measure sets, then
$m(F(u(C)))=0$, and since $(F\circ u)'$ exists a.e. on $C$, again by Corollary~3.12,
$(F\circ u)' =0$ a.e.  on~$C$. Thus, $(F\circ u)'=0 = (f\circ u)u'$ a.e. on~$C$.
\\ \smallskip\\
{\bf Observation 3.8:} With $B$ the subset of~$[c,d]$ of measure zero on which~$F'\not=f$,
and $C=\{t\in [a,b]:u(t)\in B\}$, then for the case $mC\not=0$ in the proof of
Proposition~3.26 above it was proved that $u'=0$ a.e. on~$C$, and since $f\circ u$ is defined
on~$C$, then $(f\circ u)u'=0$ a.e. on $C$ as well. Since it was also proved
that $(F\circ u)'=0$ a.e. on~$C$, then it was concluded that $(F\circ u)'=0=(f\circ u)u'$
a.e. on~$C$. However, it can happen that when using the chain rule as described in
Proposition~3.26 above and in Corollary 3.13 below, although $F'=f$ on $[c,d]\setminus B$,
$f$ may not be defined everywhere on~$B$. Thus, $f\circ u$ may not be defined on $C$ as
required in the proof of Proposition~3.26 above when $mC\not=0$. However, what matters here
is that both $u'$ and $(F\circ u)'$ are zero a.e. on~$C$. Therefore, assuming $mC\not=0$,
when computing $(F\circ u)'$ with the chain rule as suggested in Proposition~3.26 above and
in Corollary~3.13 below, if $(F\circ u)'$ is set to zero at any point in~$[a,b]$ at which $u'$
is zero (whether or not $f\circ u$ is defined there), and computed or left undefined according
to the chain rule elsewhere, then $(F\circ u)'$ so obtained will be correct a.e. on~$[a,b]$.
Accordingly, assuming $mC\not=0$, one should keep in mind that if $(F\circ u)'$ is not
computed as just suggested, so that $(F\circ u)'$ might be left undefined at points where
$f\circ u$ is not defined although $u'$ is zero, one could end up with $(F\circ u)'$ not
defined on a set of nonzero measure in~$[a,b]$. Actually, instead of computing $(F\circ u)'$
as just suggested, we do something simpler. Since, as mentioned above, what matters here
is that both $u'$ and $(F\circ u)'$ are zero a.e. on~$C$, without any loss of generality, we
simply set $f$ equal to 1 at points in $B$ where it is not defined and proceed with the
chain rule to compute $(F\circ u)'$, as $f\circ u$ is then defined on~$C$ so that
$(f\circ u)u'$ is zero at points in~$C$ where $u'$ is zero, thus zero a.e.  on~$C$.
More precisely, we define a new function $\hat{f}$ on~$[c,d]$  by setting $\hat{f}$ equal to
$f$ at points in $[c,d]$ where $f$ exists, and to 1 where it does not.
In what follows, we will refer to $\hat{f}$ as $f$ {\bf extended to all of~$[c,d]$}.
This function is then defined everywhere in $[c,d]$, equals $f$ a.e. on~$[c,d]$, and takes
the place of $f$ in the chain rule although it is still called $f$ there.
Finally, note that above when we say anything about computing $(F\circ u)'$ with the chain rule,
it is not $(F\circ u)'$ that is necessarily computed but a function that happens to be equal to
$(F\circ u)'$ a.e. on~$[a,b]$. \\ \smallskip\\
{\bf Corollary 3.13 (Chain rule (Alternate form) \cite{serrin}):} Given an absolutely
continuous function $F$ on $[c,d]$, and a real-valued function $f$ on $[c,d]$,
$F'=f$ a.e. on~$[c,d]$, if $u:[a,b]\rightarrow [c,d]$ is a function such that
$u$ and $F\circ u$ are differentiable a.e. on~$[a,b]$, then
\[(F\circ u)'=  (f\circ u)u' \mathrm{\ a.e.\ on\ } [a,b].\]
{\bf Proof:} From Proposition 3.26 and Proposition 3.21. \\ \smallskip\\
{\bf Proposition 3.27 (Change of variable for Lebesgue integral \cite{serrin}):} Given a
function $f$, Lebesgue integrable over~$[c,d]$, and a function $u:[a,b]\rightarrow [c,d]$,
differentiable a.e. on~$[a,b]$, then the following two conditions are equivalent,
where $F(x)=\int_c^x f(t)dt$, $x\in [c,d]$:\\
i. $F\circ u$ is absolutely continuous on $[a,b]$.\\
ii. $(f\circ u)u'$ is Lebesgue integrable over~$[a,b]$ and
for all $\alpha,\,\beta\in [a,b]$ it must be that
\[ \int_{u(\alpha)}^{u(\beta)} f(x)dx = \int_{\alpha}^{\beta} f(u(t))u'(t)dt, \]
with $\int_{u(\alpha)}^{u(\beta)} f(x)dx = -\int_{u(\beta)}^{u(\alpha)}
f(x)dx = -\int_{[u(\beta),u(\alpha)]} f(x)dx$ if $u(\beta) < u(\alpha)$,
and 
 $\int_{\alpha}^{\beta} f(u(t))u'(t)dt=-\int_{\beta}^{\alpha} f(u(t))u'(t)dt
 =-\int_{[\beta,\alpha]} f(u(t))u'(t)dt$ if $\beta<\alpha$. 
\\ \smallskip\\
{\bf Proof:}\\
i $\Rightarrow$ ii:\\
Since $F$ is absolutely continuous on~$[c,d]$ (Proposition 3.11), $F'=f$ a.e. on $[c,d]$
(Proposition~3.7), $u$ is differentiable a.e. on~$[a,b]$, and $F\circ u$, being absolutely
continuous on~$[a,b]$, must be differentiable a.e. on~$[a,b]$ (Corollary~3.5), then by
Corollary~3.13 (chain rule), $(F\circ u)'=  (f\circ u)u'$  a.e. on~$[a,b]$ (here and in the
corollaries that follow, without any loss of generality, $f$ is interpreted as $f$
extended to all of~$[c,d]$ (Observation~3.8 about the chain rule)).
Note, by Corollary~3.5, since $F\circ u$ is absolutely
continuous on~$[a,b]$, then $(f\circ u)u'$ is Lebesgue integrable over~$[a,b]$, and
by Corollary~3.6 (Fundamental Theorem of Lebesgue integral calculus), applied to the
absolutely continuous function $F\circ u$, for all $\alpha,\,\beta\in [a,b]$
it must be that
\begin{eqnarray*}
\int_{\alpha}^{\beta} f(u(t))u'(t)dt &=& \int_{\alpha}^{\beta} (F\circ u)'(t)dt
                                 =  (F\circ u)(\beta) - (F\circ u)(\alpha)\\
                                &=& F(u(\beta))- F(u(\alpha))
                                 = \int_{u(\alpha)}^{u(\beta)} f(x)dx.
\end{eqnarray*}
ii $\Rightarrow$ i:\\
Since $(f\circ u)u'$ is Lebesgue integrable over~$[a,b]$, and, in particular,
for $x\in [a,b]$
\[ F(u(x))-F(u(a)) = \int_{u(a)}^{u(x)} f(s)ds =
\int_a^{x} f(u(t))u'(t)dt, \]
by Proposition 3.12, $F\circ u = F(u)$ must be absolutely continuous on~$[a,b]$.\\ \smallskip\\
{\bf Corollary 3.14 (Change of variable for Lebesgue integral (Alternate form I)
\cite{serrin}):} Given a function $f$, Lebesgue integrable
over~$[c,d]$, and a function $u:[a,b]\rightarrow [c,d]$, monotonic and absolutely
continuous on~$[a,b]$, then $(f\circ u)u'$ is Lebesgue integrable over~$[a,b]$ and
for all $\alpha,\,\beta\in [a,b]$ it must be that
\[ \int_{u(\alpha)}^{u(\beta)} f(x)dx = \int_{\alpha}^{\beta} f(u(t))u'(t)dt. \]
{\bf Proof:} Since $u$ is clearly differentiable a.e. on $[a,b]$, and $F$,
$F(x)=\int_c^x f(t)dt$, $x\in [c,d]$, is absolutely continuous so that the
composition $F\circ u$ is absolutely continuous on~$[a,b]$ by Proposition~3.25,
then by Proposition~3.27, $(f\circ u)u'$ is Lebesgue integrable over~$[a,b]$
and for all $\alpha,\,\beta\in [a,b]$ it must be that
\[ \int_{u(\alpha)}^{u(\beta)} f(x)dx = \int_{\alpha}^{\beta} f(u(t))u'(t)dt. \]
{\bf Corollary 3.15 (Change of variable for Lebesgue integral (Alternate form II) 
\cite{serrin}):} Given a function $f$, bounded and measurable on~$[c,d]$, and a function
$u:[a,b]\rightarrow [c,d]$, absolutely continuous on~$[a,b]$, then $(f\circ u)u'$
is Lebesgue integrable over~$[a,b]$ and for all $\alpha,\,\beta\in [a,b]$ it must be that
\[ \int_{u(\alpha)}^{u(\beta)} f(x)dx = \int_{\alpha}^{\beta} f(u(t))u'(t)dt. \]
{\bf Proof:} By Proposition~2.32, $f$ is Lebesgue integrable over~$[c,d]$. Since $u$ is
clearly differentiable a.e. on~$[a,b]$, and by Observation~3.6, with $F(x)=\int_c^x f(t)dt$,
$x\in [c,d]$, it must be that $F\circ u$ is absolutely continuous on~$[a,b]$,
then by Proposition~3.27, $(f\circ u)u'$ is Lebesgue integrable over~$[a,b]$
and for all $\alpha,\,\beta\in [a,b]$ it must be that
\[ \int_{u(\alpha)}^{u(\beta)} f(x)dx = \int_{\alpha}^{\beta} f(u(t))u'(t)dt. \]
{\bf Corollary 3.16 (Change of variable for Lebesgue integral over a measurable
set~\cite{klassen}:)} Given $A$, a measurable subset of~$[0,1]$, and a function
$\gamma: [0,1]\rightarrow [0,1]$, $\gamma$ absolutely continuous on~$[0,1]$,
$\dot{\gamma}>0$ a.e. on~$[0,1]$, $\gamma(0)=0$, $\gamma(1)=1$, then $\gamma^{-1}$
exists and is absolutely continuous on~$[0,1]$, and $\tilde{A}=\gamma^{-1}(A)$ is a measurable
subset of~$[0,1]$. Accordingly, given a function~$f$, Lebesgue integrable over~$[0,1]$,
then $(f\circ \gamma)\dot{\gamma}$ is Lebesgue integrable over~$\tilde{A}$ and
\[ \int_A f(x)dx = \int_{\tilde{A}} f(\gamma(t))\dot{\gamma}(t)dt. \]
{\bf Proof:} Clearly $\gamma$ is strictly increasing and thus $\gamma^{-1}$ exists and is
absolutely continuous on~$[0,1]$ by Proposition~3.24. By Proposition 3.21,
$\tilde{A}=\gamma^{-1}(A)$ is then a measurable subset of~$[0,1]$.\\
Define $I_A:[0,1]\rightarrow {\bf R}$ by $I_A(t)=1$ if $t\in A$, $I_A(t)=0$
if $t\in [0,1]\setminus A$, and $I_{\tilde{A}}:[0,1]\rightarrow {\bf R}$ by
$I_{\tilde{A}}(t)=1$ if $t\in \tilde{A}$, $I_{\tilde{A}}(t)=0$
if $t\in [0,1]\setminus\tilde{A}$.\\
Note $I_{\tilde{A}}= I_A\circ\gamma$.\\
Also note $I_A\cdot f$ is Lebesgue integrable over~$[0,1]$, since it equals $f$ on~$A$
and~$0$ on~$[0,1]\setminus A$. It follows then by Corollary~3.14 that
\[ \int_0^1 I_A(x)f(x)dx = \int_{\gamma(0)}^{\gamma(1)} I_A(x)f(x)dx=
\int_0^1 I_A(\gamma(t))f(\gamma(t))\dot{\gamma}(t)dt, \]   i.e.,
\[ \int_Af(x)dx = \int_0^1 I_{\tilde{A}}(t)f(\gamma(t))\dot{\gamma}(t)dt
= \int_{\tilde{A}}f(\gamma(t))\dot{\gamma}(t)dt. \]

%%%%%%%%%%%%%%%%%%%%%%%%%%%
%%%%%%%%%%%%%%%%%%%%%%%%%%%%%%%%%%%%%
%%%%%%%%%%%%%%%%%%%
\section{\large Functional Data and Shape Analysis and its\\ Connections to Lebesgue
Integration and\\ Absolute Continuity}
{\bf Observation 4.1:} In what follows, we review some important aspects of functional
data and shape analysis of the type in~\cite{srivastava}, while at the same time
pointing out its dependence on Lesbesgue integration, absolute continuity and the
connections between them. As in~\cite{srivastava} where absolutely continuous functions
on~$[0,1]$ are generalized to functions of range ${\bf R^n}$, {\bf R} the set of real
numbers, $n$~a positive integer, we consider absolutely continuous functions
on~$[0,1]$ but restrict ourselves to those with range in~${\bf R}^1= {\bf R}$.
We denote by~$AC[0,1]$ the set of such functions. With two absolutely continuous
functions on~$[0,1]$ considered equal if they differ by a constant, we note that
the principal goal in~\cite{srivastava} is essentially that of presenting tools for
analyzing the shapes of absolutely continuous functions and defining a distance metric
for computing a distance between any two of them. Specializing to~$AC[0,1]$, a crucial
aspect of the approach in~\cite{srivastava} is then that of identifying a bijective
correspondence between functions in~$AC[0,1]$ and functions in $L^2[0,1]$, and taking
advantage of this correspondence to compute easily the distance between functions
in~$AC[0,1]$ (the definition of this distance in terms of the so-called Fisher-Rao
metric appears below) by computing the distance between the corresponding functions
in~$L^2[0,1]$. Actually, as mentioned above, the goal of this approach is not so much
that of computing the distance between functions in~$AC[0,1]$ but of computing the
distance between their shapes. More precisely, in this approach, each function
in~$AC[0,1]$ is associated with its unique (a.e. on~$[0,1]$) so-called square-root
slope function (SRSF) in~$L^2[0,1]$, and vice versa, and a distance metric is defined
for computing the distance between the shapes of any two functions in~$AC[0,1]$ in
terms of the $L^2$ distances between SRSF's of reparametrizations of the two
functions. This distance, although computed in~$L^2[0,1]$, is a measure of how
much one of the absolutely continuous functions must be reparametrized (with
so-called warping functions) to align as much as possible with the other one.
Since given two functions in~$AC[0,1]$ that are not equal, the possibility
exists that one function can be reparametrized to align exactly with the other one,
i.e., become exactly the other one, the set of reparametrization functions or warping
functions then induces a quotient space of~$L^2[0,1]$.
Accordingly, a distance metric is defined in~\cite{srivastava}
that computes the distance between any two equivalence classes in the quotient
space of $L^2[0,1]$ by the set of warping functions, thus computing the distance
between the shapes of the two corresponding functions in~$AC[0,1]$.
\\ \smallskip\\
{\bf Definition 4.1 (SRSF representation of functions \cite{srivastava}):} Given
$f\in AC[0,1]$, the real-valued {\bf square-root slope function (SRSF)} $q$ of $f$,
is defined for each $t\in [0,1]$ at which $f'$ exists by
\[ q(t)= \mathrm{sign} (f'(t))\sqrt{|f'(t)|}. \]
{\bf Observation 4.2:} By Corollary 3.5, $f'$ exists a.e. on $[0,1]$. Thus $q$ is
defined a.e. on~$[0,1]$. We note that $q$, the SRSF of~$f$, is the $1-$dimensional
version of the {\bf square-root velocity function (SRVF)} $q$ of an absolutely
continuous function $f$, $f:[0,1]\rightarrow {\bf R}^n$, defined as
follows. Let $F: {\bf R}^n\rightarrow {\bf R}^n$ be the continuous map defined by
$F(v)=v/\sqrt{|v|}$ if $|v|\not=0$, $F(v)=0$ otherwise, $|\cdot|$ the
Euclidean norm. Then the SRVF~$q$ of~$f$, $q:[0,1]\rightarrow {\bf R}^n$, is defined
for each $t\in [0,1]$ at which $f'$ exists by
\[ q(t)= F(f'(t)) = f'(t)/\sqrt{|f'(t)|} \]
if $|f'(t)|\not=0$, $0$ (in ${\bf R}^n$) otherwise.\\
See \cite{lahiri}, \cite{srivastava} for a rigorous development of the SRVF.
\\ \smallskip\\
{\bf Proposition 4.1 (Square integrability of SRSF \cite{srivastava}):} Given
$f\in AC[0,1]$, the SRSF~$q$ of~$f$ is square-integrable
over $[0,1]$, i.e.,~$q\in L^2[0,1]$, with $\int_0^1 |q(t)|^2 dt = \int_0^1 |f'(t)| dt$,
i.e., $||q||_2^2$ = length of~$f$.
\\ \smallskip\\
{\bf Proof:} By Corollary 3.5, $f'$ is measurable and Lebesgue integrable
over~$[0,1]$. Note $h(t) = |q(t)|^2 = |\mathrm{sign} (f'(t))\sqrt{|f'(t)|}|^2 =
|f'(t)|$ for each $t\in [0,1]$ at which $f'$ exists. Thus $h$ is measurable
and  Lebesgue integrable over $[0,1]$ (Proposition~2.31) so that $q\in L^2[0,1]$
and $\int_0^1 |q(t)|^2 dt = \int_0^1 |f'(t)| dt$.
\\ \smallskip\\
{\bf Observation 4.3:} As noted in Observation 2.21, a Lebesgue integrable function over
a measurable set $E$ can be undefined on a subset of $E$ of measure zero. That can be the
case above for functions $f'$ and~$|q|^2$ with $E=[0,1]$ which we know exist a.e. on $[0,1]$.
However, without any loss of generality, in the spirit of Observation~3.8 about the chain
rule, $f'$ and $q$ will eventually be interpreted below as $f'$ and $q$ extended to all
of~$[0,1]$. Finally, note the length of $f$ above is measured in~{\bf R},
a $1-$dimensional space. \\ \smallskip\\
{\bf Proposition 4.2 (Reconstruction of an absolutely continuous function from
its SRSF \cite{srivastava}):} Given $f\in AC[0,1]$, let $q$ be the SRSF of~$f$.
Then for each $t\in [0,1]$ it must be that $f(t) = f(0) +\int_0^t q(s)|q(s)|ds$.
\\ \smallskip\\
{\bf Proof:} Note that for each $s\in [0,1]$ at which $f'$ exists, then
\begin{eqnarray*}
q(s)|q(s)|&=& \mathrm{sign} (f'(s))\sqrt{|f'(s)|}
|\mathrm{sign} (f'(s))\sqrt{|f'(s)|}|\\
 &=& \mathrm{sign} (f'(s))\sqrt{|f'(s)|} \sqrt{|f'(s)|}\\
 &=& \mathrm{sign} (f'(s))|f'(s)| = f'(s).
\end{eqnarray*}
By Proposition 3.12, for each $t\in [0,1]$ it must be that 
$f(t) = f(0) +\int_0^t f'(s)ds$.
Thus, for each $t\in [0,1]$ it must be that $f(t) = f(0) +\int_0^t q(s)|q(s)|ds$.
\\ \smallskip\\
{\bf Proposition 4.3 ($L^2[0,1]$'s equivalence with the set of all SRSF's
\cite{srivastava}):} Let $q$ be in $L^2[0,1]$ and $C$ any real number.
Let $h(t)=q(t)|q(t)|$ for each $t\in [0,1]$ at which $q$ exists.
Then $h$ is defined a.e. on $[0,1]$, $h$ is measurable and Lebesgue integrable
over $[0,1]$, and the function $f$ defined for each $t\in [0,1]$ by
$f(t) = C +\int_0^t h(s)ds$ is absolutely continuous on~$[0,1]$ with $q$ equal
to the SRSF of~$f$ a.e. on~$[0,1]$. \\ \smallskip\\
{\bf Proof:} As $q$ is defined a.e.  on $[0,1]$, then so is $h$. In addition,
since $|h|=|q|^2$ is measurable and Lebesgue integrable over $[0,1]$,
then so is $h$ (Proposition~2.31). By Proposition~3.12, $f$ is then absolutely continuous
on~$[0,1]$.\\ Let $\hat{q}$ be the SRSF of $f$. Then for each $t\in [0,1]$ at which $f'$
exists it must be that $\hat{q}(t) = \mathrm{sign}(f'(t))\sqrt{|f'(t)|}$ and $\hat{q}$ is
defined a.e. on~$[0,1]$. Since by Proposition 3.7, $f'=h$ a.e. on $[0,1]$,
then it must also be that $\hat{q}(t) = \mathrm{sign}(h(t))\sqrt{|h(t)|}$ for
almost all $t\in [0,1]$.\\
But $\mathrm{sign}(h(t))\sqrt{|h(t)|} =
\mathrm{sign}(q(t)|q(t)|)\sqrt{q(t)^2} = \mathrm{sign}(q(t))|q(t)|=q(t)$
for each $t\in [0,1]$ at which $q$ exists and therefore for almost all $t\in [0,1]$.
Thus, $q=\hat{q}$ a.e. on~$[0,1]$. \\ \smallskip\\
{\bf Definition 4.2:} Under the composition of functions operation, the {\bf admissible
class $\Gamma$ of warping functions} is a semigroup of functions (not every element
has an inverse) defined~by
\begin{eqnarray*}
\Gamma &=& \{\gamma| \gamma: [0,1]\rightarrow [0,1],
\ \gamma \ \mathrm{absolutely\ continuous\ on}\ [0,1],\\
&\ &\ \ \dot{\gamma}\geq 0 \mathrm{\ a.e.\ on}\ [0,1],\ \gamma(0)=0,\ \gamma(1)=1\},
\end{eqnarray*}
where $\dot{\gamma}$ is the derivative of~$\gamma$.\\ \smallskip\\
The {\bf group $\Gamma_0$ of invertible warping functions}, $\Gamma_0\subset\Gamma$,
is defined by
\[ \Gamma_0 = \{\gamma|\gamma\in\Gamma,\ \dot{\gamma}>0 \mathrm{\ a.e.\ on\ } [0,1]\}.\]
{\bf Observation 4.4:} The functions in $\Gamma$ and $\Gamma_0$ play an important role in
functional data and shape analysis as they are used to reparametrize an absolutely continuous
function by warping its domain during the process of aligning its shape to the shape of
another absolutely continuous function. As demonstrated in \cite{lahiri}, \cite{srivastava},
it is $\Gamma$ and $\Gamma_0$ that induce a quotient space of~$L^2[0,1]$ with a well-defined
distance metric. More on this below. We note, given $\gamma\in\Gamma$, since $\gamma$ is
continuous, $\gamma(0)=0$, $\gamma(1)=1$, $\gamma:[0,1]\rightarrow [0,1]$, then by the
intermediate value theorem \cite{rudin}, $\gamma([0,1])=[0,1]$.
We note, given $\gamma\in\Gamma_0$, $\gamma$ is
strictly increasing, thus has an inverse $\gamma^{-1}$ which
is also in~$\Gamma_0$ as $\gamma^{-1}(0)=0$, $\gamma^{-1}(1)=1$,
$\gamma^{-1}$ is absolutely continuous on~$[0,1]$ by Proposition~3.24,
and $(\gamma^{-1})'>0$ a.e. on~$[0,1]$, also from Proposition~3.24, since $\gamma$ (the
inverse of $\gamma^{-1}$) is absolutely continuous on~$[0,1]$.\\
Note, given $\gamma$, $\tilde{\gamma}\in\Gamma$, then $\gamma\circ\tilde{\gamma}$ is
absolutely continuous on~$[0,1]$ by Proposition~3.25 and clearly
$(\gamma\circ\tilde{\gamma})(0)=0$, $(\gamma\circ\tilde{\gamma})(1)=1$. Accordingly, if
$\gamma$, $\tilde{\gamma}\in\Gamma_0$, in order to conclude that
$\gamma\circ\tilde{\gamma}\in\Gamma_0$, we prove $(\gamma\circ\tilde{\gamma})'>0$ a.e.
on~$[0,1]$. For this purpose let $A=\{t\in [0,1]:\dot{\gamma}(t)>0\}$, $B=[0,1]\setminus A$,
$C=\tilde{\gamma}^{-1}(B)$. Clearly, $mB=0$ and since $\tilde{\gamma}^{-1}$ is absolutely
continuous on~$[0,1]$ as just proved above, then $mC=0$ by Proposition~3.21. Let
$D=[0,1]\setminus C$. Accordingly, we only need to prove $(\gamma\circ\tilde{\gamma})'>0$
a.e. on~$D$. Clearly, $\dot{\tilde{\gamma}}$ exists (and is positive) a.e. on~$D$.
Since $\tilde{\gamma}(D)\subseteq A$, and $\dot{\gamma}$ exists (and is positive) on~$A$,
then $(\gamma\circ\tilde{\gamma})'$ exists a.e. on~$D$. Indeed it exists exactly at the points
in~$D$ where $\dot{\tilde{\gamma}}$ exists. Thus, by the usual chain rule of calculus,
$(\gamma\circ\tilde{\gamma})'(t) = \dot{\gamma}(\tilde{\gamma}(t))\dot{\tilde{\gamma}}(t)$
for $t\in D$ at which $\dot{\tilde{\gamma}}$ exists. Since as mentioned above
$\dot{\gamma}(\tilde{\gamma}(t))$ exists and is positive for all $t\in D$, and
$\dot{\tilde{\gamma}}$ exists and is positive a.e. on~$D$, then
$(\gamma\circ\tilde{\gamma})'=(\dot{\gamma}\circ\tilde{\gamma})\dot{\tilde{\gamma}}>0$
a.e. on~$D$.\\
Finally, given $\gamma$, $\tilde{\gamma}\in\Gamma$, we show $\gamma\circ\tilde{\gamma}\in\Gamma$.
It suffices to show $(\gamma\circ\tilde{\gamma})'\geq 0$ a.e. on~$[0,1]$.
For this purpose let $A=\{t\in [0,1]:\dot{\gamma}(t) \geq 0\}$, $B=[0,1]\setminus A$,
$C=\{t\in [a,b]:\tilde{\gamma}(t)\in B\}$. Clearly, $mB=0$. Letting $G=[a,b]\setminus C$,
then we can show that $(\gamma\circ\tilde{\gamma})'\geq 0$ a.e. on G, in the same manner
we showed above for $\gamma$, $\tilde{\gamma}\in\Gamma_0$ that $(\gamma\circ\tilde{\gamma})'>0$
a.e. on~$D$. Thus, in order to complete the proof we show $(\gamma\circ\tilde{\gamma})'=0$ a.e.
on~$C$.  Since $\tilde{\gamma}(C)\subseteq B$, then $m(\tilde{\gamma}(C)) = mB = 0$, and since
$\gamma$ is absolutely continuous on~$[0,1]$, then $m(\gamma(\tilde{\gamma}(C)))=0$ by
Proposition~3.21.  Note $(\gamma\circ\tilde{\gamma})'$ exists a.e. on~$C$ as
$\gamma\circ\tilde{\gamma}$ is absolutely continuous on~$[0,1]$.
That $(\gamma\circ\tilde{\gamma})'=0$ a.e. on~$C$ now follows follows from Corollary~3.12.
\\ \smallskip\\
{\bf Proposition 4.4 (SRSF of a warped absolutely continuous function \cite{srivastava}):}
Given $f\in AC[0,1]$ and $\gamma\in\Gamma$, then $f\circ\gamma\in AC[0,1]$ and
$(f\circ\gamma)(0) = f(0)$. With $q$ the SRSF of~$f$, without any loss of generality,
in the spirit of Observation~3.8 about the chain rule, interpreting $f'$ and $q$ as $f'$
and $q$ extended to all of $[0,1]$, it then follows that the SRSF of $f\circ\gamma$ equals
$(q\circ\gamma)\sqrt{\dot{\gamma}}$ a.e. on $[0,1]$.
\\ \smallskip\\
{\bf Proof:} Clearly, $(f\circ\gamma)(0)=f(\gamma(0))=f(0)$. That $f\circ\gamma\in AC[0,1]$
follows directly from Proposition~3.25. Accordingly, it then follows from Corollary~3.13
(chain rule) that $(f\circ \gamma)' = (f'\circ \gamma)\dot{\gamma}$ a.e. on~$[0,1]$. Thus,
the SRSF of $f\circ\gamma$ which is defined for each $t\in[0,1]$ at which $(f\circ\gamma)'$
exists as $\mathrm{sign}((f\circ\gamma)'(t))\sqrt{|(f\circ\gamma)'(t)|}$ must equal
\begin{eqnarray*}
\mathrm{sign}((f'\circ\gamma)(t)\dot{\gamma}(t))
\sqrt{|(f'\circ\gamma)(t)\dot{\gamma}(t)|}
&=&\mathrm{sign}(f'(\gamma(t)))\sqrt{|f'(\gamma(t))|}\sqrt{\dot{\gamma}(t)}\\
 &=& q(\gamma(t))\sqrt{\dot{\gamma}(t)}
\end{eqnarray*}
for almost all $t\in [0,1]$.
\\ \smallskip\\
{\bf Observation 4.5:} Note the SRSF $q$ of $f$ is defined (Definition~4.1) for each
$t\in [0,1]$ at which $f'$ exists. However, although the SRSF of $f\circ\gamma$
equals $(q\circ\gamma)\sqrt{\dot{\gamma}}$ a.e. on~$[0,1]$, it is not necessarily true
that the SRSF of $f\circ\gamma$ exists at each $t\in [0,1]$ for which
$q(\gamma(t))\sqrt{\dot{\gamma}(t)}$ exists or if it exists it is equal to it.
\\ \smallskip\\
{\bf Observation 4.6:} An isometry is a distance-preserving transformation
between two metric spaces. Here we describe in a nonrigorous manner isometries
(and differentials as well) in the context of differential and Riemannian geometry.
Let $M$, $N$ be spaces and let $\varphi$ be a mapping from $M$ into $N$, with $M$,
$N$, $\varphi$ satisfying certain smoothness properties (in the language of
differential geometry, $M$ and $N$ are {\bf smooth} or {\bf differentiable
manifolds} which are spaces that locally resemble Euclidean, Hilbert or Banach
spaces, and $\varphi$ is {\bf differentiable (generalized to smooth manifolds)};
here and in what follows, differentiability in the context of smooth manifolds
is assumed to be of all orders). Given $\epsilon>0$, $p\in M$, assume
$\alpha:(-\epsilon,\epsilon)\rightarrow M$ can be defined, $\alpha(0)=p$,
$\alpha$ a curve in $M$, differentiable (generalized to smooth manifolds)
so that $\alpha'(0)$ makes sense. Then $\alpha'(0)$ is considered to be a
tangent vector to the curve $\alpha$ at $t=0$, and to $M$ at~$p$. Accordingly,
the set of all tangent vectors to $M$ at $p$ is called the {\bf tangent space}
of $M$ at~$p$ and denoted by~$T_pM$. Similarly, given $q\in N$, the set of all
tangent vectors to $N$ at $q$ is called the tangent space of $N$ at $q$ and
denoted by~$T_qN$. With $\alpha$ as above, define
$\beta:(-\epsilon,\epsilon)\rightarrow N$ by $\beta=\varphi\circ\alpha$. Then
$\beta(0)=\varphi(p)$, $\beta$ is a curve in~$N$, and we assume it is
differentiable (generalized to smooth manifolds) so that $\beta'(0)$ makes
sense and is then in $T_{\varphi(p)}N$. The mapping
$d\varphi_p:T_pM\rightarrow T_{\varphi(p)}N$ given by
$d\varphi_p(\alpha'(0)) = \beta'(0)$ is a linear mapping called the
{\bf differential} of~$\varphi$ at~$p$. Finally, assume there is a correspondence
on $M$, smooth in some manner (see below), that associates
to each point $p$ in $M$ an inner product $<,>_p$ on the tangent space $T_pM$.
Similarly, assume there is a correspondence on $N$, smooth in the same manner,
that associates to each point $q$ in $N$ an inner product $<,>_q$ on the tangent
space~$T_qN$. If $\varphi$ as above is bijective and satisfies certain smoothness
properties (in the language of differential geometry, $\varphi$ is a
{\bf diffeomorphism}: $\varphi$ is bijective, and $\varphi$ and $\varphi^{-1}$
are differentiable (generalized to smooth manifolds)), then $\varphi$ is called
an {\bf isometry} if $<u,v>_p = <d\varphi_p(u),d\varphi_p(v)>_{\varphi(p)}$, for
all $p\in M$ and all $u,v\in T_pM$, where the inner product on the left is the one
on $T_pM$ and the inner product on the right is the one on~$T_{\varphi(p)}N$.
On the other hand, if $\varphi$ is differentiable and satisfies
$<u,v>_p = <d\varphi_p(u),d\varphi_p(v)>_{\varphi(p)}$, for all $p\in M$ and all
$u,v\in T_pM$, but is not a diffeomorphism, then $\varphi$ is called
a~{\bf semi-isometry}.\\
In the language of Riemannian geometry, the smooth correspondence above between points
in a smooth manifold and inner products on tangent spaces of the space at
the points is called a {\bf Riemannian metric} or {\bf structure}. Smooth
manifolds equipped with such a structure are called {\bf Riemannian manifolds}.
Using the Riemannian structure, the length of a curve in a Riemannian manifold~$M$
is computed as follows. Given $\alpha:[0,1]\rightarrow M$, a curve or path in~$M$,
differentiable (generalized to smooth manifolds) on~$[0,1]$ so that~$\alpha'(t)$
makes sense for $t\in [0,1]$ and is then in the tangent space~$T_{\alpha(t)}M$,
the length $L(\alpha)$ of the the path $\alpha$ is then given by
\[ L(\alpha)=\int_0^1\sqrt{<\alpha'(t),\alpha'(t)>_{\alpha(t)}}dt,\]
where $<,>_{\alpha(t)}$ is the inner product on the tangent space $T_{\alpha(t)}M$,
and the smoothness of the Riemannian structure is such that
$\sqrt{<\alpha'(t),\alpha'(t)>_{\alpha(t)}}$, $t\in [0,1]$, is integrable
over~$[0,1]$ so that $L(\alpha)$ is well defined. In addition, given $p$, $q\in M$,
the {\bf geodesic distance} $d(p,q)$ between them is defined as the
minimum of the lengths of all paths $\alpha$ in $M$, $\alpha:[0,1]\rightarrow M$,
differentiable (generalized to smooth manifolds)
with $\alpha(0)=p$ and $\alpha(1)=q$, i.e.,
\[ d(p,q) = \min_{\alpha:[0,1]\rightarrow M,\,\alpha
\,\mathrm{differentiable\,on\,}[0,1],\,\alpha(0)\,=\,p,\,\alpha(1)\,=\,q}
L(\alpha).\]
If a path $\alpha$ exists such that $d(p,q)$ achieves its minimum at $\alpha$,
then $\alpha$ is called a {\bf geodesic} in~$M$ between~$p$ and~$q$.
We note that geodesics in Euclidean spaces and $L^2[0,1]$ are given
by straight lines.
Thus, for example, given $p$, $q$ in $L^2[0,1]$, then
$\alpha:[0,1]\rightarrow L^2[0,1]$ defined by $\alpha(t)=(1-t)p+tq$ for $t\in [0,1]$,
is the geodesic between $p$ and $q$ and the distance $d(p,q)$ is
\[ \int_0^1 (\int_0^1 |p(s)-q(s)|^2ds)^{1/2}dt = \int_0^1 ||p-q||_2dt =
||p-q||_2.\]
Finally, we note that with the distance $d(p,q)$ as defined
above, it then follows that an isometry also as defined above is indeed
a distance-preserving transformation. We show this in a nonrigorous manner.
Given Riemannian manifolds~$M$, $N$, let $\alpha:[0,1]\rightarrow M$ be a path
from $p$ to $q$ in $M$, $\alpha(0)=p$, $\alpha(1)=q$,
$\varphi:M\rightarrow N$ an isometry. Let $\beta = \varphi\circ\alpha$. Then
$\beta:[0,1]\rightarrow N$,
$\beta(0)=\varphi(p)$, $\beta(1)=\varphi(q)$, and $\beta$ is a path from
$\varphi(p)$ to $\varphi(q)$ in~$N$. Since for any $t\in [0,1]$ the differential
$d\varphi_{\alpha(t)}:T_{\alpha(t)}M\rightarrow T_{\varphi(\alpha(t))}N$ is given by
$d\varphi_{\alpha(t)}(\alpha'(t)) = \beta'(t)$, then 
\begin{eqnarray*}
L(\beta)&=&\int_0^1\sqrt{<\beta'(t),\beta'(t)>_{\beta(t)}}dt\\
&=& \int_0^1\sqrt{<d\varphi_{\alpha(t)}(\alpha'(t)),
d\varphi_{\alpha(t)}(\alpha'(t))>_{\varphi(\alpha(t))}}dt\\
&=&\int_0^1\sqrt{<\alpha'(t),\alpha'(t)>_{\alpha(t)}}dt=L(\alpha)
\end{eqnarray*}
since $\varphi$ is an isometry. Similarly, given a path $\beta:[0,1]\rightarrow N$
from $\varphi(p)$ to $\varphi(q)$ in~$N$, there is a path $\alpha:[0,1]\rightarrow M$
from $p$ to~$q$ in~$M$ with $L(\alpha)=L(\beta)$. Thus, $d(p,q)=d(\varphi(p),\varphi(q))$.
\\See \cite{carmo1}, \cite{carmo2}, \cite{lee1}, \cite{lee2}, \cite{srivastava}
for a more rigorous development of the concepts of smooth manifolds, Riemannian
manifolds, differentials, isometries,~etc. \\ \smallskip\\
{\bf Observation 4.7:} In what follows, given $q\in L^2[0,1]$, $\gamma\in\Gamma$,
we use $(q,\gamma)$~as short notation for $(q\circ\gamma)\sqrt{\dot{\gamma}}$.
Here again, without any loss of generality, in the spirit of Observation~3.8 about the
chain rule, $q$ is interpreted as $q$ extended to all of $[0,1]$. As it will be shown
below, $(q,\gamma)\in L^2[0,1]$ so that without any loss of generality, again in the
spirit of Observation~3.8 about the chain rule, $(q,\gamma)$ can be interpreted as
$(q,\gamma)$ extended to all of~$[0,1]$, and given $\overline{\gamma}\in\Gamma$,
$((q,\gamma),\overline{\gamma})$ can be interpreted as $((q,\gamma),\overline{\gamma})$
extended to all of~$[0,1]$.
\\ \smallskip\\
{\bf Proposition 4.5:} Given $q\in L^2[0,1]$ and $\gamma\in\Gamma$, then
$(q,\gamma)\in L^2[0,1]$. In addition, given $\overline{\gamma}\in\Gamma$, then
$((q,\gamma),\overline{\gamma})=(q,\gamma\circ\overline{\gamma})$ a.e. on~$[0,1]$,
and if $\gamma\in\Gamma_0$, then $((q,\gamma),\gamma^{-1}) = q$ a.e. on~$[0,1]$.
\\ \smallskip\\
{\bf Proof:} Let $h(t)=|q(t)|^2$ for each $t\in[0,1]$ at which $q$ exists.
Then $h$ is Lebesgue integrable over~$[0,1]$ and by Corollary~3.14,
$(h\circ\gamma)\dot{\gamma}= |q\circ\gamma|^2\dot{\gamma}$ is Lebesgue integrable
over~$[0,1]$. Thus, $(q,\gamma)= (q\circ\gamma)\sqrt{\dot{\gamma}}\in L^2[0,1]$.\\
Now, if $\overline{\gamma}\in\Gamma$, then
\begin{eqnarray*}
((q,\gamma),\overline{\gamma})&=&((q\circ\gamma)\sqrt{\dot{\gamma}},\overline{\gamma})
=(((q\circ\gamma)\sqrt{\dot{\gamma}})\circ\overline{\gamma})
\sqrt{\dot{\overline{\gamma}}}\\
&=& ((q\circ\gamma\circ\overline{\gamma})\sqrt{\dot{\gamma}\circ\overline{\gamma}})
\sqrt{\dot{\overline{\gamma}}}
= (q\circ(\gamma\circ\overline{\gamma}))\sqrt{(\dot{\gamma}\circ\overline{\gamma})
\dot{\overline{\gamma}}}\\
&=& (q\circ(\gamma\circ\overline{\gamma}))\sqrt{(\gamma\circ\overline{\gamma})'}
=(q,\gamma\circ\overline{\gamma})
\end{eqnarray*}
a.e. on $[0,1]$ using Corollary 3.13 (chain rule) as $\gamma\circ\overline{\gamma}\in
AC[0,1]$, by interpreting $\dot{\gamma}$
as $\dot{\gamma}$ extended to all of~$[0,1]$ (Observation~3.8 about the chain rule).\\
Finally, if $\gamma\in\Gamma_0$, then
\begin{eqnarray*}
((q,\gamma),\gamma^{-1}) &=& ((q\circ\gamma)\sqrt{\dot{\gamma}},\gamma^{-1})
=(((q\circ\gamma)\sqrt{\dot{\gamma}})\circ\gamma^{-1})
\sqrt{(\gamma^{-1})'}\\
&=& ((q\circ\gamma\circ\gamma^{-1})\sqrt{\dot{\gamma}\circ\gamma^{-1}})
\sqrt{(\gamma^{-1})'}
= q\sqrt{(\dot{\gamma}\circ\gamma^{-1})(\gamma^{-1})'}\\
&=& q\sqrt{(\dot{\gamma}\circ\gamma^{-1})/(\dot{\gamma}\circ\gamma^{-1})}=q
\end{eqnarray*}
a.e. on $[0,1]$ using Proposition 3.24.
\\ \smallskip\\
{\bf Definition 4.3:} The {\bf action of $\Gamma_0$ on $L^2[0,1]$} is
the operation that takes any element $\gamma\in\Gamma_0$ and any element $q$ of
$L^2[0,1]$, and computes $(q,\gamma)=(q\circ\gamma)\sqrt{\dot{\gamma}}$.
The {\bf action of $\Gamma$ on $L^2[0,1]$} is similarly defined.
\\ \smallskip\\
{\bf Proposition 4.6 (Action of $\Gamma$ on $L^2[0,1]$ is by semi-isometries.
Action of $\Gamma_0$ on $L^2[0,1]$ is by isometries \cite{srivastava}):}
For each $\gamma\in\Gamma$, let $\varphi^{\gamma}: L^2[0,1]\rightarrow L^2[0,1]$
be defined for $q\in L^2[0,1]$ by
$\varphi^{\gamma}(q)=(q,\gamma)=(q\circ\gamma)\sqrt{\dot{\gamma}}$. Then
$\varphi^{\gamma}$ is differentiable and
\[ <d\varphi^{\gamma}(u),d\varphi^{\gamma}(v)> = <(u,\gamma),(v,\gamma)> = <u,v> \]
for any $u$, $v\in L^2[0,1]$,
where $<,>$ is the $L^2[0,1]$ inner product and $d\varphi^{\gamma}$ is the
differential of $\varphi^{\gamma}$, with $<,>$ and $d\varphi^{\gamma}$ the same at
every $q\in L^2[0,1]$. Thus, $\varphi^{\gamma}$ is a semi-isometry and the action
of $\Gamma$ on $L^2[0,1]$ is said to be by semi-isometries. If $\gamma\in\Gamma_0$,
then $\varphi^{\gamma}$ is a diffeomorphism. Thus, $\varphi^{\gamma}$ is an isometry
and the action of $\Gamma_0$ on $L^2[0,1]$ is said to be by isometries.
\\ \smallskip\\
{\bf Proof:} If $\gamma\in\Gamma$, from Proposition 4.5 it follows that the range
of $\varphi^{\gamma}$ is indeed in~$L^2[0,1]$.
Since the tangent space of $L^2[0,1]$ at any point is $L^2[0,1]$ itself, it follows
that $<,>$ is the same at every~$q\in L^2[0,1]$.
Given $u$, $v\in L^2[0,1]$, from Proposition~2.41 (H\"{o}lder's
inequality), $u\cdot v\in L^1$. Let $h(s)= u(s)v(s)$ for each
$s\in [0,1]$ at which $u$ and $v$ exist. By Corollary~3.14 (change of variable),
\begin{eqnarray*}
<u,v> &=& \int_0^1u(s)v(s)ds
= \int_{\gamma(0)}^{\gamma(1)}h(s)ds =\int_0^1 h(\gamma(t))\dot{\gamma}(t)dt\\
&=& \int_0^1 u(\gamma(t))\sqrt{\dot{\gamma}(t)} v(\gamma(t))\sqrt{\dot{\gamma}(t)}dt\\
&=& <(u,\gamma),(v,\gamma)> = <d\varphi^{\gamma}(u),d\varphi^{\gamma}(v)>
\end{eqnarray*}
as $\varphi^{\gamma}$ is linear so that it is differentiable and $d\varphi^{\gamma}$
acts on an element of $L^2[0,1]$ the same way $\varphi^{\gamma}$~does.
Thus, $d\varphi^{\gamma}$ is the same at every $q\in L^2[0,1]$,
and in addition, $\varphi^{\gamma}$ is a semi-isometry and the action
of $\Gamma$ on $L^2[0,1]$ is by semi-isometries. If $\gamma\in\Gamma_0$,
then $\varphi^{\gamma}$ is a bijection and its inverse~$\varphi^{\gamma^{-1}}$
is linear so that it is differentiable. Thus, $\varphi^{\gamma}$ is a diffeomorphism and
therefore an isometry, and the action of $\Gamma_0$ on $L^2[0,1]$ is by~isometries.
\\ \smallskip\\
{\bf Corollary 4.1 (Action of $\Gamma_0$ on $L^2[0,1]$ is distance preserving
\cite{srivastava}):} Given $q_1$, $q_2\in L^2$, and $\gamma\in\Gamma_0$, then
$||q_1-q_2||_2 = ||(q_1,\gamma)-(q_2,\gamma)||_2$. \\ \smallskip\\
{\bf Proof:} $d(q_1,q_2)= ||q_1-q_2||_2$ and
$d((q_1,\gamma),(q_2,\gamma))=||(q_1,\gamma)-(q_2,\gamma)||_2$ (Observation~4.6).
The action of $\Gamma_0$ on~$L^2[0,1]$ is by isometries (Proposition~4.6). Thus,
$d(q_1,q_2)=d((q_1,\gamma),(q_2,\gamma))$ (Observation 4.6), and hence,
$||q_1-q_2||_2 = ||(q_1,\gamma)-(q_2,\gamma)||_2$. \\ \smallskip\\
{\bf Corollary 4.2 (Action of $\Gamma_0$ on $L^2[0,1]$ is norm preserving
\cite{srivastava}):} Given $q\in L^2$, and $\gamma\in\Gamma_0$, then
$||q||_2 = ||(q,\gamma)||_2$.
\\ \smallskip\\
{\bf Observation 4.8 (Action of $\Gamma$ on $L^2[0,1]$ is distance and norm
preserving):} Corollary 4.1 and Corollary 4.2 can be shown to hold for all of
$\Gamma$ as follows. Given $q_1$, $q_2\in L^2$, and $\gamma\in\Gamma$, then
$q_1-q_2\in L^2$ and  by Corollary~3.14 (change of variable),
\begin{eqnarray*}
||(q_1,\gamma)-(q_2,\gamma)||_2^2 &=&
\int_0^1 |q_1(\gamma(t))\sqrt{\dot{\gamma}(t)}-
q_2(\gamma(t))\sqrt{\dot{\gamma}(t)}|^2dt\\
&=& \int_0^1 (q_1(\gamma(t))-q_2(\gamma(t)))^2\dot{\gamma}(t)dt\\
&=& \int_{\gamma(0)}^{\gamma(1)} (q_1(s)-q_2(s))^2ds
= \int_0^1 (q_1(s)-q_2(s))^2ds\\
&=& ||q_1-q_2||_2^2.
\end{eqnarray*}
{\bf Definition 4.4:}
Let $AC^0[0,1]=\{f: f\in AC[0,1],\, f'>0 \mathrm{\ a.e.\ on\ } [0,1]\}$. The
{\bf Fisher-Rao metric} at any $f\in AC^0[0,1]$ is defined as the inner product
\[\ll u,v\gg_f\, = \frac{1}{4}\int_0^1 \dot{u}(t)\dot{v}(t)\frac{1}{f'(t)}dt. \]
for any $u, v\in T_fAC^0[0,1]$.\\ \smallskip\\
{\bf Observation 4.9:} The integral in the definition of the Fisher-Rao metric at
$f\in AC^0[0,1]$ is well defined as $u$, $v$ are functions on~$[0,1]$ that are absolutely
continuous \cite{schmeding}, hence $\dot{u}$, $\dot{v}$, $f'$ exist a.e. on~$[0,1]$, $f'>0$
a.e. on~$[0,1]$, thus $\dot{u}/\sqrt{f'}$, $\dot{v}/\sqrt{f'}$ exist a.e. on~$[0,1]$
and are in $L^2[0,1]$ (see below), so that $\dot{u}\dot{v}/f'$ is Lebesgue
integrable over~$[0,1]$ by Proposition~2.41 (H\"{o}lder's inequality). In addition,
this metric, as defined at elements of $AC^0[0,1]$, is known to have the behavior
of a Riemannian metric~\cite{srivastava}. In what follows, we assume $AC^0[0,1]$
is endowed with this metric.
\\ \smallskip\\
{\bf Proposition 4.7:} Given $f\in AC^0[0,1]$ and $\gamma\in\Gamma_0$,
then $f\circ\gamma\in AC^0[0,1]$ and $(f\circ\gamma)(0)=f(0)$.\\ \smallskip\\
{\bf Proof:} That $f\circ\gamma\in AC[0,1]$ and $(f\circ\gamma)(0)=f(0)$ was
established in Proposition~4.4.
Accordingly, in order to conclude that $f\circ\gamma\in AC^0[0,1]$ we prove
$(f\circ\gamma)'>0$ a.e. on~$[0,1]$. For this purpose let
$A=\{t\in [0,1]:f'(t)>0\}$, $B=[0,1]\setminus A$, $C=\gamma^{-1}(B)$. Clearly,
$mB=0$ and since $\gamma^{-1}$ is absolutely continuous on~$[0,1]$ (Observation~4.4),
then $mC=0$ by Proposition~3.21. Let $D=[0,1]\setminus C$. Accordingly, we only need to
prove $(f\circ\gamma)'>0$ a.e. on~$D$. Clearly, $\dot{\gamma}$ exists (and is positive)
a.e. on~$D$. Since $\gamma(D)\subseteq A$, and $f'$ exists (and is positive) on~$A$,
then $(f\circ\gamma)'$ exists a.e. on~$D$. Indeed it exists exactly at the points
in~$D$ where $\dot{\gamma}$ exists. Thus, by the usual chain rule of calculus,
$(f\circ\gamma)'(t) = f'(\gamma(t))\dot{\gamma}(t)$ for $t\in D$ at which $\dot{\gamma}$
exists. Since as mentioned above $f'(\gamma(t))$ exists and is positive for all $t\in D$,
and $\dot{\gamma}$ exists and is positive a.e. on~$D$, then
$(f\circ\gamma)'=(f'\circ\gamma)\dot{\gamma}>0$ a.e. on~$D$.
\\ \smallskip\\
{\bf Definition 4.5:} The {\bf action of $\Gamma_0$ on $AC^0[0,1]$} is
the operation that takes any element $\gamma\in\Gamma_0$ and any element $f$ of
$AC^0[0,1]$, and computes~$f\circ\gamma$.\\ \smallskip\\
{\bf Observation 4.10:} In what follows, two functions in $AC[0,1]$ are considered
equal if they differ by a constant. Simpler yet, we assume all functions in $AC[0,1]$ have
the same value at zero. Since by Proposition~4.4, if $f\in AC[0,1]$, $\gamma\in\Gamma$, then
$f\circ\gamma\in AC[0,1]$ and $(f\circ\gamma)(0)=f(0)$, and since in addition
the SRSF of $(f+C)$ is the same for any constant~$C$, the latter is a reasonable
assumption. \\ \smallskip\\
{\bf Proposition 4.8 (Action of $\Gamma_0$ on $AC^0[0,1]$ with Fisher-Rao metric
is by isometries \cite{srivastava}):} For each $\gamma\in\Gamma_0$, let
$\varphi^{\gamma}: AC^0[0,1]\rightarrow AC^0[0,1]$ be defined for
$f\in AC^0[0,1]$ by $\varphi^{\gamma}(f)=f\circ\gamma$. Then $\varphi^{\gamma}$
is a diffeomorphism and
\[ \ll d\varphi_f^{\gamma}(u),d\varphi_f^{\gamma}(v)\gg_{f\circ\gamma}
= \ll u\circ \gamma,v\circ\gamma\gg_{f\circ \gamma}= \ll u,v\gg_f\]
for any $u$, $v\in T_fAC^0[0,1]$, where $\ll ,\gg_f$ is the inner product that defines
the Fisher-Rao metric at~$f$, $\ll ,\gg_{f\circ\gamma}$ is the inner product that
defines it at~$f\circ\gamma$, and $d\varphi_f^{\gamma}$ is the differential of
$\varphi^{\gamma}$ at~$f$.  Thus, $\varphi^{\gamma}$ is an isometry
and the action of $\Gamma_0$ on $AC^0[0,1]$ is said to be by~isometries.
\\ \smallskip\\
{\bf Proof:} If $\gamma\in\Gamma_0$, from Proposition 4.7, the range of $\varphi^{\gamma}$
is indeed in~$AC^0[0,1]$. As noted in Observation~4.9, given $f\in AC^0[0,1]$, $u$,
$v\in T_fAC^0[0,1]$, then $u$, $v$ are functions on~$[0,1]$ that are absolutely continuous
\cite{schmeding}, and $h=\dot{u}\dot{v}/f'$ is Lebesgue integrable over~$[0,1]$.  By
Corollary~3.14 (change of variable), Corollary~3.13 (chain rule) as $u\circ\gamma$,
$v\circ\gamma$, $f\circ\gamma\in AC[0,1]$, by interpreting $h$, $\dot{u}$, $\dot{v}$, $f'$
as $h$, $\dot{u}$, $\dot{v}$, $f'$ extended to all of~$[0,1]$ (Observation~3.8 about the
chain rule), and noting that all denominators below are greater than zero a.e. on~$[0,1]$
(Proposition~4.7 and its proof), then
\begin{eqnarray*}
\ll u,v\gg_f &=& \frac{1}{4}\int_0^1\dot{u}(s)\dot{v}(s)\frac{1}{f'(s)}ds
= \frac{1}{4}\int_{\gamma(0)}^{\gamma(1)}h(s)ds\\
&=& \frac{1}{4}\int_0^1 h(\gamma(t))\dot{\gamma}(t)dt
= \frac{1}{4}\int_0^1 \dot{u}(\gamma(t))\dot{v}(\gamma(t))\frac{1}{f'(\gamma(t))}
\dot{\gamma}(t)dt\\
&=&\frac{1}{4}\int_0^1 \dot{u}(\gamma(t))\dot{\gamma}(t)\dot{v}(\gamma(t))\dot{\gamma}(t)
\frac{1}{f'(\gamma(t))\dot{\gamma}(t)}dt\\
&=&\frac{1}{4}\int_0^1 (u(\gamma(t)))'(v(\gamma(t)))'\frac{1}{(f(\gamma(t)))'}dt\\
&=& \ll u\circ\gamma,v\circ\gamma\gg_{f\circ\gamma}
= \ll d\varphi^{\gamma}_f(u),d\varphi^{\gamma}_f(v)\gg_{f\circ\gamma}
\end{eqnarray*}
as $\varphi^{\gamma}$ is linear so that it is differentiable and $d\varphi^{\gamma}_f$
acts on an element of $T_fAC^0[0,1]$ the same way $\varphi^{\gamma}$~does on an
element of~$AC^0[0,1]$. Since $\gamma\in\Gamma_0$, then $\varphi^{\gamma}$ is
a bijection and its inverse~$\varphi^{\gamma^{-1}}$ is linear so that it is differentiable.
Thus, $\varphi^{\gamma}$ is a diffeomorphism and therefore an isometry, and the action
of $\Gamma_0$ on $AC^0[0,1]$ is by~isometries.\\ \smallskip\\
{\bf Proposition 4.9 (Fisher-Rao metric on $AC^0[0,1]$ under SRSF representation becomes
$L^2[0,1]$ metric \cite{srivastava}):} Given $f\in AC^0[0,1]$ and $q$ the SRSF of~$f$,
define a mapping $F:AC^0[0,1]\rightarrow L^2[0,1]$ by~$F(f)=q=\sqrt{f'}$. Then $F$ is
differentiable, and for any $v \in T_fAC^0[0,1]$, it must be that
$F_f^*(v)=\dot{v}/(2\sqrt{f'})\in T_qL^2[0,1]= L^2[0,1]$,
where $F_f^*$ is the differential of~$F$ at~$f$. Given $v_1$, $v_2\in T_fAC^0[0,1]$,
then $<F_f^*(v_1),F_f^*(v_2)>=\ll v_1,v_2\gg_f$, where $<,>$ is
the $L^2[0,1]$ inner product and $\ll ,\gg_f$ is the inner product that defines the
Fisher-Rao metric at~$f$. \\ \smallskip\\
{\bf Proof:} Let
$L_0^1[0,1]=\{\hat{f}: \hat{f}\in L^1[0,1],\, \hat{f}>0 \mathrm{\ a.e.\ on\ } [0,1]\}$.\\
Given $\hat{f}\in L_0^1[0,1]$, define a mapping $S: L_0^1[0,1]\rightarrow L^2[0,1]$ by
$S(\hat{f}) = \sqrt{\hat{f}}$.\\
In addition, given $f\in AC^0[0,1]$, define a mapping (the derivative mapping)
$D:AC^0[0,1]\rightarrow L_0^1[0,1]$ by $D(f)=f'$.\\
With $F$ as defined above, then $F=S\circ D$.\\
Given $v\in T_fAC^0[0,1]$, then $D_f^*(v)=\dot{v}\in T_{f'}L_0^1[0,1]$, where $D_f^*$
is the differential of $D$ at~$f$, as $D$ is linear so that it is differentiable and
$D_f^*$ acts on an element of $T_fAC^0[0,1]$ the same way $D$ acts on an element
of~$AC^0[0,1]$.\\
Let $s:{\bf R}\rightarrow {\bf R}$ be the mapping defined by $s(x)=\sqrt{x}$,
$x\in {\bf R}$, $x>0$. Then $s$ is differentiable for $x>0$, and
$s^*(y)=s'(x)y=y/(2\sqrt{x})$ for any $y\in {\bf R}$, where $s^*$ is the differential
of~$s$. From this, following closely the definition of the differential of a differentiable
function \cite{lee1,srivastava}, it then follows that $S$ is differentiable and
given $w\in T_{\hat{f}}L_0^1[0,1]$, then
$S_{\hat{f}}^*(w) = w/(2\sqrt{\hat{f}})\in T_{\sqrt{\hat{f}}}L^2[0,1]$, where
$S_{\hat{f}}^*$ is the differential of $S$ at~$\hat{f}$.\\
Thus, $F=S\circ D$ is differentiable and its differential
$F_f^*:T_fAC^0[0,1]\rightarrow T_qL^2[0,1]$ at $f\in AC^0[0,1]$ is
$S_{D(f)}^*\circ D_f^*=S_{f'}^*\circ D_f^*$~\cite{lee1}.\\
Accordingly, given $v\in T_fAC^0[0,1]$, then
$F_f^*(v)=S_{f'}^*(D_f^*(v))=S_{f'}^*(\dot{v})=\dot{v}/(2\sqrt{f'})\in
T_{\sqrt{f'}}L^2[0,1]=T_qL^2[0,1]=L^2[0,1]$.\\
Finally, given $v_1$, $v_2\in T_fAC^0[0,1]$, then
\begin{eqnarray*}
<F_f^*(v_1),F_f^*(v_2)>&=&<\dot{v}_1/(2\sqrt{f'}),\dot{v}_2/(2\sqrt{f'})>
= \frac{1}{4}\int_0^1 \dot{v}_1(t)\dot{v}_2(t)\frac{1}{f'(t)}dt\\
&=&\ll v_1,v_2 \gg_f.
\end{eqnarray*}
{\bf Observation 4.11 (Distance between functions in $AC[0,1]$):} Given
$f_1$, $f_2 \in AC[0,1]$, let $q_1$, $q_2$ be the SRSF's of $f_1$, $f_2$, respectively.
We note that computing the distance between $f_1$ and $f_2$ with the Fisher-Rao metric
as defined above may not be possible as a path in $AC[0,1]$ from $f_1$ to $f_2$ might
contain functions whose derivatives are not positive a.e. on~$[0,1]$. Even if this was
not the case, the minimization involved would be nontrivial. Accordingly, motivated by
Proposition~4.9 above, the convention is to say that the {\bf Fisher-Rao
distance} between $f_1$ and $f_2$ is $d_{FR}(f_1,f_2)=||q_1-q_2||_2$, i.e.,
the $L^2$ distance between $q_1$ and $q_2$. In addition, since the geodesic from $q_1$
to $q_2$ is a straight line, given $s\in [0,1]$, then $q=(1-s)q_1+sq_2$ is a function in
this geodesic, and by Proposition~4.3, a function $f\in AC[0,1]$ can be computed for each
$t\in [0,1]$ by $f(t) = C +\int_0^t q(x)|q(x)|dx$, where~$C=f_1(0)=f_2(0)$, with the
SRSF of $f$ equal to $q$ a.e. on~$[0,1]$. Doing this for enough functions on the straight
line joining $q_1$ and $q_2$, a collection of functions can be obtained in~$AC[0,1]$ that
are then said to approximate a geodesic (based on the Fisher-Rao metric) from $f_1$ to~$f_2$. 
\\ \smallskip\\
{\bf Definition 4.6:} Given $q\in L^2[0,1]$, define the {\bf orbit $[q]_{\Gamma_0}$
of $q$ under $\Gamma_0$} by $[q]_{\Gamma_0}=\{\overline{q}:
\overline{q}=(q,\gamma)=(q\circ\gamma)\sqrt{\dot{\gamma}}\ \mathrm{a.e.\ on\ }$[0,1]$,
\ \mathrm{some}\ \gamma\in\Gamma_0\}$. Denote by $cl([q]_{\Gamma_0})$ the closure in $L^2[0,1]$
of~$[q]_{\Gamma_0}$. \\ \smallskip\\
{\bf Observation 4.12:} In what follows, given $q_1$, $q_2\in L^2[0,1]$,
$q_1\in [q_2]_{\Gamma_0}$ so that $q_1= (q_2,\gamma)$ a.e. on~$[0,1]$ for some
$\gamma\in\Gamma_0$, without any loss of generality we may simply say $q_1=(q_2,\gamma)$.
Accordingly, given $q_1$, $q_2\in L^2[0,1]$, $q_1\in [q_2]_{\Gamma_0}$ so that
$q_1= (q_2,\gamma)$ for some $\gamma\in\Gamma_0$, then it follows (Proposition~4.5)
that $[q_1]_{\Gamma_0}\subseteq [q_2]_{\Gamma_0}$, and $q_2=(q_1,\gamma^{-1})$ so that
$[q_2]_{\Gamma_0}\subseteq [q_1]_{\Gamma_0}$ and thus $[q_1]_{\Gamma_0}=[q_2]_{\Gamma_0}$.
Using similar arguments,
given $q_1$, $q_2\in L^2[0,1]$, an equivalence relation $\,\sim\,$ can be defined and
justified on~$L^2[0,1]$ for which $q_1\,\sim\,q_2$ if $q_1$ and $q_2$ are in the same orbit
under~$\Gamma_0$. Accordingly, with this equivalence relation a quotient space is obtained
which is the set of all orbits of elements of $L^2[0,1]$ under $\Gamma_0$ and which we
denote by~$L^2[0,1]/\Gamma_0$. An attempt then can be made as follows to define a distance
function~$d$ between elements of $L^2[0,1]/\Gamma_0$ that would make $L^2[0,1]/\Gamma_0$
a metric space. Given $q_1$, $q_2\in L^2[0,1]$, let
\begin{eqnarray*}
 d([q_1]_{\Gamma_0},[q_2]_{\Gamma_0}) &=& \inf_{\gamma_1,\gamma_2\in\Gamma_0}
||(q_1,\gamma_1)-(q_2,\gamma_2)||_2\\
&=& \inf_{\gamma\in\Gamma_0} ||q_1-(q_2,\gamma)||_2 =
\inf_{\gamma\in\Gamma_0} ||(q_1,\gamma)-q_2||_2,
\end{eqnarray*}
where the bottom equations follow from Corollary~4.1 (action of $\Gamma_0$ on $L^2[0,1]$
is distance preserving) again using Proposition~4.5 where appropriate.
Of the properties that $d$ must satisfy to be a distance function all have been
established \cite{srivastava} except one: $d([q_1]_{\Gamma_0},[q_2]_{\Gamma_0})=0$ if and
only if~$[q_1]_{\Gamma_0}=[q_2]_{\Gamma_0}$. Unfortunately, as demonstrated in \cite{lahiri},
the orbits as defined are not closed in $L^2[0,1]$, which allows for examples with
$[q_1]_{\Gamma_0}\not=[q_2]_{\Gamma_0}$ but $d([q_1]_{\Gamma_0},[q_2]_{\Gamma_0})=0$.
\\ \smallskip\\
{\bf Proposition 4.10 (\cite{lahiri}):} Given $q_1$, $q_2\in L^2[0,1]$, then 
$d([q_1]_{\Gamma_0},[q_2]_{\Gamma_0})=0$ if and only
if~$cl([q_1]_{\Gamma_0})=cl([q_2]_{\Gamma_0})$. In particular, if
$q_1\in cl([q_2]_{\Gamma_0})$ so that $d([q_1]_{\Gamma_0},[q_2]_{\Gamma_0})=0$,
then $cl([q_1]_{\Gamma_0})=cl([q_2]_{\Gamma_0})$. \\ \smallskip\\
{\bf Proof:} If $d([q_1]_{\Gamma_0},[q_2]_{\Gamma_0})=0$, fix
$\overline{q}_1\in [q_1]_{\Gamma_0}$ and note $[\overline{q}_1]_{\Gamma_0}=[q_1]_{\Gamma_0}$
(Observation~4.12) so that
$d([\overline{q}_1]_{\Gamma_0},[q_2]_{\Gamma_0}) = d([q_1]_{\Gamma_0},[q_2]_{\Gamma_0})=0$.
Then~given integer $n>0$, there is $\gamma_n\in\Gamma_0$ such that
$||\overline{q}_1-(q_2,\gamma_n)||_2<1/n$. Thus, $\overline{q}_1\in cl([q_2]_{\Gamma_0})$.
Since $\overline{q}_1$ is arbitrary in $[q_1]_{\Gamma_0}$ then
$[q_1]_{\Gamma_0}\subseteq cl([q_2]_{\Gamma_0})$, thus
$cl([q_1]_{\Gamma_0})\subseteq cl([q_2]_{\Gamma_0})$. Similarly, 
$cl([q_2]_{\Gamma_0})\subseteq cl([q_1]_{\Gamma_0})$, thus
\mbox{$cl([q_1]_{\Gamma_0})=cl([q_2]_{\Gamma_0})$}.\\
Assume $cl([q_1]_{\Gamma_0})=cl([q_2]_{\Gamma_0})$.
Then, in particular, $q_1\in cl([q_2]_{\Gamma_0})$ so that given integer $n>0$,
there is $\gamma_n\in\Gamma_0$ with $||q_1-(q_2,\gamma_n)||_2<1/n$. Thus,
$d([q_1]_{\Gamma_0},[q_2]_{\Gamma_0}) = \inf_{\gamma\in\Gamma_0} ||q_1-(q_2,\gamma)||_2=0$.
\\ \smallskip\\
{\bf Observation 4.13:} Using arguments similar to those in the proof of Proposition~4.10
above, given $q_1$, $q_2\in L^2[0,1]$, an equivalence relation $\,\sim\,$ can be defined and
justified on~$L^2[0,1]$ for which $q_1\,\sim\,q_2$ if $q_1$ and $q_2$ are in the closure
of the same orbit under~$\Gamma_0$. Accordingly, with this equivalence relation a quotient
space is obtained which is the set of all closures of orbits of elements of $L^2[0,1]$ under
$\Gamma_0$ and which we denote by~$L^2[0,1]/\sim\,$. In what follows, we extend the
function~$d$ above to the quotient space~$L^2[0,1]/\sim\,$.
\\ \smallskip\\
{\bf Corollary 4.3 (Distance between equivalence classes in $L^2[0,1]/\sim\,$ \cite{lahiri},
\cite{srivastava}):} Given $q_1$, $q_2\in L^2[0,1]$, let
\[d(cl([q_1]_{\Gamma_0}),cl([q_2]_{\Gamma_0})) =
\inf_{\overline{q}_1\in cl([q_1]_{\Gamma_0}),\overline{q}_2\in cl([q_2]_{\Gamma_0})}
||\overline{q}_1-\overline{q}_2||_2.\]
Then $d(cl([q_1]_{\Gamma_0}),cl([q_2]_{\Gamma_0})) =
\inf_{\gamma_1,\gamma_2\in\Gamma_0} ||(q_1,\gamma_1)-(q_2,\gamma_2)||_2=
\inf_{\gamma\in\Gamma_0} ||q_1-(q_2,\gamma)||_2 =
\inf_{\gamma\in\Gamma_0} ||(q_1,\gamma)-q_2||_2 =
d([q_1]_{\Gamma_0},[q_2]_{\Gamma_0})$, and
$d$ is a distance function between elements of $L^2[0,1]/\sim\,$, so that $L^2[0,1]/\sim\,$
is a metric space with this distance function.\\ \smallskip\\
{\bf Proof:} Note
\[\inf_{\overline{q}_1\in cl([q_1]_{\Gamma_0}),\overline{q}_2\in cl([q_2]_{\Gamma_0})}
||\overline{q}_1-\overline{q}_2||_2 = 
\inf_{\gamma_1,\gamma_2\in\Gamma_0} ||(q_1,\gamma_1)-(q_2,\gamma_2)||_2.\]
Thus, $d(cl([q_1]_{\Gamma_0}),cl([q_2]_{\Gamma_0}))
= \inf_{\gamma_1,\gamma_2\in\Gamma_0} ||(q_1,\gamma_1)-(q_2,\gamma_2)||_2
=\inf_{\gamma\in\Gamma_0} ||q_1-(q_2,\gamma)||_2
=\inf_{\gamma\in\Gamma_0} ||(q_1,\gamma)-q_2||_2
= d([q_1]_{\Gamma_0},[q_2]_{\Gamma_0})$,
as previously noted in Observation~4.12.\\
That $d$ is a distance function follows from Proposition 4.10 and results about
properties of this distance function in~\cite{srivastava}.
\\ \smallskip\\
{\bf Observation 4.14:} Given $f_1$, $f_2\in AC[0,1]$, and $q_1$, $q_2$, the SRSF's
of $f_1$, $f_2$, respectively, we note that $q_1$ and $q_2$ remain unchanged after translations
of $f_1$ and $f_2$ (by translations we mean $f_1$ and $f_2$ become $f_1+c_1$ and $f_2+c_2$,
respectively, for constants $c_1$, $c_2$) so that the distance between the equivalence classes
of $q_1$ and $q_2$, defined by $d$ above, is the same before and after the translations. That this
is true follows from the definition of the SRSF. For scalar multiplications of $f_1$ and $f_2$,
the distance between the equivalence classes of $q_1$ and $q_2$ before and after the scalar
multiplications can be approximated or computed exactly, if possible, by the same elements of
$\Gamma_0$ as the following proposition shows. Accordingly, it is customary to normalize $q_1$
and $q_2$ so that $||q_1||_2=||q_2||_2=1$ and then compute the distance between their equivalence
classes with $d$ as above, as from the comments just made doing so is compatible with the
requirement that the shapes of $f_1$ and $f_2$ be invariant under translation and scalar
multiplication.
\\ \smallskip\\
{\bf Proposition 4.11 (\cite{srivastava}):} Given $q_1$, $q_2\in L^2[0,1]$, and $\gamma^*$,
$\gamma\in\Gamma_0$ for which $||q_1-(q_2,\gamma^*)||_2\leq ||q_1-(q_2,\gamma)||_2$, then
$||bq_1-(cq_2,\gamma^*)||_2\leq ||bq_1-(cq_2,\gamma)||_2$, for any $b,\,c,\,bc>0$.
\\ \smallskip\\
{\bf Proof:} With $<,>$ as the $L^2$ inner product, note
\[ ||q_1-(q_2,\gamma^*)||_2^2 = ||q_1||_2^2 -2 <q_1,(q_2,\gamma^*)> +
||(q_2,\gamma^*)||_2^2, \]
and
\[ ||q_1-(q2,\gamma)||_2^2 = ||q_1||_2^2 -2 <q_1,(q_2,\gamma)> +
||(q_2,\gamma)||_2^2.\]
Thus, $||q_1-(q_2,\gamma^*)||_2\leq ||q_1-(q_2,\gamma)||_2$ and
$||q_2||_2=||(q_2,\gamma^*)||_2=||(q_2,\gamma)||_2$, $bc>0$, imply
$-2bc <q_1,(q_2,\gamma^*)>\,\,\,\ \leq\ -2bc <q_1,(q_2,\gamma)>$.\\
Accordingly, since $||cq_2||_2=||(cq_2,\gamma^*)||_2=||(cq_2,\gamma)||_2$, then
\begin{eqnarray*}
||bq_1||_2^2 -2 <bq_1,(cq_2,\gamma^*)> + ||(cq_2,\gamma^*)||_2^2 & &\\
\leq \mathrm{\ \ } ||bq_1||_2^2 -2 <bq_1,(cq_2,\gamma)> + ||(cq_2,\gamma)||_2^2, & &
\end{eqnarray*}
so that $||bq_1-(cq_2,\gamma^*)||_2\leq ||bq_1-(cq_2,\gamma)||_2$.
\\ \smallskip\\
{\bf Observation 4.15:} Figure~1 illustrates an instance of approximately computing
$d(cl([q_1]_{\Gamma_0}),cl([q_2]_{\Gamma_0}))$ as expressed in Corollary~4.3 above.
Here $q_1$, $q_2$ are the SRSF's of functions $f_1$, $f_2$,
respectively, plotted in the lefmost diagram, $f_1$ in red, $f_2$ in blue, $q_1$, $q_2$
normalized so that
$||q_1||_2=||q_2||_2=1$. The distance (about 0.1436) was approximately computed
(in about 154 seconds) with {\em adapt-DP}~\cite{bernal}, a fast linear Dynamic
Programming algorithm. The resulting warping function
$\gamma\in\Gamma_0$ that approximately minimizes $||q_1-(q_2,\gamma)||_2$ is plotted in the
rightmost diagram, and $f_1$ and $f_2\circ\gamma$ are plotted in the middle diagram in which
they appear essentially aligned. The functions $f_1$ and $f_2$ were given in the form of
sets of 19,693 and 19,763 points, respectively, with nonuniform domains in~$[0,1]$.
A copy of {\em adapt-DP} with usage instructions and data files for the same example
in Figure~1 can be obtained using~links:
\verb+https://doi.org/10.18434/T4/1502501+
\verb+http://math.nist.gov/~JBernal+ \verb+/Fast_Dynamic_Programming.zip+
\begin{figure}
\begin{center}
\begin{tabular}{ccc}
\includegraphics[width=0.3\textwidth]{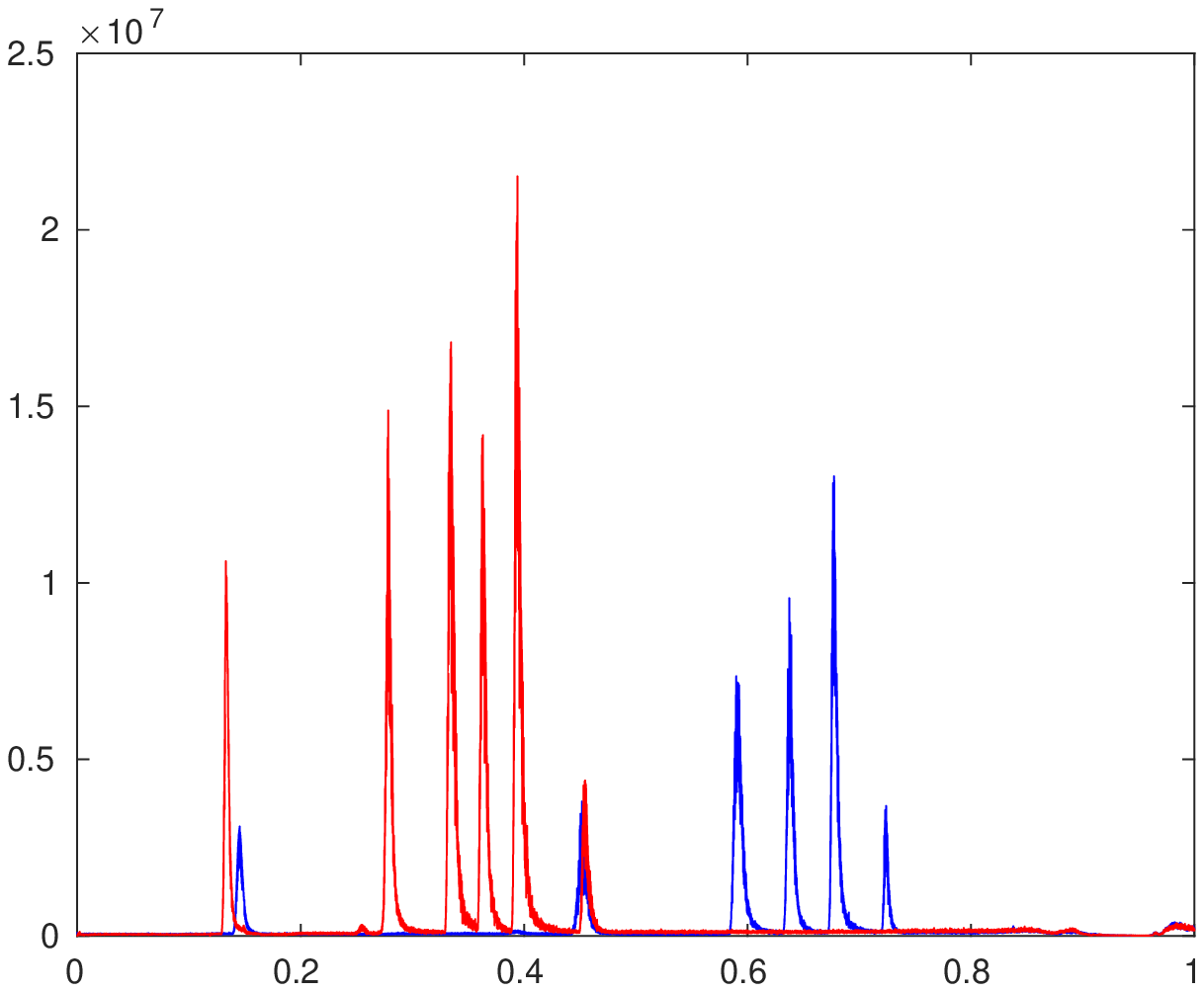}
&
\includegraphics[width=0.3\textwidth]{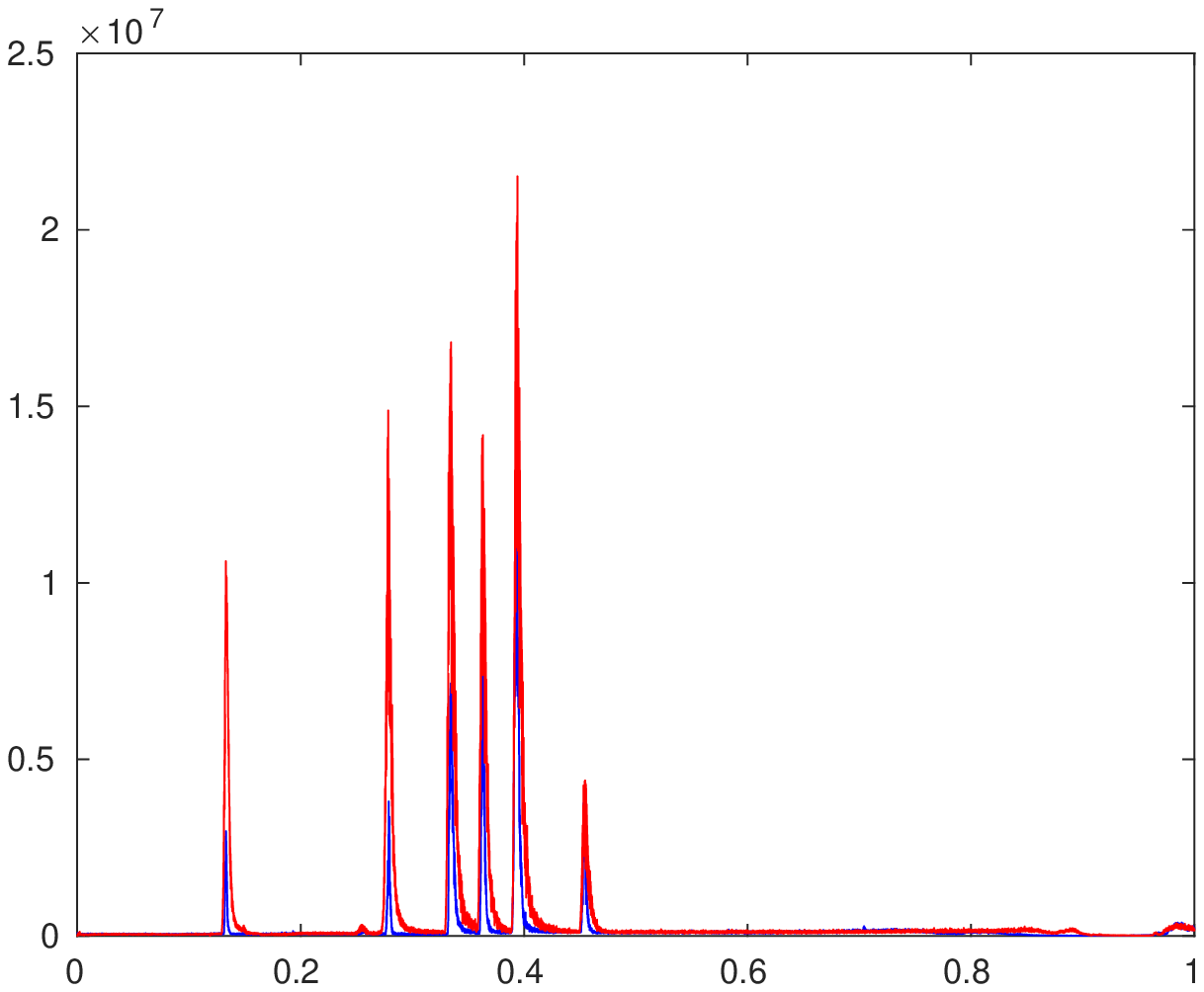}
&
\includegraphics[width=0.3\textwidth]{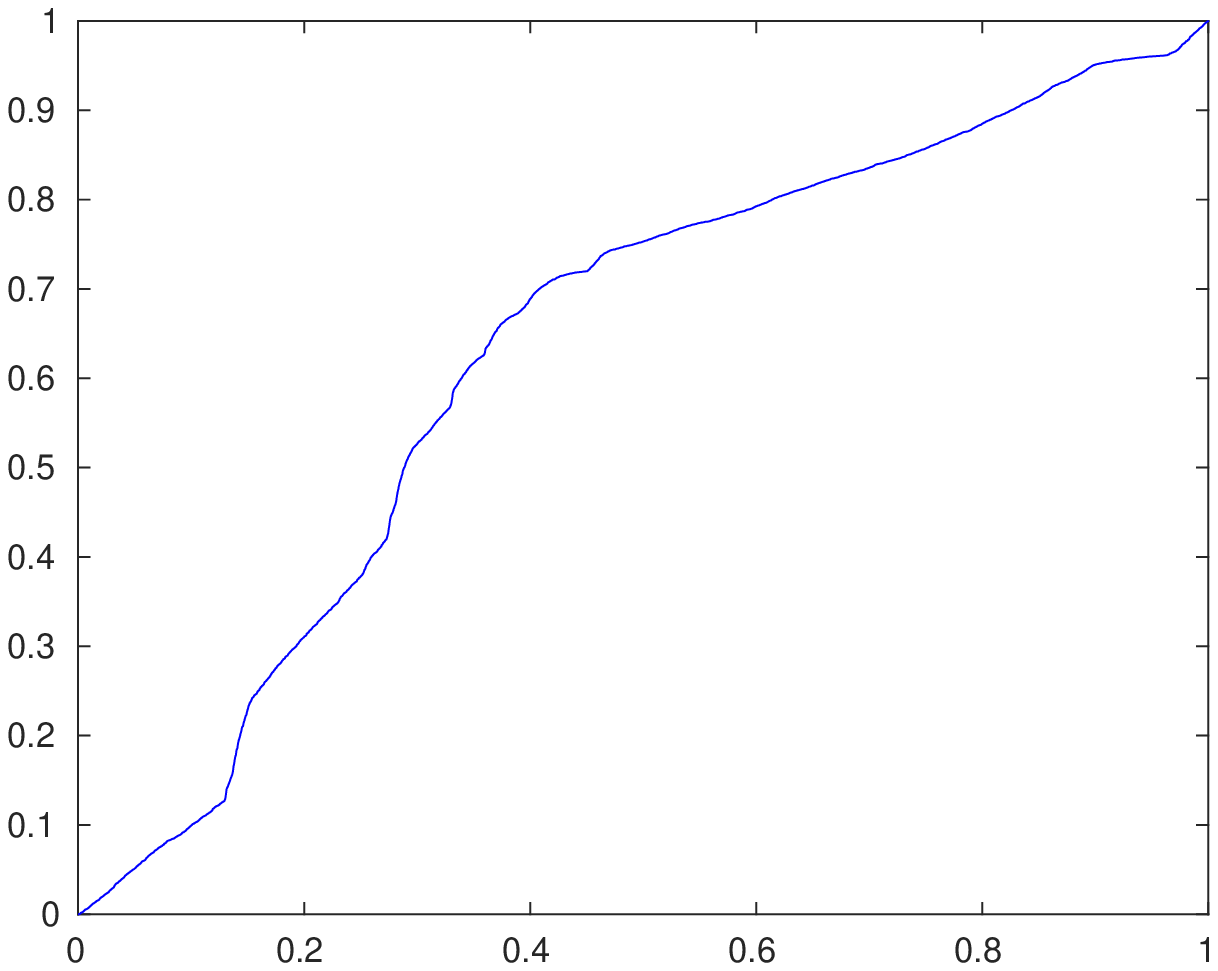}
\end{tabular}
\end{center}
\caption{Function alignment by warping that is obtained from the computation of
$d(cl([q_1]_{\Gamma_0}),cl([q_2]_{\Gamma_0}))$.}
\end{figure}
\\ \smallskip\\
{\bf Definition 4.7:} Given $q\in L^2[0,1]$, define the {\bf orbit $[q]_{\Gamma}$
of $q$ under $\Gamma$} by $[q]_{\Gamma}=\{\overline{q}:
\overline{q}=(q,\gamma)=(q\circ\gamma)\sqrt{\dot{\gamma}}\ \mathrm{a.e.\ on\ }$[0,1]$,
\ \mathrm{some}\ \gamma\in\Gamma\}$.  \\ \smallskip\\
{\bf Observation 4.16:} In what follows, we present results found mostly in~\cite{lahiri}
for the purpose of showing that given $q\in L^2[0,1]$, then there exist $w\in L^2[0,1]$,
$\gamma\in\Gamma$, such that $q=(w,\gamma)$, $|w|$ constant a.e. on~$[0,1]$,
$cl([w]_{\Gamma_0}) = [w]_{\Gamma}$ so that $q\in cl([w]_{\Gamma_0})$, and thus
$cl([q]_{\Gamma_0}) = cl([w]_{\Gamma_0}) = [w]_{\Gamma}$.
We note that this result doesn't change how
$d(cl([q_1]_{\Gamma_0}),cl([q_2]_{\Gamma_0}))$ in Corollary~4.3 is computed for $q_1$,
$q_2\in L^2[0,1]$. It should still be done by computing
$\inf_{\gamma\in\Gamma_0} ||q_1-(q_2,\gamma)||_2$ or
$\inf_{\gamma\in\Gamma_0} ||(q_1,\gamma)-q_2||_2$ as implied by Corollary~4.3.
\\ \smallskip\\
{\bf Proposition 4.12:} $A_0=\{q\in L^2[0,1]: ||q||_2=1, \ q>0\ \mathrm{a.e.\ on}
\ [0,1]\}$ has closure in $L^2[0,1]$ equal to $A=\{q\in L^2[0,1]: ||q||_2=1, \ q\geq 0
\ \mathrm{a.e.\ on}\ [0,1]\}$.
\\ \smallskip\\
{\bf Proof:} Clearly $A_0\subset A$. Let $\epsilon>0$ be given. Given $q\in A\setminus A_0$,
then a measurable subset~$B$ of $[0,1]$ exists, $m(B)>0$, on which~$q=0$.
Let $C=[0,1]\setminus B$. Then $q>0$ a.e. on~$C$ and $\int_C q(t)^2dt = \int_0^1 q(t)^2dt=1$.\\
Choose $b>0$, $b<1$ with $1-b<\epsilon/2$, and set $a=((1-b^2)/m(B))^{1/2}$.\\
Define a function $\hat{q}$ on $[0,1]$ by $\hat{q}=bq$ on~$C$ and $\hat{q}=a$ on~$B$.
Then $\hat{q}>0$ a.e. on~$[0,1]$ and $\int_0^1\hat{q}(t)^2dt=\int_C(bq(t))^2dt +\int_B a^2 dt$
$=b^2\cdot 1 + a^2m(B) = b^2+(1-b^2) = 1$ so that $\hat{q}\in A_0$.\\
Note $\int_0^1(q(t)-\hat{q}(t))^2dt = \int_C(q(t)-bq(t))^2dt+\int_B a^2dt=(1-b)^2 +a^2m(B)$=
$(1-b)^2+(1-b^2)= 1-2b+b^2+1-b^2=2-2b=2(1-b)<2\epsilon/2=\epsilon$ so that $q$ is in the
closure of~$A_0$ in~$L^2[0,1]$ and this is true for every $q$ in $A\setminus A_0$.\\
Finally, if $q\not\in A$, we show $q$ is not in the closure of~$A_0$ in~$L^2[0,1]$.
If $||q||_2\not =1$, then clearly $q$ is not in the closure. Thus, assume $||q||_2=1$.
Since $q\not\in A$, then a measurable subset $B$ of~$[0,1]$ exists, $m(B)>0$, on which
$q<0$. Thus, $\int_B q(t)^2 >0$ (i of Proposition~2.34) so that $\int_B (q(t)-\hat{q}(t))^2dt$
$>\int_B q(t)^2dt$ for any $\hat{q}\in A_0$. Thus, $q$ can not be in the closure of $A_0$
in~$L^2[0,1]$ and $A$ must then be the closure of $A_0$ in~$L^2[0,1]$.\\ \smallskip\\
{\bf Corollary 4.4 (SRSF's of functions in $\Gamma_0$ and $\Gamma$; orbit of the
constant function equal to 1 \cite{lahiri}):} $Q_0=\{\overline{q}: \overline{q}=
\sqrt{\dot{\gamma}} \ \mathrm{a.e.\ on}\ [0,1],\ \mathrm{some}\ \gamma\in\Gamma_0\}$ has closure
in $L^2[0,1]$ equal to $Q=\{\overline{q}: \overline{q}= \sqrt{\dot{\gamma}}\ \mathrm{a.e.\ on}
\ [0,1],\ \mathrm{some}\ \gamma\in\Gamma\}$. In addition,
with $A_0=\{q\in L^2[0,1]: ||q||_2=1, \ q>0\ \mathrm{a.e.\ on} \ [0,1]\}$,
$A=\{q\in L^2[0,1]: ||q||_2=1, \ q\geq 0 \ \mathrm{a.e.\ on}\ [0,1]\}$,
and $q_0$ the constant function equal to 1 on~$[0,1]$, then $[q_0]_{\Gamma_0}= Q_0=A_0$ and
$cl([q_0]_{\Gamma_0})=Q=A$. \\ \smallskip\\
{\bf Proof:} If $\gamma\in\Gamma_0$ and $\overline{q}=\sqrt{\dot{\gamma}}$ a.e. on~$[0,1]$,
then $\int_0^1 \overline{q}(t)^2dt = \int_0^1 \dot{\gamma}(t)dt=\gamma(1)-\gamma(0)=1$. Also
$\overline{q}>0$ a.e. on $[0,1]$ since $\gamma\in\Gamma_0$ so that $\overline{q}\in A_0$.
Thus, $Q_0\subseteq A_0$. On the other hand, if $q\in A_0$, then $q>0$ a.e. on~$[0,1]$,
$q\in L^2[0,1]$, $||q||_2=1$. By Proposition~4.3, $\gamma$ defined for each $t\in [0,1]$
by $\gamma(t) = \int_0^t q(s)|q(s)|ds$ $=\int_0^t q(s)^2ds$ is absolutely
continuous on~$[0,1]$ with $q$ equal to the SRSF of $\gamma$ a.e. on~$[0,1]$. Clearly
$\gamma(0)=0$, $\gamma(1)= ||q||_2=1$, $\dot{\gamma} = q^2$ a.e. on~$[0,1]$, thus
$\gamma\in\Gamma_0$ and $\sqrt{\dot{\gamma}}=q$ a.e. on~$[0,1]$ so that $q\in Q_0$.
Thus~$Q_0=A_0$.  Similarly,~$Q=A$ so that the closure of $Q_0$ in $L^2[0,1]$ is~$Q$ by
Proposition~4.12. Finally, $[q_0]_{\Gamma_0}=\{\overline{q}:\overline{q} = (q_0,\gamma)=
(q_0\circ \gamma)\sqrt{\dot{\gamma}}=\sqrt{\dot{\gamma}}\ \mathrm{a.e.\ on}\ [0,1],
\ \mathrm{some}\ \gamma\in\Gamma_0 \} = Q_0 = A_0$, and since, as just proved, the closure
of $Q_0$ in $L^2[0,1]$ is $Q$, then $cl([q_0]_{\Gamma_0})= Q = A$.
\\ \smallskip\\
{\bf Corollary 4.5 (\cite{lahiri}):} With $q_0$ the constant function equal to 1
on~$[0,1]$, given $q_1$, $q_2 \in L^2[0,1]$, $||q_1||_2=||q_2||_2=1$, if either
(i) $q_1\geq 0$ a.e. on $[0,1]$ and $q_2\geq 0$ a.e. on $[0,1]$,
or (ii) $q_1\leq 0$ a.e. on $[0,1]$ and $q_2\leq 0$ a.e. on $[0,1]$, then in the case of
(i) it must be that $cl([q_1]_{\Gamma_0}) = cl([q_2]_{\Gamma_0}) = cl([q_0]_{\Gamma_0})$,
and in the case of (ii) it must be that $cl([q_1]_{\Gamma_0}) = cl([q_2]_{\Gamma_0}) =
cl([-q_0]_{\Gamma_0})$. In both cases a sequence $\{\gamma_n\}$ exists in $\Gamma_0$ with
$(q_1,\gamma_n)\rightarrow q_2$ in~$L^2[0,1]$.\\ \smallskip\\
{\bf Proof:} With $A=\{q\in L^2[0,1]: ||q||_2=1, \ q\geq 0 \ \mathrm{a.e.\ on}\ [0,1]\}$,
if (i) is true, then $q_1$, $q_2\in A = cl([q_0]_{\Gamma_0})$
(Corollary~4.4) so that by Proposition~4.10,
$cl([q_1]_{\Gamma_0}) = cl([q_2]_{\Gamma_0}) = cl([q_0]_{\Gamma_0})$.
On the other hand, if (ii) is true, then using similar arguments as above with
$-q_0=-1$ taking the place of~$q_0$, it then follows that
$cl([q_1]_{\Gamma_0}) = cl([q_2]_{\Gamma_0}) = cl([-q_0]_{\Gamma_0})$. In both cases
$cl([q_1]_{\Gamma_0}) = cl([q_2]_{\Gamma_0})$ implies the existence of~$\{\gamma_n\}$.
\\ \smallskip\\
{\bf Corollary 4.6 (\cite{lahiri}):} Given $q_1$, $q_2\in L^2[0,1]$,
and a sequence of numbers $t_0=0 <t_1<\ldots <t_n=1$, such that for each $i$, $i=1,\ldots,n$,
$\int_{t_{i-1}}^{t_i}q_1(t)^2dt = \int_{t_{i-1}}^{t_i}q_2(t)^2dt$, and either
$q_1\geq 0$ and $q_2\geq 0$ a.e. on $[t_{i-1},t_i]$, or $q_1\leq 0$ and $q_2\leq 0$ a.e. on
$[t_{i-1},t_i]$, then $cl([q_1]_{\Gamma_0})= cl([q_2]_{\Gamma_0})$.\\ \smallskip\\
{\bf Proof:} Given $i$, $1\leq i\leq n$, then a sequence
$\{\lambda_n^i\}$ exists of absolutely continuous functions,
$\lambda_n^i:[t_{i-1},t_i]\rightarrow [t_{i-1},t_i]$, $\dot{\lambda}_n^i>0$ a.e.
on~$[t_{i-1},t_i]$, $\lambda_n^i(t_{i-1})=t_{i-1}$, $\lambda_n^i(t_i)=t_i$ for each $n$ such
that $(q_1,\lambda_n^i)\rightarrow q_2$ in~$L^2[t_{i-1},t_i]$.
Here $(q_1,\lambda_n^i)$ is understood to be $(q_1\circ\lambda_n^i)\sqrt{\dot{\lambda}_n^i}$
and $L^2[t_{i-1},t_i]$ the set of square-integrable functions over~$[t_{i-1},t_i]$. Proof of
the existence of $\{\lambda_n^i\}$ along the lines of that of Corollary 4.5 with
$[t_{i-1},t_i]$ taking the place of $[0,1]$ and the value of 
$\int_{t_{i-1}}^{t_i} q_1(t)^2dt = \int_{t_{i-1}}^{t_i} q_2(t)^2dt$
not necessarily equal to~1.\\
Finally, define a sequence of functions $\{\gamma_n\}$, $\gamma_n:[0,1]\rightarrow [0,1]$,
by setting $\gamma_n(t)=\lambda_n^i(t)$ if $t\in [t_{i-1},t_i]$ for each~$n$. It follows
$\gamma_n$ is absolutely continuous, $\gamma_n(0)=0$, $\gamma_n(1)=1$, $\dot{\gamma}_n>0$ a.e.
on~$[0,1]$ for each~$n$. Thus $\{\gamma_n\}\subset \Gamma_0$ and since
$(q_1,\lambda_n^i)\rightarrow q_2$ in~$L^2[t_{i-1},t_i]$ for each~$i$, then
$(q_1,\gamma_n)\rightarrow q_2$ in $L^2[0,1]$. Thus, $q_2\in cl([q_1]_{\Gamma_0})$ and by
Proposition~4.10, then $cl([q_1]_{\Gamma_0})=cl([q_2]_{\Gamma_0})$.
\\ \smallskip\\
{\bf Corollary 4.7 (\cite{lahiri}):} Given $q_1$, $q_2\in L^2[0,1]$,
and two sequences of numbers $t_0=0 <t_1<\ldots <t_n=1$, $t'_0=0 <t'_1<\ldots <t'_n=1$,
such that for each $i$, $i=1,\ldots,n$,
$\int_{t_{i-1}}^{t_i}q_1(t)^2dt = \int_{t_{i-1}'}^{t_i'}q_2(t)^2dt$, and either
$q_1\geq 0$ a.e. on $[t_{i-1},t_i]$ and $q_2\geq 0$ a.e. on $[t_{i-1}',t_i']$,
or $q_1\leq 0$ a.e. on $[t_{i-1},t_i]$ and $q_2\leq 0$ a.e. on $[t_{i-1}',t_i']$,
then $cl([q_1]_{\Gamma_0})= cl([q_2]_{\Gamma_0})$. \\ \smallskip\\
{\bf Proof:} Let $\gamma$ be the piecewise linear element of
$\Gamma_0$ for which $\gamma(t_i')=t_i$, $i=0,\ldots,n$, and let~$w=(q_1,\gamma)$.
It then follows by Corollary~3.14 (change of variable) that for each $i$, $i=1,\ldots,n$,
$\int_{t_{i-1}'}^{t_i'} w(t)^2dt = \int_{t_{i-1}'}^{t_i'}
(q_1(\gamma(t)))^2\dot{\gamma}(t)dt = \int_{t_{i-1}}^{t_i} q_1(s)^2ds
= \int_{t_{i-1}'}^{t_i'}q_2(t)^2dt$, and since $w\geq 0$ a.e. on~$[t_{i-1}',t_i']$
if $q_1\geq 0$ a.e. on~$[t_{i-1},t_i]$, and $w\leq 0$ a.e. on~$[t_{i-1}',t_i']$
if $q_1\leq 0$ a.e. on~$[t_{i-1},t_i]$, then $w$ and $q_2$ satisfy the hypothesis of
Corollary~4.6 for the sequence $\{t_n'\}$ so that $cl([w]_{\Gamma_0})=cl([q_2]_{\Gamma_0})$.
Since $w=(q_1,\gamma)$, $\gamma\in\Gamma_0$, then $w\in [q_1]_{\Gamma_0}$ so that by
Observation~4.12, $[w]_{\Gamma_0}=[q_1]_{\Gamma_0}$ and therefore
$cl([q_1]_{\Gamma_0})= cl([q_2]_{\Gamma_0})$.\\ \smallskip\\
{\bf Proposition 4.13 (\cite{lahiri}):} Given $q\in L^2[0,1]$, then
$[q]_{\Gamma}\subseteq cl([q]_{\Gamma_0})$.
\\ \smallskip\\
{\bf Proof:} The proposition is first proved for step functions on~$[0,1]$.  Accordingly,
we assume $q$ is a step function and $\gamma\in\Gamma$.\\
Let $t_0=0<t_1<\ldots <t_n=1$ be the set of numbers that define the partition
associated with $q$ as a step function. For each $i$, $i=0,\ldots,n$, let $t_i'\in [0,1]$
be such that $\gamma(t_i')=t_i$ with $t_0'=0$ and~$t_n'=1$. Note $t_0'=0<t_1'<\ldots <t_n'=1$,
as $\gamma$ is a nondecreasing function from $[0,1]$ onto~$[0,1]$.\\
Let $w=(q,\gamma)$. It then follows by Corollary~3.14 (change of variable) that for each $i$,
$i=1,\ldots,n$, $\int_{t_{i-1}'}^{t_i'} w(t)^2dt = \int_{t_{i-1}'}^{t_i'}
(q(\gamma(t)))^2\dot{\gamma}(t)dt = \int_{t_{i-1}}^{t_i} q(s)^2ds$, and either
$q\geq 0$ a.e. on $[t_{i-1},t_i]$ and $w\geq 0$ a.e. on $[t_{i-1}',t_i']$,
or $q\leq 0$ a.e. on $[t_{i-1},t_i]$ and $w\leq 0$ a.e. on $[t_{i-1}',t_i']$.
Thus, by Corollary~4.7, $cl([q]_{\Gamma_0})= cl([w]_{\Gamma_0})$ so that, in particular,
$w=(q,\gamma)\in cl([q]_{\Gamma_0})$ and therefore, since $\gamma$ is arbitrary in $\Gamma$,
then $[q]_{\Gamma}\subseteq cl([q]_{\Gamma_0})$.\\
Finally, we assume $q$ is any function in $L^2[0,1]$ and $\gamma\in\Gamma$. Given
$\epsilon>0$, by Proposition~2.44 (density of step functions in $L^p$), there is a
step function $v$ on $[0,1]$ such that~$||q-v||_2<\epsilon/3$. As just proved above,
$(v,\gamma)\in cl([v]_{\Gamma_0})$ so that for some $\overline{\gamma}\in\Gamma_0$ it must
be that $||(v,\gamma)-(v,\overline{\gamma})||_2<\epsilon/3$. Thus, by Corollary~4.1 and
Observation~4.8 (action of $\Gamma_0$ and $\Gamma$ is distance preserving)
\begin{eqnarray*}
||(q,\gamma)-(q,\overline{\gamma})||_2 &=& ||(q,\gamma)-(v,\gamma)||_2 +
||(v,\gamma)-(v,\overline{\gamma})||_2\\
 &+& ||(v,\overline{\gamma})-(q,\overline{\gamma})||_2\\
 &=& ||q-v||_2 + ||(v,\gamma)-(v,\overline{\gamma})||_2 + ||v-q||_2\\
 &<& \epsilon/3 + \epsilon/3 + \epsilon/3 = \epsilon
\end{eqnarray*}
so that $(q,\gamma)\in cl([q]_{\Gamma_0})$ and therefore, since $\gamma$ is arbitrary in
$\Gamma$, then $[q]_{\Gamma}\subseteq cl([q]_{\Gamma_0})$.\\ \smallskip\\
{\bf Proposition 4.14 (Constant-speed parametrization of an absolutely continuous function
\cite{stein}):} Given $f\in AC[0,1]$, then there exist $h\in AC[0,1]$,
$\gamma\in\Gamma$, such that $|h'| = L = \int_0^1 |f'(t)|dt$ (the length of $f$)
a.e. on~$[0,1]$ and~$f=h\circ\gamma$ on~$[0,1]$. \\ \smallskip\\
{\bf Proof:} Given $f\in AC[0,1]$, let $L=\int_0^1 |f'(t)|dt$. If $L=0$ then $f$ is
constant on~$[0,1]$ (i of Proposition~2.34, Proposition~3.8). Otherwise,
define $\gamma:[0,1]\rightarrow [0,1]$ by $\gamma(t)=(1/L)\int_0^t |f'(s)|ds$ for each
$t\in [0,1]$. Accordingly, $\gamma(0)=0$, $\gamma(1)=(1)$, $\gamma\in AC[0,1]$ by
Proposition~3.11, and $\dot{\gamma}(t)=(1/L)|f'(t)|$ a.e. on~$[0,1]$ by Proposition~3.7
so that $\dot{\gamma}\geq 0$ a.e. on~$[0,1]$. Thus~$\gamma\in\Gamma$.\\
Given $s\in [0,1]$, then for some $t\in [0,1]$ it must be that $\gamma(t)=s$. Define
$h:[0,1]\rightarrow {\bf R}$ by~$h(s)=f(t)$. The function $h$ is well defined for if
$s=\gamma(t_1)=\gamma(t_2)$, $t_1<t_2 \in [0,1]$, then $0 = \int_0^{t_2}|f'(x)|dx -
\int_0^{t_1}|f'(x)|dx = \int_{t_1}^{t_2} |f'(x)|dx$. Thus, by i of Proposition~2.34,
$f'=0$ a.e. on~$[t_1,t_2]$ so that by Proposition~3.8 $f$~is constant on~$[t_1,t_2]$
and, in particular,~$f(t_1)=f(t_2)$.\\
Clearly $h(\gamma(t)) = f(t)$ for each $t\in [0,1]$. Note for $s_1, s_2\in [0,1]$,
$s_1<s_2$, then $s_1=\gamma(t_1)$, $s_2=\gamma(t_2)$, $t_1$, $t_2\in [0,1]$, $t_1<t_2$, and
\[ |h(s_2)-h(s_1)| = |f(t_2)-f(t_1)|=|\int_{t_1}^{t_2} f'(x)dx|\leq \int_{t_1}^{t_2}|f'(x)|dx
=L\cdot (s_2-s_1). \]
From this inequality it follows clearly that $h\in AC[0,1]$ (Definition~3.5). Accordingly,
$h$ is differentiable a.e. on~$[0,1]$ and $|h'|\leq L$ a.e. on~$[0,1]$ also from the
inequality. Note that by Corollary~3.14 (change of variable) and Corollary~3.13 (chain rule),
then
\[ \int_0^1 |h'(s)|ds = \int_0^1 |h'(\gamma(t))|\dot{\gamma}(t)dt = \int_0^1|f'(t)|dt = L.\]
By i of Proposition~2.34, then $|h'|=L$ a.e. on~$[0,1]$. \\ \smallskip\\
{\bf Corollary 4.8:} Given $q\in L^2[0,1]$, then there exist $w\in L^2[0,1]$, $\gamma\in\Gamma$,
such that $|w|=\sqrt{L}$ a.e. on~$[0,1]$ and $q=(w,\gamma)$ a.e. on~$[0,1]$, where
$L=\int_0^1|f'(t)|dt$ (the length of~$f$), $f\in AC[0,1]$, the SRSF of~$f$ equal
to~$q$ a.e. on~$[0,1]$. In particular, if $||q||_2=1$ so that $L=1$, then
$|w|=1$ a.e. on~$[0,1]$.
\\ \smallskip\\
{\bf Definition 4.8:} A function $q\in L^2[0,1]$, is said to be in {\bf standard form}
if for measurable subsets $A$, $B$ of~$[0,1]$, with $A\cap B=\emptyset$, $A\cup B=[0,1]$, then
\[ q(t) = \left\{ \begin{array}{ll}
            1 & \mathrm{for\ } t \in A\\
           -1 & \mathrm{for\ } t \in B.
                    \end{array}
            \right. \]
Clearly, if $q$ is in standard form, then $||q||_2=1$.\\
Let $SF[0,1]=\{ q: q\in L^2[0,1],\ q\ \mathrm{in\ standard\ form}\}.$
\\ \smallskip\\
{\bf Proposition 4.15 (\cite{lahiri}):} Given $q$, $w\in SF[0,1]$, if $q\not=w$ in $L^2$,
i.e., if $m(\{t:t\in [0,1], q(t)\not=w(t)\})>0$, then $w\not\in cl([q]_{\Gamma_0})$.
Thus $cl([w]_{\Gamma_0})\cap cl([q]_{\Gamma_0}) = \emptyset$.\\
Proof in \cite{lahiri} using Corollary~3.16 (Change of variable for Lebesgue integral over
a measurable set) and Observation~2.25 (Schwarz's inequality over a measurable set).
\\ \smallskip\\
{\bf Corollary 4.9: (Uniqueness of constant-speed parametrization):}
Given $\tilde{q}\in L^2[0,1]$, $||\tilde{q}||_2=1$, if for $\gamma$, $\tilde{\gamma}\in\Gamma$,
and $q$, $w\in SF[0,1]$, $\tilde{q}= (q,\gamma)$ and $\tilde{q}=(w,\tilde{\gamma})$,
then $q=w$ a.e. on~$[0,1]$. \\ \smallskip\\
{\bf Proof:} By Proposition 4.13, $[q]_{\Gamma}\subseteq cl([q]_{\Gamma_0})$
and $[w]_{\Gamma}\subseteq cl([w]_{\Gamma_0})$. Thus, $\tilde{q}\in 
cl([q]_{\Gamma_0})\cap cl([w]_{\Gamma_0})$ so that by Proposition~4.15, $q=w$ a.e. on~$[0,1]$.
\\ \smallskip\\
{\bf Proposition 4.16 (\cite{lahiri}):} Given $w\in SF[0,1]$, then
$cl([w]_{\Gamma_0})= [w]_{\Gamma}$. \\ \smallskip\\
{\bf Proof:} From Proposition 4.13, we know $[w]_{\Gamma}\subseteq cl([w]_{\Gamma_0})$.
Thus, it suffices to show $cl([w]_{\Gamma_0}) \subseteq [w]_{\Gamma}$. For this purpose,
let $\tilde{q}$ be in $cl([w]_{\Gamma_0})$. Clearly $||\tilde{q}||_2=1$, and by Corollary~4.8,
for some $q\in SF[0,1]$, and some $\gamma\in\Gamma$, it must be that $\tilde{q}=(q,\gamma)$ a.e.
on~$[0,1]$. By Proposition~4.13, $[q]_{\Gamma}\subseteq cl([q]_{\Gamma_0})$. Thus, $\tilde{q}\in 
cl([q]_{\Gamma_0})\cap cl([w]_{\Gamma_0})$ so that by Proposition~4.15, $q=w$ a.e. on~$[0,1]$,
and therefore $\tilde{q}$ is in~$[w]_{\Gamma}$.
Thus $cl([w]_{\Gamma_0}) \subseteq [w]_{\Gamma}$. \\ \smallskip\\
{\bf Corollary 4.10 (\cite{lahiri}):} Given $q\in L^2[0,1]$, if $q\not= 0$ a.e. on~$[0,1]$, then
$cl([q]_{\Gamma_0})= [q]_{\Gamma}$. \\ \smallskip\\
{\bf Proof:} If $||q||_2=1$, then by Corollary~4.8, for some $w\in SF[0,1]$, and some
$\gamma\in\Gamma$, it must be that $q=(w,\gamma)$ a.e. on~$[0,1]$. Since $q\not= 0$ a.e
on~$[0,1]$, then $(w,\gamma)= (w\circ\gamma)\sqrt{\dot{\gamma}}\not= 0$ a.e. on~$[0,1]$,
and therefore $\dot{\gamma}\not=0$ a.e. on~$[0,1]$. Thus, $\gamma\in\Gamma_0$ so that
$[q]_{\Gamma_0}=[w]_{\Gamma_0}$ and $(q,\gamma^{-1})= ((w,\gamma),\gamma^{-1})= w$
a.e. on~$[0,1]$ (Proposition~4.5).  Accordingly, $cl([q]_{\Gamma_0}) = cl([w]_{\Gamma_0})=
[w]_{\Gamma} = [(q,\gamma^{-1})]_{\Gamma}$ (Proposition~4.16). Given
$\tilde{q}\in [(q,\gamma^{-1})]_{\Gamma}$, then for some $\tilde{\gamma}\in\Gamma$,
$\tilde{q}=((q,\gamma^{-1}),\tilde{\gamma})=(q,\gamma^{-1}\circ\tilde{\gamma})$ a.e. on~$[0,1]$
(Proposition~4.5), and since $\gamma^{-1}\circ\tilde{\gamma}\in\Gamma$ (Observation~4.4),
then~$\tilde{q}\in [q]_{\Gamma}$. On the other hand, given $\tilde{q}\in [q]_{\Gamma}$, then for
some $\tilde{\gamma}\in\Gamma$, $\tilde{q}=(q,\tilde{\gamma})=
(q,\gamma^{-1}\circ\gamma\circ\tilde{\gamma})=((q,\gamma^{-1}),\gamma\circ\tilde{\gamma})$
a.e. on~$[0,1]$ (Proposition~4.5), and since $\gamma\circ\tilde{\gamma}\in\Gamma$
(Observation~4.4), then~$\tilde{q}\in [(q,\gamma^{-1})]_{\Gamma}$. Thus 
$[(q,\gamma^{-1})]_{\Gamma}=[q]_{\Gamma}$ and therefore $cl([q]_{\Gamma_0})= [q]_{\Gamma}$.\\
If $||q||_2\not=1$, then clearly $||q||_2\not=0$, $q/||q||_2\not=0$ a.e. on~$[0,1]$, and as
just proved $cl([q/||q||_2]_{\Gamma_0})=[q/||q||_2]_{\Gamma}$. Given
$\tilde{q}\in cl([q]_{\Gamma_0})$, then for a sequence $\{\gamma_n\}\subset\Gamma_0$,
$(q,\gamma_n)\rightarrow\tilde{q}$ in $L^2$. Thus
$(q/||q||_2,\gamma_n)\rightarrow\tilde{q}/||q||_2$ in $L^2$ implying
$\tilde{q}/||q||_2=(q/||q||_2,\gamma)$ for some~$\gamma\in\Gamma$, and therefore
$\tilde{q}=(q,\gamma)$ so that $\tilde{q}\in [q]_{\Gamma}$. On the other hand, given
$\tilde{q}\in [q]_{\Gamma}$, then for some $\gamma\in\Gamma$, $\tilde{q}=(q,\gamma)$. Thus
$\tilde{q}/||q||_2=(q/||q||_2,\gamma)$ implying for a sequence $\{\gamma_n\}\subset\Gamma_0$,
$(q/||q||_2,\gamma_n)\rightarrow\tilde{q}/||q||_2$ in $L^2$, and therefore 
$(q,\gamma_n)\rightarrow\tilde{q}$ in $L^2$ so that $\tilde{q}\in cl([q]_{\Gamma_0})$. Thus
$cl([q]_{\Gamma_0})= [q]_{\Gamma}$.
\\ \smallskip\\
{\bf Observation 4.17:} As noted in \cite{lahiri}, given $f_1$, $f_2\in AC[0,1]$, and their
SRSF's $q_1$, $q_2\in L^2[0,1]$, respectively, if $\tilde{q}_1\in cl([q_1]_{\Gamma_0})$,
$\tilde{q}_2\in cl([q_2]_{\Gamma_0})$ exist such that
$d(cl([q_1]_{\Gamma_0}),cl([q_2]_{\Gamma_0})) = ||\tilde{q}_1-\tilde{q}_2||_2$,
assuming without any loss of generality that $q_1\not=0$ a.e. on~$[0,1]$, $q_2\not=0$ a.e.
on~$[0,1]$ (Corollary~4.8), then by Corollary~4.10 above there exist $\gamma_1$,
$\gamma_2 \in\Gamma$, such that $\tilde{q}_1=(q_1,\gamma_1)$ and
$\tilde{q}_2=(q_2,\gamma_2)$. The pair $\gamma_1$, $\gamma_2$ is called an {\bf optimal
matching} for~$f_1$,~$f_2$. In particular, it is proved in~\cite{lahiri} that if
%$q_1$, $q_2$ are such that
at least one of $cl([q_1]_{\Gamma_0})$, $cl([q_2]_{\Gamma_0})$ contains the
SRSF of a piecewise linear function, then $\tilde{q}_1$, $\tilde{q}_2$ exist as above 
and therefore there is an optimal matching for~$f_1$,~$f_2$. This is actually proved
in~\cite{lahiri} for absolutely continuous functions $f_1$, $f_2$ with range~${\bf R}^n$.
\\ \smallskip\\
%\\ \pagebreak\\
{\bf\large Summary}\\ \smallskip\\
In order to understand the theory of functional data and shape analysis as presented in
Srivastava and Klassen's textbook ``Functional and Shape Data Analysis" \cite{srivastava},
it is important to understand the basics of Lebesgue integration and absolute continuity,
and the connections between them. In this paper of the survey type, we have tried to provide
a way to do exactly that. We have reviewed fundamental concepts and results about Lebesgue
integration and absolute continuity, some results connecting the two notions, most of the
material borrowed from Royden's ``Real Analysis" \cite{royden} and Rudin's ``Principles of
Mathematical Analysis" \cite{rudin}. Additional important material was obtained from
Saks' \cite{saks}, and Serrin and Varberg's \cite{serrin} seminal papers. In addition,
we have presented fundamental concepts and results about functional data and shape
analysis in 1-dimensional space, in the process shedding light on its dependence on
Lebesgue integration and absolute continuity, and the connections between them, most of
the material borrowed from Srivastava and Klassen's aforementioned textbook.
Additional material presented at the end of the paper was obtained from Lahiri, Robinson
and Klassen's outstanding manuscript~\cite{lahiri}.
\\ \smallskip\\
{\bf\large Acknowledgements} \\ \smallskip\\
I am most grateful to Professor James F. Lawrence of George Mason University and the National
Institute of Standards and Technology for the many insightful conversations on the subjects of
Lebesgue integration and absolute continuity, and to Professor Eric Klassen of Florida State
University for his generosity in always providing answers to my questions about his remarkable
work on shape analysis.

\begin{center}
{\bf \Large Index of Terms}
\end{center}
\noindent
Each entry crossreferenced to one of the following: Dx.y (Definition x.y),\\
Px.y (Proposition x.y), Cx.y (Corollary x.y), Ox.y (Observation x.y).\\
Given x.y, y is then the number of the Definition, Proposition, Corollary or
Observation in Section~x in which the entry can be found. For example, the entry
``Borel" is associated in the index with D2.3 which means it can be found
in Definition~3 of Section~2.
\\ \smallskip \\
Absolute continuity of the indefinite integral, P3.11\\
Absolute continuity of the Lebesgue integral, P2.38\\
absolutely continuous, D3.5\\
Absolutely continuous $f$ implies $f$ is of bounded variation, P3.10\\
Absolutely continuous $f$ is constant if $f'$ is zero a.e., P3.8\\
Absolutely continuous $f$ maps measurable sets to measurable sets, P3.21\\
Absolutely continuous $f$ maps zero-measure sets to zero-measure sets, P3.21\\
$AC[0,1]$, O4.1\\
$AC^0[0,1]$, D4.4\\
action of $\Gamma$ on $L^2[0,1]$, D4.3\\
Action of $\Gamma$ on $L^2[0,1]$ is by semi-isometries, P4.6\\
Action of $\Gamma$ on $L^2[0,1]$ is distance and norm preserving, O4.8\\
action of $\Gamma_0$ on $AC^0[0,1]$, D4.5\\
Action of $\Gamma_0$ on $AC^0[0,1]$ with Fisher-Rao metric is by isometries, P4.8\\
action of $\Gamma_0$ on $L^2[0,1]$, D4.3\\
Action of $\Gamma_0$ on $L^2[0,1]$ is by isometries, P4.6\\
Action of $\Gamma_0$ on $L^2[0,1]$ is distance preserving, C4.1\\
Action of $\Gamma_0$ on $L^2[0,1]$ is norm preserving, C4.2\\
{\em adapt-DP}, O4.15\\
admissible class $\Gamma$ of warping functions, D4.2\\
a.e., D2.15\\
algebra, D2.1\\
almost everywhere, D2.15\\
Approximation of a measurable function by simple functions, P2.23\\
Banach-Zarecki Theorem, P3.22\\
Borel set, D2.3\\
Bounded Convergence Theorem, C2.6\\
bounded variation, D3.3\\
canonical representation, O2.11\\
Cantor function, O3.3\\
Cantor set, O2.4\\
Carath\'{e}odory's criterion, D2.9\\
Cauchy sequence, D2.23\\
Chain rule, P3.26\\
Chain rule (Alternate form), C3.13\\
Change of variable for Lebesgue integral, P3.27\\
Change of variable for Lebesgue integral (Alternate form I), C3.14\\
Change of variable for Lebesgue integral (Alternate form II), C3.15\\
Change of variable for Lebesgue integral over a measurable set, C3.16\\
Change of variable for Riemann integral, P3.14\\
characteristic function, D2.13\\
closed set, D2.3\\
$cl([q]_{\Gamma_0})$, D4.6\\
compact set, D2.4\\
complete normed linear space, D2.23\\
Composition of absolutely continuous functions, P3.25\\
Constant-speed parametrization of an absolutely continuous function, P4.14\\
continuous, O2.10\\
converge in norm, D2.23\\
Countable additivity of $m$, P2.12\\
Countable additivity of the Lebesgue integral, P2.36\\
Countable subadditivity of $m$, P2.12\\
Countable subadditivity of $m^*$, P2.7\\
Density of simple and step functions in $L^p$ space, P2.44\\
$d(cl([q_1]_{\Gamma_0}),cl([q_2]_{\Gamma_0}))$, C4.3\\
$d([q_1]_{\Gamma_0},[q_2]_{\Gamma_0})$, O4.12\\
derivative, D3.1\\
diffeomorphism, O4.6\\
Differentiability of the indefinite integral, P3.7\\
Differentiability of the Riemann integral, C3.1\\
differentiable, D3.1\\
differentiable (generalized to smooth manifolds), O4.6\\
differentiable manifolds, O4.6\\
differential, O4.6\\
Distance between equivalence classes in $L^2[0,1]/\sim\,$, C4.3\\
Distance between functions in $AC[0,1]$, O4.11\\
Egoroff's Theorem, P2.22\\
Equivalent conditions for a measurable function, P2.17\\
Equivalent conditions for a measurable set, P2.15\\
Equivalent conditions for an absolutely continuous function, P3.12\\
extended real numbers, D2.7\\
extended to all of $[c,d]$, O3.8\\
$f$ extended to all of $[c,d]$, O3.8\\
Fatou's Lemma, P2.39\\
Fisher-Rao distance, O4.11\\
Fisher-Rao metric, D4.4\\
Fisher-Rao metric on $AC^0[0,1]$ under SRSF representation becomes $L^2[0,1]$ metric, P4.9\\
Fundamental Theorem of calculus I, P3.1\\
Fundamental Theorem of calculus II, P3.2\\
Fundamental Theorem of calculus for continuous functions, C3.1\\
Fundamental Theorem of Lebesgue integral calculus, C3.6\\
Fundamental Theorem of Lebesgue integral calculus (Alternate form), P3.13\\
$\Gamma$, D4.2\\
$\Gamma_0$, D4.2\\
geodesic, O4.6\\
geodesic distance, O4.6\\
group $\Gamma_0$ of invertible warping functions, D4.2\\
Heine-Borel, P2.3\\
H\"{o}lder's inequality, P2.41\\
indefinite integral, D3.4\\
Indefinite integral of $f$ zero everywhere, then $f$ is zero a.e., P3.6\\
Integrable equivalent to measurable, P2.32\\
Inverse function theorem, P3.23\\
invertible warping functions, D4.2\\
isometry, O4.6\\
Jordan decomposition, P3.5\\
$L^2[0,1]/\sim\,$, O4.13\\
$L^2[0,1]/\Gamma_0$, O4.12\\
$L^2[0,1]$'s equivalence with the set of all SRSF's, P4.3\\
$L^p[0,1]$ or $L^p$ space, D2.22\\
$L^p$ norm, D2.22\\
$L^p$ norm of a function, D2.22\\
$L^{\infty}[0,1]$ or $L^{\infty}$ space, D2.22\\
$L^{\infty}$ norm, D2.22\\
$L^{\infty}$ norm of a function, D2.22\\
Lebesgue integrable, D2.21\\
Lebesgue integral of a measurable function, D2.21\\
Lebesgue integral of a measurable nonnegative function, D2.20\\
Lebesgue integral of a simple function, D2.19\\
(Lebesgue) measurable function, D2.11\\
(Lebesgue) measurable set, D2.9\\
Lebesgue measure $m$, D2.10\\
Lebesgue's criterion for Riemann integrability, P2.35\\
Lebesgue's Dominated Convergence Theorem, P2.40\\
Lebesgue's Monotone Convergence Theorem, P2.37\\
left-hand limit, D3.2\\
limit point, D2.5\\
Lindel\"{o}f, P2.2\\
lower Riemann integral, D2.17\\
Lusin's Theorem, P2.24\\
$m$, D2.10\\
$m^*$, D2.8\\
Measurability of the derivative of a measurable function, P3.19\\
measurable function, D2.11\\
measurable set, D2.9\\
mesh, D2.18\\
Minkowski's inequality, P2.42\\
Monotonic functions: continuity, P3.3\\
Monotonic functions: differentiability, P3.4\\
Monotonic surjective $f$ implies $f$ is continuous, C3.2\\
negative part $f^-$ of $f$, D2.16\\
Nested sequences of measurable sets Lemma, P2.13\\
normed linear space with norm~$||\cdot||$, D2.23\\
open set, D2.3\\
optimal matching, O4.17\\
orbit of the constant function equal to 1, C4.4\\
orbit $[q]_{\Gamma}$ of $q$ under $\Gamma$, D4.7\\
orbit $[q]_{\Gamma_0}$ of $q$ under $\Gamma_0$, D4.6\\
outer measure $m^*$, D2.8\\
partition, D2.12\\
perfect set, D2.6\\
positive part $f^+$ of $f$, D2.16\\
$[q]_{\Gamma}$, D4.7\\
$[q]_{\Gamma_0}$, D4.6\\
Reconstruction of an absolutely continuous function from its SRSF, P4.2\\
Riemannian manifolds, O4.6\\
Riemannian metric, O4.6\\
Riemannian structure, O4.6\\
Riemann integrable, D2.17\\
Riemann integrable implies Lebesgue integrable, P2.33\\
Riemann integral, D2.17\\
Riemann sum, D2.18\\
Riemann sums of $f$ converge if and only $f$ is Riemann integrable, P2.27\\
Riemann sums of $f$ that converge implies $f$ is bounded, P2.26\\
Riesz-Fischer, P2.43\\
right-hand limit, D3.2\\
Saks' inequality, P3.17\\
Saks' Theorem, C3.10\\
Schwarz's inequality, O2.25\\
semigroup $\Gamma$ of warping functions, D4.2\\
semi-isometry, O4.6\\
Serrin-Varberg's Theorem, P3.18\\
Serrin-Varberg's Theorem (Alternate form), C3.12\\
$SF[0,1]$, D4.8\\
$\sigma$-algebra, D2.2\\
simple function, D2.14\\
smooth manifolds, O4.6\\
Square integrability of SRSF, P4.1\\
square-root slope function (SRSF), D4.1\\
square-root velocity function (SRVF), O4.2\\
SRSF of a warped absolutely continuous function, P4.4\\
SRSF representation of functions, D4.1\\
SRSF's of functions in $\Gamma_0$ and $\Gamma$, C4.4\\
SRVF, O4.2\\
standard form, D4.8\\
step function, D2.12\\
subdivision, D2.12\\
Substitution rule for Lebesgue integral, P3.16\\
Substitution rule for Riemann integral, P3.15\\
$T_pM$, 04.6\\
tangent space $T_pM$, O4.6\\
uniformly continuous, O3.3\\
Uniqueness of constant-speed parametrization, C4.9\\
upper Riemann integral, D2.17\\
Zarecki's criterion for an absolutely continuous inverse, P3.24
%
%\\ \pagebreak\\
%{\bf \large Appendix: Change Log}\\
%{\bf Revision 2 Release 2 - June 27, 2019}\,\medskip\\
%{July 18, 2018 updates:}\\
%$\bullet$ Added one sentence for case $L=0$ in proof of Proposition 4.14.\\
%$\bullet$ Added Observation 4.17.\\
%{June 27, 2019 updates:}\\
%$\bullet$ Added DOI that points to data in Observation 4.15.

\begin{thebibliography}{5}

\bibitem{apostol}
Apostol, T. M.:
Mathematical Analysis, 2nd edition.
Reading, Massachusetts: Addison-Wesley Publishing Company. (1974)

\bibitem{bernal}
Bernal, J., Dogan, G., Hagwood, C. R.:
Fast Dynamic Programming for Elastic Registration of Curves.
Proceedings of DIFF-CVML workshop, CVPR 2016, Las Vegas, Nevada.
(2016)

\bibitem{bruckner}
Bruckner, A. M., Bruckner, J. B., Thomson, B. S.:
Real Analysis, 1st edition.
Upper Saddle River, New Jersey: Prentice-Hall. (1997)

\bibitem{burk}
Burk, F.:
Lebesgue Measure and Integration. An Introduction, 1st edition.
New York: John Wiley \& Sons. (1998)

\bibitem{cabada}
Cabada, A., Pouso, R. L.:
On First Order Discontinuous Scalar Differential Equations.
Nonlinear Studies 2 (1999) 161--170.

\bibitem{carmo1}
do Carmo, M. P.:
Differential Geometry of Curves and Surfaces.
Upper Saddle River, New Jersey: Prentice-Hall. (1976)

\bibitem{carmo2}
do Carmo, M. P.:
Riemannian Geometry.
Boston: Birkh\"{a}user. (1992)

\bibitem{halmos}
Halmos, P. R.:
Measure Theory.
New York: Springer-Verlag. (1974)

\bibitem{klassen}
Klassen, E.:
Private communication. (2018)

\bibitem{lahiri}
Lahiri, S., Robinson, D., Klassen, E.:
Precise Matching of PL curves in ${\bf R}^n$ in the
Square Root Velocity Framework.
Geometry, Imaging and Computing 2(3) (2015) 133--186.

\bibitem{lee1}
Lee, J. M.:
Introduction to Smooth Manifolds.
Graduate Texts in Mathematics, Vol. 218.
New York: Springer-Verlag. (2003)

\bibitem{lee2}
Lee, J. M.:
Riemannian Manifolds. An Introduction to Curvature.
Graduate Texts in Mathematics, Vol. 176.
New York: Springer-Verlag. (1997)

\bibitem{pugh}
Pugh, C. C.:
Real Mathematical Analysis, 2nd edition.
New York: Springer. (2015)

\bibitem{rana}
Rana, I. K.:
An Introduction to Measure and Integration, 2nd edition.
Graduate Studies in Mathematics, Vol. 45.
Providence, Rhode Island: American Mathematical Society. (2002)

\bibitem{richardson}
Richardson, L. F.:
Advanced Calculus. An introduction to Linear Analysis, 1st edition.
Hoboken, New Jersey: John Wiley \& Sons. (2008)

\bibitem{royden}
Royden, H. L.:
Real Analysis, 2nd edition.
New York: Macmillan. (1968)

\bibitem{royden2}
Royden, H. L., Fitzpatrick P. M.:
Real Analysis, 4th edition.
Upper Saddle River, New Jersey: Pearson Education. (2010)

\bibitem{rudin}
Rudin, W.:
Principles of Mathematical Analysis, 2nd edition.
New York: McGraw-Hill. (1964)

\bibitem{rudin2}
Rudin, W.:
Real and Complex Analysis, 3rd edition.
New York: McGraw-Hill. (1987)

\bibitem{saks}
Saks, S.:
Theory of the Integral. Monografie Matematyczne Tom VII, 2nd revised edition.
New York: Hafner Publishing Co. (1937)

\bibitem{schmeding}
Schmeding, A.:
Manifolds of Absolutely Continuous Curves and the Square Root Velocity
Framework.
arXiv preprint arXiv:1612.02604 (2016)

\bibitem{serrin}
Serrin, J., Varberg, D. E.:
A General Chain Rule for Derivatives and the Change of Variables
Formula for the Lebesgue Integral.
The American Mathematical Monthly 76(5) (1969) 514--520.

\bibitem{srivastava}
Srivastava, A., Klassen, E. P.:
Functional and Shape Data Analysis.
New York: Springer. (2016)

\bibitem{stein}
Stein, E. M., Shakarchi, R.:
Real Analysis: Measure Theory, Integration, and Hilbert Spaces.
Princeton, New Jersey: Princeton University Press. (2005)

\bibitem{wade}
Wade, W. R.:
An Introduction to Analysis, 3rd edition.
Upper Saddle River, New Jersey: Pearson Prentice Hall. (2004)

\bibitem{yeh}
Yeh, J.
Real Analysis. Theory of Measure and Integration, 3rd edition.
Hackensack, New Jersey: World Scientific Publishing Co. (2014)

\end{thebibliography}
\end{document}